\documentclass[12pt, reqno]{amsart}
\usepackage{mathrsfs}
\usepackage[active]{srcltx}
\usepackage{amsmath,amsthm}    
\usepackage{amssymb}    
\usepackage{mathtools}
\usepackage{longtable}
\usepackage{todonotes}
\usepackage[]{mdframed}
\usepackage[autostyle]{csquotes}
\usepackage{bm}
\usepackage{wrapfig}
\allowdisplaybreaks
\usepackage{chngcntr}

\usetikzlibrary{shapes.geometric, arrows,positioning}
\tikzstyle{mybox} = [rectangle, rounded corners, minimum width=5.5cm, text width=5cm, minimum height=1cm,align=left, draw=black, fill=black!10]
\usepackage{geometry}

\DeclareMathOperator{\dd}{d}

\UseRawInputEncoding

\usepackage{xcolor}
\definecolor{bluecite}{HTML}{0875b7}
\usepackage[unicode=true,
bookmarksopen={true},
pdffitwindow=true,
colorlinks=true,
linkcolor=bluecite,
citecolor=bluecite,
urlcolor=bluecite,
hyperfootnotes=false,
pdfstartview={FitH},
pdfpagemode= UseNone]{hyperref}

\newcommand{\ds}{\displaystyle}

\newtheorem{proposition}{Proposition}[section]
\newtheorem{theorem}{Theorem}[section]

\newtheorem{remark}{Remark}[section]
\numberwithin{equation}{section}

\usepackage{bbm}  
\usepackage{dsfont}

\author{Alexandru Krist\'aly}	

\address{Department of Economics, Babe\c s-Bolyai University, Str. Teodor Mihali 58-60, 400591 Cluj-Napoca,
	Romania \& Institute of Applied Mathematics, \'Obuda
	University, B\'ecsi \'ut 96/B, 1034
	Budapest, Hungary}

\email{alexandru.kristaly@ubbcluj.ro; kristaly.alexandru@uni-obuda.hu}
%{26D15; 35B35; 47J20; 34K20}

\title[Finsler--Sobolev inequalities: optimality \& rigidity]
{Finsler--Sobolev inequalities: optimality \& rigidity}
\date{} 

\begin{document}
	\maketitle
	
%\vspace{-0.4cm}
	\begin{abstract} We
		establish   optimal  Sobolev and logarithmic-Sobolev inequalities on non-compact,  possibly non-rever\-sible  Finsler manifolds  with non-negative weighted Ricci curvature. 
		 The optimal constants exhibit no dependence on reversibility   -- confirming Ohta's conjectural principle -- being  governed solely by the forward and backward asymptotic volume ratios.   
%				The optimal constants exhibit no dependence on reversibility   -- confirming Ohta's conjectural principle -- being  governed solely by the forward and backward asymptotic volume ratios. 
%				Two approaches are used in the proofs. First, we show the power of the optimal mass transport theory to prove optimal Sobolev inequalities on Finsler manifolds   by using  displacement convexity of certain entropy functionals. Second, we apply rearrangement techniques to  describe the rigidity in the equality cases of the above optimal inequalities, expressed through fine geometric properties of the flag and $S$-curvatures as well as the precise profiles of the corresponding extremizers. 
	We then apply the logarithmic-Sobolev inequality to establish the optimal hypercontractivity estimate for the heat flow in the same geometric setting.  
	%,  characterizing   the  equality case as well.
	The equality cases of the above optimal inequalities imply strong rigidity in the underlying manifold, characterized  by fine   properties of the flag and $S$-curvatures as well as the precise profiles of the corresponding extremizers.
		We illustrate our results by examples  including Randers metrics obtained by Zermelo navigation with Killing vector fields, and weighted convex cones in $\mathbb{R}^n$ endowed with arbitrary Minkowski norms, the latter  showing  that the forward and backward asymptotic volume ratios may differ.

%			The aim of this paper is threefold. First, by using displacement convexity of suitable entropy functionals, we prove sharp isoperimetric, Sobolev and logarithmic Sobolev inequalities on (not necessarily reversible) Finsler manifolds with non-negative $N$-Ricci curvature. This argument shows the applicability of the optimal mass transport theory, developed by S. Ohta, to derive sharp Sobolev-type inequalities on non-compact Finsler manifolds. 
%		Second, by means of rearrangement techniques, we characterize the equality cases in the above inequalities in terms of fine properties of flag and $S$-curvatures as well as precise shapes of extremizers. Third, we provide optimal hypercontractivity estimates for the heat flow on Finsler manifolds with non-negative $N$-Ricci curvature, characterizing also the equality cases. 

	\end{abstract}
\vspace{-0.4cm}
\
	
	\section{Introduction}

	\subsection{Motivation.}
		
%		The aim of this paper is to establish qualitative properties of \textit{Sobolev inequalities} on not necessarily reversible \textit{Finsler manifolds}, with particular emphasis on their \textit{optimality}.

		Sobolev inequalities -- controlling the Dirichlet energy of a function by its Lebesgue norm -- play a fundamental role in the analysis of partial differential equations,  their validity and optimality being deeply influenced by the geometry of the underlying space. A systematic investigation of these questions on Riemannian manifolds was initiated in the early 1970s with the seminal work of Aubin \cite{Aubin}; for more recent developments, see Hebey \cite{Hebey}.

	As a natural extension of Riemannian geometry, Finsler geometry provides a realistic framework to model  non-reversible  phenomena, such as the Zermelo navigation problem, see Bao--Robles--Shen \cite{BRS};
	 billiard problems on the Funk space, see Faifman \cite{Faifman};   Matsumoto mountain slope metric describing the walking law under gravity, see Matsumoto  \cite{Matsumoto};  and Wulff shapes in anisotropic crystal structures,  see Wulff \cite{Wulff}. 
	 
A Finsler manifold is  modeled over a smooth manifold  $M$ (possibly   with Lipschitz  boundary $\partial M$) endowed with a Minkowski norm $F(x,\cdot):T_xM\to [0,\infty)$ for every $x\in M$, which are required to  smoothly depend on the base point $x$.   Riemannian manifolds appear when the Minkowski norm is Hilbertian.   
	 
Non-reversible Finsler structures manifest even when the underlying space is the standard Euclidean space $\mathbb R^n$.  Indeed, if we consider an open, bounded, convex set $K\subset \mathbb R^n$  containing the origin, we may consider the weight function on directions through the Euclidean
	scalar product $\langle \cdot , \cdot  \rangle$ as 
	$$\|x\|_K=\sup\{\langle x, v \rangle: v\in K \},$$
which is reversible if and only if $K $ is centrally symmetric.    These weight functions play a crucial role in describing anisotropic surface energy equilibria in terms of  Wulff shapes. In particular, Gromov \cite{Gromov} established the optimal anisotropic $L^1$-Sobolev inequality 
\begin{equation}\label{Gromov-inequality}
\int_{\mathbb R^n}\|-\nabla u(x)\|_K {\rm d}x \geq n |K|^\frac{1}{n}\|u\|_{L^{\frac{n}{n-1}}(\mathbb R^n,\mathcal L^n)}, \quad \forall u\in C_0^1(\mathbb R^n),  
\end{equation}
where $|K|$ is the measure of $K\subset \mathbb R^n$ with respect to the Lebesgue measure $\mathcal L^n.$ The equality case, alternative proofs and quantitative versions for  \eqref{Gromov-inequality} can be found in the papers by Cabr\'e--Ros-Oton--Serra \cite{Cabre-JEMS}, 
 Figalli--Maggi--Pratelli \cite{FMP-Inventiones}   and  van Schaftingen \cite{vanSch}. The objects $(\mathbb R^n,\|\cdot\|_K,\mathcal L^n),$  with $K\subset \mathbb R^n$ having smooth boundary,  constitute the simplest Finsler structures, called as Minkowski spaces, see   Bao--Chern--Shen \cite{BCS}; indeed, they possess  zero flag and $S$-curvatures, while their geodesics are straight lines.
	  
 Driven by the study of nonlinear PDEs and entropy estimates, Sobolev inequalities have become essential not only on the spaces    $(\mathbb R^n,\|\cdot\|_K,\mathcal L^n),$ but also on more general Finsler manifolds under certain \textit{curvature} bounds.
Such problems include heat and gradient flows, see Ohta--Sturm \cite{Ohta-Sturm} and Ohta--Zhao \cite{Ohta-Zhao};   optimal mass transport, see Ohta \cite{Ohta-2018} and Krist\'aly--Zhao \cite{Kristaly-Zhao-JMPA}; and eigenvalue problems, see Ge--Shen \cite{Ge-Shen} and Wang--Xia \cite{Wang-Xia}.
 %		Many results from Riemannian geometry have their natural Finslerian counterparts, as Hopf--Rinow, Cartan--Hadamard, Gauss--Bonnet and Bonnet--Myers theorems as well as Rauch and Bishop--Gromov  comparison principles;  see e.g.\ Bao--Chern--Shen \cite{BCS}. However, %over the last two decades 

Despite the above   achievements, it would be na\"ive to think that Rie\-man\-nian  results can simply be extended to the Finsler setting. Indeed, for example, no classification of Finsler manifolds with constant flag curvature is available, see  Shen \cite{Shen-S=0,Shen-TAMS}; Katok's celebrated example  provides non-reversible Finsler metrics of constant positive flag curvature on spheres with only finitely many geometrically distinct closed geodesics, see Katok \cite{Katok} and  Bryant--Foulon--Ivanov--Matveev--Ziller \cite{Bryantetal};
 the failure of Busemann NPC inequality on general non-positively curved Finsler manifolds, see Ivanov--Lytchak \cite{Ivanov-Lyt}; 
the failure of vector space structure of Sobolev spaces  on certain non-reversible Finsler manifolds, see Krist\'aly--Rudas \cite{KR-Nonlinear}, etc.

Motivated by the above phenomena, our aim is to   study  \textit{optimal $L^p$-Sobolev-type inequalities} $(p\geq 1)$ on not necessarily rever\-sible  Finsler manifolds, possibly with boundary, together with some applications. In fact, our aim is twofold. 
	
%	First, we prove that  on the Funk space  -- which is a non-reversible Finsler manifold with negative constant flag curvature --   no   Sobolev inequality can  hold (see Theorem \ref{Thm-Sobolev-Funk}), i.e., the Rayleigh quotient vanishes.   This result refutes the natural expectation that Riemannian--Sobolev inequalities admit a direct extension to the Finsler setting. 
	
	First, on not necessarily reversible Finsler manifolds with non-negative weighted $N$-Ricci curvature, we prove the optimal $L^p$-Sobolev and  $L^p$-logarithmic-Sobolev inequalities (see Theorems \ref{Thm-Sobolev-1} \& \ref{log-Sobolev-main}). Although one could naturally expect the optimal constants to depend explicitly on the \textit{reversibility} constant, our results reveal that non-reversibility is   encoded into the forward and backward asymptotic volume ratios. 
	We also investigate the existence of extremizers for these inequalities, which yields strong rigidity of the underlying Finsler manifold: along specific geodesic rays, the flag curvature vanishes,  the $S$-curvature either vanishes or exhibits a hyperbolic behavior, while the extremizers admit an explicit profile (see Theorems \ref{Thm-Sobolev-2-equality} \& \ref{log-Sobolev-main}).

%		 Having the   optimality in Sobolev inequalities, we also focus on the existence of extremizers, which turns out to imply strong rigidity on the underlying Finsler manifold, reflected in  specific behavior of the $S$- and flag curvatures. 
		 		
		 	Second, by applying the   logarithmic-Sobolev inequality, we   establish the optimal hypercontractivity estimate  for the heat flow on Finsler manifolds with non-negative weighted  $N$-Ricci curvature; moreover, the equality case  is also characterized  (see Theorem  \ref{hypercontractivity-main-theorem}). 
		 	
%		We provide various examples to illustrate our results, including Zermelo navigation problem (belonging to the class of Randers spaces), proper convex cones in $\mathbb R^n$ and smooth Wulff sets (which are Minkowski spaces). 
	% \newpage
	
\subsection{Main results}

%	More precisely, we consider Sobolev and logarithmic Sobolev inequalities in such geometric settings, providing also applications in the theory of 
%	Modelling various nonlinear phenomena often lead us to functional inequalities. 
	
	Let $\mathcal M=(M,F,{\sf m})$ be an $n$-dimensional Finsler space $(M,F)$, possibly with boundary,   with a smooth measure ${\sf m}$.  If $p,q\geq  1$, we consider  the space  
	 $$W^{p}_q(\mathcal M)=\{u\in L^q(M,{\sf m}):F^*(Du)\in L^p(M,{\sf m})\},$$
	 where $F^*$ is the dual of $F$ and $Du$ is the differential of $u:M\to \mathbb R$,  see \S \ref{section-preliminaries}. The   \textit{optimal Sobolev inequality on $\mathcal M$} is equivalent to determine  the \textit{Rayleigh quotient on}  $\mathcal M$, i.e., 
	 $$\mathcal R^{p}_q(\mathcal M):=\inf_{u\in W^{p}_q(\mathcal M)\setminus \{0\}} \frac{\|F^*(Du)\|_{L^p(M,{\sf m})}}{\|u\|_{L^q(M,{\sf m})}}, $$
		where $\|\cdot\|_{L^s(M,{\sf m})}$ is the usual $L^s$-Lebesgue norm on $M$ with respect to the measure ${\sf m}.$ 
		
		When $\mathcal M$ is compact, constant functions on $M$ belong to $W^{p}_q(\mathcal M)$, thus   one has that $\mathcal R^{p}_q(\mathcal M)=0$ for every $p,q\geq 1$; accordingly, from now on, $\mathcal M$ will be always non-compact.    

In the simplest Euclidean case $\mathcal M_{\rm Euc}=(\mathbb R^n,|\cdot|,\mathcal L^n)$ and $1<p<n$, by scaling reasoning, the only meaningful value for $q$ is the Sobolev critical exponent $p^*_n=\frac{pn}{n-p}$. Note   that
\begin{equation}\label{aubin-talenti-constant}
	\mathcal R^{p}_{p^*_n}(\mathcal M_{\rm Euc})={\sf AT}(p,n):=\pi^{\frac{1}{2}}n^{\frac{1}{p}} \left(\frac{n-p}{p-1}\right)^{1/p'}\left(\frac{\Gamma(n/p)\Gamma(1+n/p')}{\Gamma(1+n/2)\Gamma(n)}
	\right)^{{1}/{n}},
\end{equation}
and $	\mathcal R_{p^*_n}^p(\mathcal M_{\rm Euc})$ is achieved -- up to a multiplicative factor and translation --  by the unique class of Talentian bubbles 
$u_\lambda(x)=(\lambda+|x|^{p'})^\frac{p-n}{p},$ $\lambda>0;$ see Aubin \cite{Aubin} and Talenti \cite{Talenti}.  Here, $p'=\frac{p}{p-1}$ stands for the conjugate of $p>1$. In the limit  $p\to 1$, as ${\sf AT}(p,n)\to n\omega_n^\frac{1}{n}$ (where $\omega_n$ is the volume of the open unit ball $\mathbb B^n$ in $\mathbb R^n$),  relation \eqref{aubin-talenti-constant} reduces to the functional form of the sharp isoperimetric inequality in $\mathbb R^n$, see \eqref{Gromov-inequality} with $K$ being the unit ball in $\mathbb R^n$
with center at the origin.  

If $\mathcal M_{\rm CH}=(M,g,{{\rm vol}_g})$ is an $n$-dimensional Cartan--Hadamard manifold\footnote{Complete, simply connected Riemannian manifold with non-positive sectional curvature, endowed with its canonical volume measure ${{\rm vol}_g}$.}  supporting the Cartan--Hadamard conjecture (i.e., the sharp isoperimetric inequality on $\mathcal M_{\rm CH}$),  one has 
\begin{equation}\label{CH-Sobolev}
	\mathcal R^p_{p^*_n}(\mathcal M_{\rm CH})={\sf AT}(p,n),
\end{equation}
for every $1< p<n$, where ${\sf AT}(p,n)$ is from \eqref{aubin-talenti-constant}; see   Hebey \cite[Chapter 8]{Hebey}. Note that the validity of the Cartan--Hadamard conjecture is known for low dimensions $n\in \{2,3,4\}$\footnote{Recent progress has been made by  Chen--Ghomi--Wang \cite{Chen-Ghomi-Wang1, Chen-Ghomi-Wang2} in dimensions up to $9,$  while an OpenAI paper \cite{OpenAI} claims that the conjecture holds in any dimension.}. 

Having  \eqref{CH-Sobolev}, one could expect the validity of a similar result on non-positively curved Finsler manifolds as well, endowed with the natural Busemann--Hausdorff measure ${\sf m}_{BH}$. 
However,    the  $n$-dimensional  Funk space   $\mathcal M_{\rm Funk}=(\mathbb B^n,F,{\sf m}_{BH})$ on the  open unit ball $\mathbb B^n$ in  $\mathbb R^n$  endowed with the Funk metric $F$, see Shen \cite{Shen-konyv}, has constant flag curvature $-1/4$  and for every $p,q\geq  1$,   
$$	  \mathcal R^{p}_q(\mathcal M_{\rm Funk})=0.
 $$
This relation follows from the fact that  constant functions belong to $W^{p}_q(\mathcal M_{\rm Funk})$ for every $p,q\geq  1$, because   on $\mathcal M_{\rm Funk}$, 
 ${\sf m}_{BH}$ is precisely the $n$-dimensional Lebesgue measure, thus ${\sf m}_{BH}(\mathbb B^n)=\omega_n$. Note that $\mathcal M_{\rm Funk}$ is a Finsler space 
 %and constant negative weighted $N$-Ricci curvature for every $N\geq n$, see Ohta \cite{Ohta-book}, 
 with infinite reversibility. 
  This example shows that  reasonable/optimal Sobolev inequalities on not necessarily reversible Finsler manifolds  require  not only a control of the flag curvature, but also a suitable behavior of the measure.

%  With respect to \eqref{CH-Sobolev},  we have:
%
%
%%
%%	Due to various applications, the most important functional inequalities %in PDEs -- with a high degree of certanity -- 
%%	are the Sobolev and logarithmic-Sobolev inequalities. 
%%	$\mathcal M=(M,F,{\sf m})$
%%	
%%	
%%
%%	
%%	
%%
%%	
%%	$$ \forall u\in \dot W^{1,p}(M,F,{\sf m})$$
%	
%	
%		\begin{theorem}\label{Thm-Sobolev-Funk} 
%		%	{\rm (Vanishing of Rayleigh quotient on Funk space)}
%	Let  $\mathcal M_{\rm Funk}=(\mathbb B^n,F,{\sf m}_{BH})$ be the %$n$-dimensional 
%	Funk space. Then for every $p,q\geq  1$, one has that
%		\begin{equation}\label{optim-Funk} \mathcal R^{p}_q(\mathcal M_{\rm Funk})=0.
%	\end{equation}
%	\end{theorem}
%\noindent 
%The  proof of Theorem \ref{Thm-Sobolev-Funk} is based on a carefully chosen test function    with specific features of the Funk space.   
% 
% The  proof of Theorem \ref{Thm-Sobolev-Funk} is constructive, combining a particular choice of   test functions   with specific features of the Funk space. 
%Note that Sobolev inequalities Cavalletti--Mondino igazoltak. 

Given a Finsler space $\mathcal M=(M,F,{\sf m})$, the {\it reversibility} $$r_F=\sup\{F(-v):F(v)=1, v\in TM\}$$ 
measures the deflection of how far  is the   metric $F$ to be reversible,  
see Rademacher \cite{Rademacher}; it is clear that $r_F\geq 1$ and $r_F=1$ if and only if $F$ is reversible. At first sight, it is natural to expect the   Sobolev constants on $\mathcal M$ to depend on the reversibility $r_F$. 
Indeed, standard arguments typically control the Dirichlet energy   $\|F^*(-D|u|)\|_{L^p(M,{\sf m})}$ rather than the full energy $\|F^*(D u )\|_{L^p(M,{\sf m})}$ of the function $u:M\to \mathbb R,$ 
 see e.g.\ Ohta--Sturm \cite{Ohta-Sturm}   and Farkas--Krist\'aly--Varga \cite{FKV-Calculus}; the two quantities are related through the reversibility, i.e.,  $\|F^*(-D|u|)\|_{L^p(M,{\sf m})}\leq r_F\|F^*(D u )\|_{L^p(M,{\sf m})}$. Note, however, that this global estimate inevitably implies the loss of optimality in   Sobolev constants.
 In particular, given a Finsler manifold
 $\mathcal M=(M,F,{\sf m})$ having finite reversibility $r_F$ and ${\sf m}(M)=1,$ with weighted Ricci curvature ${\sf Ric}_N\geq K$   for  some $K\in \mathbb R$ and $N>1$ and diameter bounded from above by $D,$
  Ohta \cite[Theorem 1.3]{Ohta-2018} proved the isoperimetric inequality 
 $$P(E)\geq r_F^{-1}\mathcal J_{K,N,D}({\sf m}(E)),\ \ \forall E\subset M,$$
 where $P$ is the perimeter, and  $\mathcal J_{K,N,D}$ is the  isoperimetric profile function of Milman \cite{Milman}; see also Cavalletti--Mondino \cite{CM-Inventiones}. While the latter inequality is sharp in the reversible setting ($r_F=1$), see Milman \cite{Milman}, in the  non-reversible case  the factor $r_F>1$  can also be eliminated when $N=D=\infty$ and $K>0$, yielding  the sharp  Bakry--Ledoux isoperimetric  inequality  on non-reversible Finsler manifolds, see Ohta \cite{Ohta-CAG}. This particular observation  led  Ohta  to formulate the conjecture that  sharp functional inequalities on non-reversible Finsler manifolds should have the same constant as in the reversible case. Using the sharp results of   Cavalletti--Mondino \cite{CM-GT} from the reversible setting, Ohta \cite{Ohta-Indian} confirmed this conjecture for  the Poincar\'e--Lichnerowicz, logarithmic-Sobolev and Sobolev inequalities for non-reversible Finsler metrics with weighted Ricci-curvature ${\sf Ric}_N\geq K>0$. Note that the latter curvature restriction implies the compactness of the manifold. 
 
Due to these results, our  aim is to  establish  optimal Sobolev inequalities on Finsler manifolds having ${\sf Ric}_N\geq  0$, which can be viewed as a limiting case of the above framework; note that these objects are usually non-compact, and ${\sf Ric}_N $ contain both the flag curvature and  measure. 

To be more precise, if $\mathcal M=(M,F,{\sf m})$ is a forward complete  $n$-dimensional    Finsler space  $($possibly with locally convex boundary$)$ and $N\geq n$, the \textit{weighted $N$-Ricci curvature} is 
$${\sf Ric}_N(v)={\sf Ric}(v)+\dot {\sf S}_{\sf m}(v)-\frac{1}{N-n}{\sf S}_{\sf m}^2(v),\ \ \ v\in TM,$$ 
with the appropriate limiting interpretations; here,  $N$ stands for the synthetic dimension, 
  ${\sf Ric}$ is the usual Ricci curvature on $(M,F)$, while ${\sf S}_{\sf m}$ is the $S$-curvature,  measuring the variation of the measure ${\sf m}$  along geodesics. For details,   see \S   \ref{section-preliminaries}. 

If ${\sf Ric}_N\geq  0$, according to the Bishop--Gromov comparison principle, the function $R\mapsto \frac{{\sf m}(B^\pm_{x_0}(R))}{R^N}$ is non-increasing on $(0,\infty)$ for every $x_0\in M$, where $B^+_{x_0}(R)$  (resp. $B^-_{x_0}(R)$) is the forward (resp. backward) ball  in $M$ with center $x_0\in M$ and radius $R>0.$ Accordingly, we consider the \textit{asymptotic volume ratio on $\mathcal M$} as 
$${\sf AVR}_{\mathcal M}=\min\{{\sf AVR}_{\mathcal M}^+,{\sf AVR}_{\mathcal M}^-\},\ \  {\rm where}\ \ {\sf AVR}_{\mathcal M}^\pm=\lim_{R\to \infty}\frac{{\sf m}(B^\pm_{x_0}(R))}{\omega_N R^N},$$
with $\omega_N=\pi^{N/2}/\Gamma(N/2+1)$.   Note that both ${\sf AVR}_{\mathcal M}^+$ and ${\sf AVR}_{\mathcal M}^-$  are well-defined, i.e., independent on $x_0\in M$, and usually $ {\sf AVR}_{\mathcal M}^+\neq {\sf AVR}_{\mathcal M}^-$, see \S \ref{subsection-cone}. 
 For a fixed $x_0\in M$, we also consider the $N$-density
 $$\theta_N(x_0)= \min\{\theta_N^+(x_0),\theta_N^-(x_0)\},\ \  {\rm where}\ \  \theta_N^\pm(x_0)=\lim_{R\to 0}\frac{{\sf m}(B_{x_0}^\pm( R))}{  R^N}.$$ 
 The  synthetic dimension  $N$   \textit{saturates} if there exists $x_0\in M$ $($or $x_0\in \partial M$ if $M$ has locally convex boundary$)$ such that   $\theta_N(x_0)<\infty$. Note that if $x_0\in M$ is an interior point of $M$, then $\theta_N(x_0)<\infty$ forces $N=n$. Thus $N\geq n$ saturates typically in singular/boundary points;  
 % In contrast, if $x_0\in \partial M$, one typically has $N>n;$ 
  see e.g.\     convex cones in $\mathbb R^n$  endowed with monomial weights (with homogeneity $\alpha\geq 0$) having synthetic dimension $N=n+\alpha$, and $x_0$ taken to be their vertex, see \S \ref{subsection-cone}.

 The following result confirms Ohta's conjectural principle in the non-compact setting: for Finsler manifolds $\mathcal M$ with ${\sf Ric}_N\geq 0$, the reversibility \textit{does not} enter the optimal Sobolev constant,   the non-reversibility being encoded into the forward and backward asymptotic volume ratios.
%
%The following result confirms Ohta's conjecture for the Rayleigh quotient on non-compact Finsler manifolds $\mathcal M$ with ${\sf Ric}_N\geq 0$; namely,   the  reversibility in the optimal Sobolev inequality can be   avoided, the non-reversibility of the Finsler manifold being encoded into the forward and backward asymptotic volume ratios. 
More precisely, we have:

		\begin{theorem}\label{Thm-Sobolev-1} 
			%{\rm (Sharpness under ${\rm Ric}_N\geq 0$)}
			% {\rm (Optimal Sobolev inequality)}
				Let $\mathcal M=(M,F,{\sf m})$ be a not necessarily reversible $n$-dimensional  complete Finsler space $($possibly with locally convex boundary$)$ such that ${\sf Ric}_N\geq 0$ for some $N\in [n,\infty)$. Let $1<p<N $ and $p^*_N=\frac{pN}{N-p}$.  If ${\sf AVR}_{\mathcal M}>0$,   then 
						\begin{equation}\label{optim-sob-ineq}	\mathcal R^{p}_q(\mathcal M)=\begin{cases}
								0, & \text{if}\ \  q< p_N^*\  or\ \text{if}\  q> p_N^*\  \mbox{and} \
								 N\    \mbox{saturates}
								;\\
								{\sf AVR}_{\mathcal M}^\frac{1}{N}\cdot {\sf AT}(p,N),& \text{if}\ \  q= p_N^*.
								\end{cases}
								\end{equation}

%			\begin{equation}\label{optim-sob-ineq}
%				\int_M F^*(Du)^p {\rm d}{\sf m}\geq ({\sf AVR}_{\mathcal M}^{\rm min})^\frac{p}{N}\cdot {\sf AT}(p,N)\cdot \left(\int_M |u|^{p^\star}{\rm d}{\sf m}\right)^{p/p^\star},\ \ \forall u\in \dot W^{1,p}(M,F,{\sf m}).
%			\end{equation}
%		Moreover, the constant ${\sf AVR}_{\rm min}^\frac{p}{N}\cdot {\sf AT}(p,N)$ is optimal. 
		\end{theorem}
	
%$\theta_N(x_0)<\infty$: 	synthetic dimension is realized by a boundary singularity
	
	The non-critical case ($q\neq p_N^*$) mimics the scaling argument through some cut-off functions, combined with the fact that ${\sf AVR}_{\mathcal M}>0$. 
	
In the   critical case ($q= p_N^*$) we apply the theory of optimal mass transportation. Indeed, the  assumption ${\sf Ric}_N\geq 0$ implies the displacement convexity of a general internal energy functional. This convexity property is used for suitable measures -- the target being the Finslerian--Talentian-bubble --  producing the appearance of ${\sf AVR}_{\mathcal M}>0$ after a careful limiting argument.  In fact, we first prove a Sobolev inequality  involving the terms $\|F^*(-Du)\|_{L^p(M,{\sf m})}$ and ${\sf AVR}^+_{\mathcal M}>0$ for every non-negative function $u\in W_{p_N^*}^p(\mathcal M)$; in tandem, a similar result can be stated for the reverse metric $\overleftarrow F(x,y)=F(x,-y )$ with ${\sf AVR}^-_{\mathcal M}>0$. These inequalities are carefully combined to obtain the  general inequality $\mathcal R^{p}_{p_N^*}(\mathcal M)\geq {\sf AVR}_{\mathcal M}^\frac{1}{N}\cdot {\sf AT}(p,N)$. The  optimality   is obtained by  blow-downing the  Talentian-bubbles.
	When $\partial M\neq \emptyset$, the completeness of $M$ is understood in the metric setting, while the local convexity  of $\partial M$  is indispensable in Theorem \ref{Thm-Sobolev-1}; see Remark \ref{remark-locally-convex}.

	 Theorem \ref{Thm-Sobolev-1} shows  the efficiency of the  optimal mass transportation to prove optimal Sobolev inequalities on (not necessarily reversible) Finsler manifolds. The latter question was formulated  by Cordero-Erausquin--Nazaret--Villani \cite{CENV} in the Riemannian setting, and recently answered  by Krist\'aly \cite{Kristaly-OMT} for Riemannian manifolds with non-negative Ricci curvature endowed with the canonical measure. 
	We stress that an adaptation of  arguments from \cite{Kristaly-OMT} to Finsler manifolds requires more care, owing to the combined difficulties arising from non-reversibility and the presence of a general (non-canonical) measure ${\sf m}$.

	As  $p\to 1$, the 'forward' part of \eqref{optim-sob-ineq} reduces to the  optimal isoperimetric inequality
	\begin{equation}\label{Manini-ineq}
			\|F^*(-Du)\|_{L^1(M,{\sf m})}\geq N (\omega_N\,{\sf AVR}^+_{\mathcal M})^{1/N}\|u\|_{L^\frac{N}{N-1}(M,{\sf m})},\quad \forall u\in W^1_{1_N^*}(\mathcal M),\ u\geq 0,
	\end{equation}
	which is the veritable Finslerian version of Gromov's inequality \eqref{Gromov-inequality} under  ${\sf Ric}_N\geq 0$. The geometric form of inequality \eqref{Manini-ineq} -- for the characteristic function $u=\mathbbm{1}_{E}$, with $E\subset M$ open and bounded -- has been recently established by Manini \cite{Manini} via localization techniques \`a la Klartag \cite{Klartag}. Under mild assumptions,  the  equality is characterized in Manini \cite{Manini} (see also Cavalletti--Manini \cite{CM-JEMS} in the symmetric setting of ${\sf CD}(0,N)$ spaces): the set $E$ should be  $ B_{x_0}^+(R)$ for some $x_0\in M$ and $R>0,$ while the measure ${\sf m}$  disintegrates along geodesic rays emanating from $x_0$, having densities $N\omega_N {\sf AVR}^+_{\mathcal M}r^{N-1}{{\mathcal L}^1}_{\llcorner{[0,\infty)}}$.  
	
	 The   saturation  in the Rayleigh quotient  is described by the following rigidity result:
%	 ; for simplicity, we  consider the case ${\sf AVR}_{\mathcal M}={\sf AVR}_{\mathcal M}^+$: 

	\begin{theorem}\label{Thm-Sobolev-2-equality} 
		%{\rm (Rigidity  under ${\rm Ric}_N\geq 0$)}
		%{\rm (Equality in the optimal Sobolev inequality)} 
		Under the same assumptions as in Theorem \ref{Thm-Sobolev-1}, with $r_F<\infty$,     the following two statements are equivalent$:$
		\begin{itemize}
			\item[(i)] $\mathcal R^p_{p^*_N}(\mathcal M)={\sf AVR}_{\mathcal M}^\frac{1}{N}\cdot {\sf AT}(p,N)$ is achieved  by some nonzero function $u\in   W^p_{p^*_N}(\mathcal M)$$; $
%			\item[(ii)] There exists $x_0\in M$ $($or $x_0\in \partial M$ if $M$ has boundary$)$ such that the following   assertions simultaneously hold$:$
%			\begin{itemize}
%				\item[(a)] if ${\sf AVR}_{\mathcal M}={\sf AVR}_{\mathcal M}^+$ $($resp. ${\sf AVR}_{\mathcal M}={\sf AVR}_{\mathcal M}^-$$),$ then $\mathcal M=(M,F,{\sf m})$ is a forward $($resp.  backward$)$ $N$-volume cone with pole $x_0$,  i.e., for every $R>0$, 
%				$$ \begin{cases}
%					{\sf m}(B_{x_0}^+(R))={\sf AVR}_{\mathcal M}^+ \omega_N R^N, & \text{if}\ \ {\sf AVR}_{\mathcal M}={\sf AVR}_{\mathcal M}^+;\\
%					{\sf m}(B_{x_0}^-(R))={\sf AVR}_{\mathcal M}^- \omega_N R^N, & \text{if}\ \ {\sf AVR}_{\mathcal M}={\sf AVR}_{\mathcal M}^-;\\
%					%					\text{any of the above functions}, & \text{if}\ \ {\sf AVR}_{\mathcal M}^+={\sf AVR}_{\mathcal M}^-.
%				\end{cases}$$
%				
%				
%				\item[(b)] for some $\lambda,\mu>0$, the extremizer has the form
%%				 $u(x)=-(\lambda+\mu d_F(x_0,x)^{p'})^\frac{p-N}{N} $ .  
%				$$u(x)=\begin{cases}
%					-(\lambda+\mu d_F(x_0,x)^{p'})^\frac{p-N}{N}, & \text{if}\ \ {\sf AVR}_{\mathcal M}={\sf AVR}_{\mathcal M}^+;\\
%					(\lambda+\mu d_F(x,x_0)^{p'})^\frac{p-N}{N}, & \text{if}\ \ {\sf AVR}_{\mathcal M}={\sf AVR}_{\mathcal M}^-.\\
%%					\text{any of the above functions}, & \text{if}\ \ {\sf AVR}_{\mathcal M}^+={\sf AVR}_{\mathcal M}^-.
%				\end{cases}$$
%			\end{itemize}
			\item[(ii)] There exists $x_0\in M$ $($or $x_0\in \partial M$ if $M$ has boundary$)$   such that the following   assertions simultaneously hold$:$
		\begin{itemize}
			\item[(a)] $\mathrm{Cut}(x_0)=\emptyset$ and if ${\sf AVR}_{\mathcal M}={\sf AVR}_{\mathcal M}^+$ $($resp. ${\sf AVR}_{\mathcal M}={\sf AVR}_{\mathcal M}^-$$)$,  then along any forward $($resp. backward$)$ geodesic $\gamma:[0,\infty)\to M$ 
			 with $\gamma(0)=x_0,$ the flag curvature vanishes, ${\sf K}(\dot \gamma(t),\cdot)=0$ $($resp. ${\sf K}(-\dot \gamma(t),\cdot)=0$$),$ and the $S$-curvature verifies ${\sf S}_{\sf m}(\gamma(t),\dot \gamma(t))=\frac{n-N}{t}$ $($resp. ${\sf S}_{\sf m}(\gamma(t),-\dot \gamma(t))=\frac{N-n}{t}$$)$ for every $t>0;$
			
			\item[(b)] for some $\lambda,\mu>0$, one of the following two
			alternatives holds$:$  \begin{itemize}
				\item[(b1)] $u\leq 0$, ${\sf AVR}_{\mathcal M}={\sf AVR}_{\mathcal M}^+$ and $u(x)=-(\lambda+\mu d_F(x_0,x)^{p'})^\frac{p-N}{p}$$;$ 
				\item[(b2)] $u\geq 0$, ${\sf AVR}_{\mathcal M}={\sf AVR}_{\mathcal M}^-$ and $u(x)=(\lambda+\mu d_F(x,x_0)^{p'})^\frac{p-N}{p}.$
			\end{itemize}
%			%				 $u(x)=-(\lambda+\mu d_F(x_0,x)^{p'})^\frac{p-N}{N} $ .  
%			$$u(x)=\begin{cases}
%				-(\lambda+\mu d_F(x_0,x)^{p'})^\frac{p-N}{p}, & \text{if}\ \ {\sf AVR}_{\mathcal M}={\sf AVR}_{\mathcal M}^+;\\
%				(\lambda+\mu d_F(x,x_0)^{p'})^\frac{p-N}{p}, & \text{if}\ \ {\sf AVR}_{\mathcal M}={\sf AVR}_{\mathcal M}^-.\\ 
%			\end{cases}$$
		\end{itemize}
		\end{itemize}
	\end{theorem}
%	We notice that the negative sign of the extremizer is not accidental; however, in the  case when  ${\sf AVR}_{\mathcal M}={\sf AVR}_{\mathcal M}^-$, the extremizer in the Rayleigh quotient is $u(x)=(\lambda+\mu d_{\overleftarrow F}(x_0,x)^{p'})^\frac{p-N}{N} $ for some $\lambda,\mu>0$ together with the corresponding curvature constraints from (ii) along any backward geodesic arriving into a specific point $x_0$.    

The limiting arguments in Theorem \ref{Thm-Sobolev-1}   prevent the identification of the equality cases in order to establish   Theorem \ref{Thm-Sobolev-2-equality}. Instead, adapting recent results of Nobili--Violo \cite{NV-Calc},   based on the optimal isoperimetric inequality \eqref{Manini-ineq}, we apply a suitable rearrangement by establishing a forward/backward P\'olya--Szeg\H o inequality; this   provides an alternative approach to proving Theorem \ref{Thm-Sobolev-1}. 
If the Rayleigh quotient $\mathcal R^p_{p^*_N}(\mathcal M)={\sf AVR}_{\mathcal M}^\frac{1}{N}\cdot {\sf AT}(p,N)$ is achieved  by some   $u\in W^p_{p^*_N}(\mathcal M)\setminus \{0\}$, it turns out  that $u$ should be sign-preserving and  if ${\sf AVR}_{\mathcal M}={\sf AVR}_{\mathcal M}^+$ $($resp.\ ${\sf AVR}_{\mathcal M}={\sf AVR}_{\mathcal M}^-$$),$ then $\mathcal M=(M,F,{\sf m})$ is a forward $($resp.\  backward$)$ $N$-volume cone with pole $x_0$,  i.e.,  ${\sf m}(B_{x_0}^\pm(R))={\sf AVR}_{\mathcal M}^\pm \omega_N R^N$  for every $R>0$; in particular, $N$ saturates as $\theta_N(x_0)=\omega_N{\sf AVR}_{\mathcal M}<\infty$.  
%				$$ \begin{cases}
%					{\sf m}(B_{x_0}^+(R))={\sf AVR}_{\mathcal M}^+ \omega_N R^N, & \text{if}\ \ {\sf AVR}_{\mathcal M}={\sf AVR}_{\mathcal M}^+;\\
%					{\sf m}(B_{x_0}^-(R))={\sf AVR}_{\mathcal M}^- \omega_N R^N, & \text{if}\ \ {\sf AVR}_{\mathcal M}={\sf AVR}_{\mathcal M}^-.\\
%					\end{cases}$$
The curvature rigidity in (ii)/(a) follows from a careful analysis of the $N$-volume cone character of $\mathcal M$, based on Ohta's proof of the Bishop--Gromov volume comparison   \cite{Ohta-Calc, Ohta-book}; see Proposition \ref{N-volume-cone-prop}. 
% which is of independent interest. 
%Moreover, the scaling factor $n-N$ in the hyperbolic form of the $S$-curvature ${\sf S}_{\sf m}$ reflects the disintegration feature of the measure ${\sf m}$  along geodesic rays emanating from the pole $x_0$ at the polynomial rate $r^{N-1}$ (see Remark \textbf{???}); this fact is in perfect concordance with Manini's result \cite{Manini} in the case $p=1$.
 The shape of the extremizer $u$ in (ii)/(b) is dictated by the equality case in the forward/backward P\'olya--Szeg\H o inequality, while the signs appearing therein are essential.
  
	We also notice that -- contrary to our initial expectations -- the rigidity provided by  Theorem \ref{Thm-Sobolev-2-equality} does not necessarily imply the vanishing of the $S$-curvature. On the one hand, this fact explains the characterization of the equality in the optimal Sobolev inequality  on $(\Sigma,H,w\mathcal L^n_{\llcorner{\Sigma}})$, where $\Sigma\subseteq \mathbb R^n$ is an open convex cone, $H$ is a Minkowski norm (i.e., $H\in C^\infty(\mathbb R^n\setminus\{0\})$, positive and strongly convex)    and $w:\Sigma \to \mathbb R_+$ is a homogeneous  weight function  of degree $\alpha$   such that $w^{1/\alpha}$ is concave,   see Balogh--Guti\'errez--Krist\'aly \cite{BGK}, as well as Cabr\'e--Ros-Oton \cite{Cabre-JDE} for   monomial weights. Indeed, it turns out that on $(\Sigma,H,w\mathcal L^n_{\llcorner{\Sigma}})$ the synthetic dimension is $N=n+\alpha$, while the $S$-curvature is $-\alpha/t$ along any (unitary parameterized) semi-line  starting from a pole of $\Sigma$; for details, see \S \ref{subsection-cone}. On the other hand, when   $\mathcal M=(M,F,{\sf m}_{BH})$ is a Berwald space endowed with the Busemann--Hausdorff measure ${\sf m}_{BH}$, the $S$-curvature vanishes, $N=n$,  and the rigidity in Theorem \ref{Thm-Sobolev-2-equality} implies that   $(M,F)$ is isometric to a Minkowski space.

Using a similar approach  as above, we can establish the \textit{optimal $L^p$-logarithmic Sobolev inequality} on not necessarily reversible Finsler manifolds (possibly with boundary) satisfying ${\sf Ric}_N\geq 0$, along with a characterization of the equality case. To state this, we recall for every $p,N>1$ the value 
		$$		\mathcal
	L_{p,N}={\small
		\frac{p}{N}\left(\frac{p-1}{e}\right)^{p-1}\left(\omega_N{\Gamma\left(\frac{N}{p'}+1\right)}\right)^{-\frac{p}{N}}}.  
	$$

		\begin{theorem}\label{log-Sobolev-main}  Let $\mathcal M=(M,F,{\sf m})$ be a not necessarily reversible $n$-dimensional  complete Finsler space $($possibly with locally convex boundary$)$ such that ${\sf Ric}_N\geq 0$ for some $N\in [n,\infty)$ and let $p>1$. 
		If ${\sf AVR}_{\mathcal M}>0$, for every  
		$ u\in W^{p}_p(\mathcal M)$ with $\|u\|_{L^p(M,{\sf m})}=1,$ one has  
		\begin{equation}\label{LSI}
			\int_{M}|u|^p\log |u|^p{\rm d} {\sf m}\leq \frac{N}{p}\log\left(\mathcal
			L_{p,N}{\sf AVR}_ {\mathcal M}^{-\frac{p}{N}}\|F^*(Du)\|_{L^p(M,{\sf m})}^p \right), 
		\end{equation}
		and  the constant $\mathcal
		L_{p,N}{\sf AVR}_ {\mathcal M}^{-\frac{p}{N}}$   is optimal. Moreover,  if  $r_F<\infty,$  the following  statements are equivalent$:$
		\begin{itemize}
			\item[(i)] Equality holds in \eqref{LSI}  for some  function $u\in   W^p_{p }(\mathcal M)$$; $
			\item[(ii)] There exists $x_0\in M$ $($or $x_0\in \partial M$ if $M$ has boundary$)$  such that the following   assertions simultaneously hold$:$
			\begin{itemize}
				\item[(a)] verbatim, as in Theorem \ref{Thm-Sobolev-2-equality} {\rm (ii)/(a)};
				
				\item[(b)] for some $\lambda>0$, if $C(\lambda):=\lambda^\frac{N}{pp'}{\sf AVR}_{\mathcal M}^{-\frac{1}{p}}\left(\omega_N    \Gamma \left(\frac{N}{p'}+1\right)\right)^{-\frac{1}{p} }$, then   one of the following two
				alternatives holds$:$
				\begin{itemize}
					\item[(b1)] $u\leq 0$, ${\sf AVR}_{\mathcal M}={\sf AVR}_{\mathcal M}^+$ and $u(x)=-C(\lambda)e^{-\frac{\lambda}{p} d_F(x_0,x)^{p'}};$
					\item[(b2)]   $u\geq 0$, ${\sf AVR}_{\mathcal M}={\sf AVR}_{\mathcal M}^-$ and $u(x)=C(\lambda)e^{-\frac{\lambda}{p} d_F(x,x_0)^{p'}}.$
				\end{itemize}
%				%				 $u(x)=-(\lambda+\mu d_F(x_0,x)^{p'})^\frac{p-N}{N} $ .  
%				$$u(x)=\lambda^\frac{N}{pp'}{\sf AVR}_{\mathcal M}^{-\frac{1}{p}}\left(\omega_N    \Gamma \left(\frac{N}{p'}+1\right)\right)^{-\frac{1}{p} }\cdot\begin{cases}
%					-e^{-\frac{\lambda}{p} d_F(x_0,x)^{p'}}, & \text{if}\ \ {\sf AVR}_{\mathcal M}={\sf AVR}_{\mathcal M}^+;\\
%					e^{-\frac{\lambda}{p} d_F(x,x_0)^{p'}}, & \text{if}\ \ {\sf AVR}_{\mathcal M}={\sf AVR}_{\mathcal M}^-.\\ 
%				\end{cases}$$
			\end{itemize}
			
		\end{itemize}
	\end{theorem}
The $L^p$-logarithmic Sobolev inequality in the Euclidean space $\mathbb R^n$  has a long history, starting with Gross \cite{Gross},  del Pino--Dolbeault \cite{Dol-delPino} and Gentil \cite{Gentil}, with the characterization of   equality only in the case   $1<p<n$. 
On convex cones of $\mathbb R^n$ with  homogeneous gauge function $w$ of degree $\alpha\geq 0$ such that $w^{1/\alpha}$ is concave,  the $L^p$-logarithmic Sobolev inequality has been established by Balogh--Don--Krist\'aly  \cite{BDK} for general $p,n>1$, jointly with the characterization of the equality case; this result will play a crucial role in the second part   of Theorem \ref{log-Sobolev-main}.   

A remarkable application of Theorem \ref{log-Sobolev-main} is the optimal hypercontractivity estimate in the heat flow  on Finsler manifolds with ${\sf Ric}_N\geq 0$, jointly with the characterization of the equality case. To state it, we consider  a not necessarily reversible $n$-dimensional  complete Finsler space $\mathcal M=(M,F,{\sf m})$ and  a (weak) solution $u\in L^2([0,T],H_0^1(M))\cap H^1([0,T],H^{-1}(M))$ to    the heat equation \begin{equation}\label{heat-flow}
	\partial_t u= \bm \Delta_{\mathcal M} u,  
\end{equation} 
	where $\bm \Delta_{\mathcal M}$ is the non-linear Finsler--Laplace operator on $\mathcal M$.   
	
	For every $2\leq p\leq q$, we introduce the number
	\begin{equation}\label{M-p-q}
		{\sf M}(p,q)=\frac{p^{\frac{1}{p}} {q' }^{\frac{1}{q'}}}{{p' }^{\frac{1}{p'}} q^{\frac{1}{q}}}\left(\frac{1}{p}-\frac{1}{q}\right)^{\frac{1}{p}-\frac{1}{q}},
	\end{equation}
		with the  limiting conventions   ${\sf M}(p,p)=1$ and  ${\sf M}(p,\infty)=(p')^{-\frac{1}{p'}}$.

	\begin{theorem}\label{hypercontractivity-main-theorem}  Let $\mathcal M=(M,F,{\sf m})$ be a not necessarily reversible $n$-dimensional  complete Finsler space $($possibly with locally convex boundary$)$ such that ${\sf Ric}_N\geq 0$ for some $N\in [n,\infty)$ and ${\sf AVR}_{\mathcal M}^->0$. Let   $T>0$ and   $u\in L^2([0,T],H_0^1(M))\cap H^1([0,T],H^{-1}(M))$  be a non-negative solution to the heat equation \eqref{heat-flow} with the initial datum $u_0=u(0,\cdot)\in L^1(M,{\sf m})\cap L^\infty(M)$, $u_0\geq 0$. Then for every $2\leq p\leq q\leq +\infty$ and $0\leq T_0<  T_1\leq T$, one has that 
		\begin{equation}\label{hyper-estimate}
			 \|u(T_1,\cdot)\|_{L^q(M,{\sf m})}\leq \frac{{\sf M}(p,q)^\frac{N}{2}({\sf AVR}_{\mathcal M}^-)^{\frac{1}{q}-\frac{1}{p}}}{(4\pi (T_1-T_0))^{\frac{N}{2}(\frac{1}{p}-\frac{1}{q})}}\|u(T_0,\cdot)\|_{L^p(M,{\sf m})}
		\end{equation}
		and  the constant ${\sf AVR}_{\mathcal M}^-$   is optimal. Moreover,  if  $r_F<\infty$,    the following  statements are equivalent$:$
		\begin{itemize}
			\item[(i)] Equality holds in \eqref{hyper-estimate}  for some $ p< q$,  $ T_0<  T_1 $ and positive function $u\in AC([T_0,T_1], L^2(M,{\sf m}))$$; $
			\item[(ii)] There exists $x_0\in M$ $($or $x_0\in \partial M$ if $M$ has boundary$)$   such that the following   assertions simultaneously hold$:$
			\begin{itemize}
				\item[(a)] $\mathrm{Cut}(x_0)=\emptyset$ and along any    backward  geodesic $\gamma:[0,\infty)\to M$ 
				with $\gamma(0)=x_0,$ one has   ${\sf K}(-\dot \gamma(t),\cdot)=0$  and  ${\sf S}_{\sf m}(\gamma(t),-\dot \gamma(t))=\frac{N-n}{t}$ for every $t>0;$ 
				
				\item[(b)] for some $\lambda>0$, one has
				%				 $u(x)=-(\lambda+\mu d_F(x_0,x)^{p'})^\frac{p-N}{N} $ .  
				\begin{equation}\label{u-t-x-precise}
					u(t,x)=\lambda (4t+K_0)^{-\frac{N}{2}}e^{-\frac{d_F(x,x_0)^2}{4t+K_0}}  ,\quad (t,x)\in [T_0,T_1]\times M,
				\end{equation}
								 where $$K_0=4\frac{(p-1)qT_1-(q-1)pT_0}{q-p}.$$
			\end{itemize}
			
		\end{itemize}
	\end{theorem}

A similar hypercontractivity estimate to Theorem \ref{hypercontractivity-main-theorem} on ${\sf RCD}(0,N)$ spaces has been established by Honda--Krist\'aly--P\^irvuceanu \cite{HKP}. While in \cite{HKP} the arguments closely follow  those from Bakry--Gentil--Ledoux    \cite[Theorem 7.1.2]{BGL}, we provide  here a direct and self-contained proof. It is worth to emphasize that Finsler manifolds with ${\sf Ric}_N \geq  0$ are ${\sf CD}(0,N)$ spaces but are not ${\sf RCD}(0,N)$ spaces unless they are Riemannian. In addition, in the ${\sf RCD}(0,N)$  setting the presence of heat kernels  is crucial role in \cite{HKP}, which is not available on   Finsler manifolds, due to the  nonlinear character of the heat equation \eqref{heat-flow}.

Before concluding this section, we recall that if $u$ is a  positive solutions to the heat
equation \eqref{heat-flow} on a Finsler space $\mathcal M=(M,F,{\sf m})$ with ${\sf Ric}_N\geq 0$,
% the famous \textit{Li--Yau gradient estimate} holds: 
%\begin{equation}\label{Li-Yau}
%	\frac{F^*(Du(t,x))^2}{u(t,x)^2}-\frac{\partial_t u(t,x)}{u(t,x)}\leq \frac{N}{2t},\quad (t,x)\in (0,\infty)\times M;
%\end{equation}
%the original form of \eqref{Li-Yau}  is due to Li--Yau \cite{LY} on Riemannian manifolds, while the present Finsler version is provided by 
%Ohta--Sturm \cite[Theorem 4.4]{Ohta-Sturm-adv}. 
then the famous \textit{Li--Yau-type Harnack inequality} holds, i.e., 
$$u(s,x)\leq u(t,y)\left(\frac{t}{s}\right)^{N/2}e^\frac{d_F(y,x)^2}{4(t-s)},\quad \forall 0<s<t\ {\rm and}\ x,y\in M,$$
see   Li--Yau \cite{LY} for Riemannian manifolds, and Ohta--Sturm \cite[Theorem 4.5]{Ohta-Sturm-adv} for the   Finsler version. 
 
 In the setting of   a complete $n$-dimensional Riemannian space  $\mathcal M_g=(M,g,{\rm vol}_g)$ with ${\sf Ric}_g\geq 0$ and ${\sf AVR}_{{\mathcal M}_g}>0$, the
 latter estimate together with the two-sided Gaussian bounds for the heat kernel $p_t(x,y)$ of $(M,g)$, cf.\ Li--Yau \cite{LY}, the Bishop--Gromov comparison principle and a suitable interpolation yield that 
 $$\|u(t,\cdot)\|_{L^q(M,{\rm vol}_g)}\leq C\frac{ {\sf AVR}_{\mathcal M_g}  ^{\frac{1}{q}-\frac{1}{p}}}{t^{\frac{n}{2}(\frac{1}{p}-\frac{1}{q})}}\|u(0,\cdot)\|_{L^p(M,{\rm vol}_g)},\quad \forall t>0\ {\rm and }\ 1\leq p\leq q\leq +\infty,$$
for some non-explicit constant $C=C(n,p,q)>0$. The advantage of  Theorem \ref{hypercontractivity-main-theorem} is that it provides the optimal form of the latter estimate for $2\leq p\leq q\leq +\infty$ (characterizing also the equality case), even  on not necessarily reversible Finsler manifolds, where -- as far as we know -- two-sided  heat kernel estimates are not available. \\ 

%In addition, in the general Finsler setting of Theorem \ref{hypercontractivity-main-theorem}, the characterization of the equality case in \eqref{hyper-estimate} suggests a rigidity phenomenon in the Li--Yau gradient estimate \eqref{Li-Yau}. 

\textit{Organization of the paper.}	In \S\ref {section-preliminaries} we recall/prove those preliminary results which are crucial in our investigations. In \S \ref{section-3} we prove the optimal Finsler--Sobolev inequality (Theorems \ref{Thm-Sobolev-1} \& \ref{Thm-Sobolev-2-equality}),  while \S \ref{section-4} is devoted to the proof of the optimal logarithmic-Sobolev inequality  (Theorem \ref{log-Sobolev-main}). In \S \ref{section-5} we prove the optimal  hypercontractivity estimate (Theorem \ref{hypercontractivity-main-theorem}). Finally, in \S \ref{section-6}, some examples are presented.

	\section{Preliminaries}\label{section-preliminaries}
	
		\subsection{Elements from Finsler geometry} 
		\subsubsection{Finsler spaces}
			Let $M$ be an $n$-dimensional connected $C^\infty$-manifold,  possibly with   a Lipschitz  boundary $\partial M$,  
		and $T{M}=\bigcup_{x \in M}T_{x} {M}$ be its tangent bundle.
		We say that $(M,F)$ is a \emph{Finsler manifold} if $F:T{M}\to [0,\infty)$ satisfies:
		\begin{enumerate} 
			\item   $F\in C^{\infty}(T{M}\setminus\{ 0 \});$
			\item  $F(x,\lambda y)=\lambda F(x,y)$ for all $\lambda\geq 0$ and $(x,y)\in T{M};$
			\item  $g_{ij}(x,y) :=[\frac12F^{2}]_{y^{i}y^{j}}(x,y)$
			is positive definite for all $(x,y)\in T{M}\setminus\{ 0 \}$,
			where $y=\sum_{i=1}^n y^i\frac{\partial}{\partial x^i}|_x$. We will denote by $g_y$ the inner product
			on $T_xM$ induced by the matrix elements $g_{ij}(x,y)$. 
		\end{enumerate}
				The  {\it reverse Finsler metric} of $(M,F)$ is defined as $\overleftarrow{F}(x,y):=F(x,-y)$. The {\it reversibility} is defined as $$r_F=\sup\{\overleftarrow F(v):F(v)=1, v\in TM\}.$$  %Clearly, a curve is a $p$-forward (resp., $p$-backward) absolutely continuous curve in $(M, \overleftarrow{F})$ if and only if it is a $p$-backward (resp., $p$-forward) absolutely continuous curve in $(M,F)$.

		Denote by $T^*{M}$ the cotangent bundle, and define the {\it dual metric} $F^*$ by
		\begin{equation*}\label{dualFinslermetric}
			F^*(x,\xi) :=\sup_{y \in T_x{M} \setminus \{0\}} \frac{\xi(y)}{F(x,y)}, \ \ \xi \in T_x^* M,.
			%\tag{2.2}\label{2.2newineq}
		\end{equation*}
	where $\xi(y)$   denotes the canonical pairing between $T_x^*M$ and $T_xM$.
				Note that
		\begin{equation}\label{Cauchy-Schwarz}
				\xi(y) \leq F^*(x,\xi) F(x,y) \quad \text{for all}\ y\in T_x{M},\, \xi \in T^*_x{M}.
		\end{equation}

		The {\it Legendre transformation} $\mathfrak{L} : T_x{M} \rightarrow T_x^*{M}$ is defined by
		\begin{equation*}
			\mathfrak{L}(x,y):=\sum_{i,j=1}^n g_{ij}(x,y) y^i \,{\dd}x^j \ \text{ if } y\neq0,
			\qquad
			\mathfrak{L}(x,y):=0 \ \text{ if } y=0.
		\end{equation*}
	Note that  $\mathfrak{L}:T{M}\setminus\{0\}\rightarrow T^*{M}\setminus\{0\}$ is a diffeomorphism
		and $F^*(\mathfrak{L}(x,y))=F(x,y)$ for any $(x,y)\in T{M}$.
		For a $C^1$-function $f : M \rightarrow \mathbb{R}$,
		the {\it gradient vector field} of $f$ is defined as $\bm \nabla  f := \mathfrak{L}^{-1}(D f)$; in particular, $F (\bm \nabla f)=F^*(Df)$.
		We remark that $\bm \nabla$ is nonlinear (unless $F$ is Riemannian), i.e.,
		$\bm \nabla(f+h)\neq\bm \nabla f+\bm \nabla h$.

		Let $C^\infty_{\rm loc}([0,1];{M})$ denote the set of piecewise smooth curves defined on $[0,1]$.
		Given $\gamma\in C^\infty_{\rm loc}([0,1];{M})$, its \textit{length} is defined as
		\[
		L_F(\gamma):=\int^1_0 F\bigl(\gamma(t),\dot \gamma (t) \bigr) \,{\dd}t.
		\]
		The {\it distance function} $d_F:M\times M \rightarrow [0,\infty)$ is given by
		\[
		d_{F}(x_0,x_1):=\inf\{L_F(\gamma) : \gamma\in C^\infty_{\rm loc}([0,1];{M}),\,
		\gamma(0)=x_0,\, \gamma(1)=x_1\}.
		\]
		 Although $d_F$ satisfy the triangle inequality and $d_F(x_0,x_1)=0$ if and only if $x_0=x_1$, usually $d_{F}(x_0,x_1)\neq d_{F}(x_1,x_0)=d_{\overleftarrow F}(x_0,x_1)$. For $x_0\in M$ and $R>0$, the \textit{forward} and \textit{backward open balls} are 
		 $$B_{x_0}^+(R)=\{y\in M:d_F(x_0,y)<R\},\quad B_{x_0}^-(R)=\{y\in M:d_F(y,x_0)<R\}.$$
		 
		 Given $x_0\in M$, the \textit{eikonal equation}  states that for almost every $x\in M$,  
	 \begin{equation}\label{eikonal}
	 	F(x,\bm \nabla d_F(x_0,\cdot)(x))=F^*(x,D d_F(x_0,\cdot)(x))=1;
	 \end{equation}
		more precisely, the latter relation is valid for every $x\in M\setminus (\{x_0\}\cup {\rm Cut}(x_0))$, where ${\rm Cut}(x_0)$ stands for the cut locus of $x_0$. 
		
		A \emph{geodesic} $\gamma:[0,1] \rightarrow M$ is understood as a constant speed, locally minimizing curve,
		in the sense that $d_F(\gamma(s),\gamma(t)) =[(t-s)/(b-a)] \cdot d_F(\gamma(a),\gamma(b))$
		for sufficiently short intervals $[a,b]$ and any $a \le s<t \le b$. 	Analytically,  a smooth curve $t\mapsto\gamma(t)$ in $(M, F)$ is a geodesic, if it satisfies
		\begin{equation}\label{geodesequ}
		D_{\dot \gamma}^{\dot \gamma} {\dot \gamma}=0,
%			\frac{{\dd}^2 \gamma^i}{{\dd}t^2} + 2 G^i \left(\gamma, \frac{{\dd} \gamma}{{\dd}t}\right) = 0,
		\end{equation}
		where $D ^{\dot \gamma}$ stand for the covariant derivative with reference vector $\dot \gamma$, see Bao--Chern--Shen \cite{BCS}. 
		%		$ G^i = G^i(x,y)$ are the {\it geodesic coefficients }(or  {\it spray coefficients})
%		given by
%		\begin{equation}\label{goedcooff}
%			G^i = \frac{1}{4} g^{ij} \left\{ \frac{\partial^2 F^2}{\partial x^l  \partial y^j} y^l -\frac{\partial F^2}{\partial x^j} \right\},
%		\end{equation}
%		$(g^{ij})$ being the inverse matrix of $(g_{ij})$.
		In the sequel, we  use $\gamma_y(t)$, $t\geq 0,$ to denote the geodesic with initial velocity $y$.
	
		A Finsler manifold $(M,F)$ is said to be {\it forward} (resp., {\it backward}) {\it complete}
		if 	so is $(M,d_F)$, i.e.,
		for every sequence $(x_k)_k\subset M$ such that, for every $\varepsilon>0$, there exists $k_\varepsilon\in\mathbb N$ with
		$
		d_F(x_m,x_l)<\varepsilon
		$ (resp. $
		d_F(x_l,x_m)<\varepsilon
		$)
		for all $l>m>k_\varepsilon$, the sequence $(x_k)_k$ converges. $(M,F)$ is \textit{complete} if it is both forward and backward complete.  When the manifold is without boundary, forward (resp. backward) completeness is equivalent to the extensibility of
		 every geodesic segment defined on $[0,1]$ (resp., $[-1, 0]$)
		to   $[0, \infty)$ (resp., $(-\infty, 0]$). 
		
		If $(M,F)$ is a Finsler manifold with boundary $\partial M$, we  say that $\partial M$  is \textit{locally
		  convex} if   for all points $x,y\in M\setminus \partial M$  and for all minimal geodesic $\gamma$
		joining $x$ to $y,$ we have that  $\gamma$ does not touch the boundary.

%		According to the Hopf--Rinow theorem,
%		the closure of a forward ball (with a finite radius) is compact if $(M,F)$ is forward complete.
	%	However, this is not the case for backward balls (see Example \ref{funkmetricsapce}).

	 A \textit{Finsler space} is a triplet $\mathcal M=(M,F,{\sf m})$, where $(M,F)$ is a Finsler manifold and ${\sf m}$ is a positive measure on $M$. 
		In local coordinates ($x^i$) in $M$, ${\sf m}$ can be expressed as
		\begin{equation*}%\label{measure_m}
			{\dd}{\sf m}(x) = \sigma(x)\, {\dd}x^1 \cdots {\dd}x^n.
		\end{equation*}
		In particular, the {\it Busemann--Hausdorff measure}
		${\sf m}_{ BH }$ is defined by
		\[
		{\sf m}_{ BH }({\dd}x):=\frac{{\rm vol}(\mathbb{B}^{n})}{{\rm vol}(B_xM)}\, {\dd}x^1 \cdots {\dd}x^n,
		\]
		where $B_xM:=\{y\in T_xM: F(x,y)<1\}$. For further use, let $S_xM:=\{y\in T_xM: F(x,y)=1\}.$ 
		
	The \textit{Finsler--Laplace operator} $\bm \Delta_{\mathcal M}$ of a function $u\in H_{\rm loc}^1(M)$ on the Finsler space $\mathcal M=(M,F,{\sf m})$ is defined as 
	\begin{equation}
		\label{Laplacian}
		\int_M v \bm \Delta_{\mathcal M}u\, {\rm d}{\sf m}= - \int_M Dv(\bm \nabla u)\, {\rm d}{\sf m},\quad \forall v\in C_0^\infty(M).
	\end{equation}
 The Alexandrov Laplacian is denoted by $\bm \Delta_A$. 
	
		\subsubsection{Weighted $N$-Ricci curvature}

		For each  $y\in  T_xM$ and any plane $  \Pi=\text{span}\{y,v\}\subset T_xM$, the
		\emph{flag curvature}  is defined by
		\begin{eqnarray}
			\nonumber
		{\sf K} (y;\Pi):=\mathbf{\sf K}(y;v):=
		\frac{g_{y} (R^y(v,y)y,v)}{g_{y}(y,y)\,g_{y}(v,v) -g_{y}(y,v)^2},
%		
%		\frac{g_{im}(x,y)\,R^i_{\ k}(x,y)\, v^kv^m}{F^2(x,y)\,g_{ij}(x,y)\,v^iv^j-[g_{ij}(x,y)\,y^iv^j]^2},
		\end{eqnarray}
		where $R^y$ is the curvature tensor, see Bao--Chern--Shen \cite{BCS}. 
%		\begin{eqnarray}
%			\nonumber R^i_{\ k}=2\frac{\partial G^i}{\partial
%				x^k}-y^j\frac{\partial^2G^i}{\partial x^j\partial
%				y^k}+2G^j\frac{\partial^2G^i}{\partial y^j\partial
%				y^k}-\frac{\partial G^i}{\partial y^j}\frac{\partial G^j}{\partial
%				y^k},
%		\end{eqnarray}
%	$G^i$ being the geodesic coefficients from \eqref{goedcooff}. 
		If $F$ is a  Riemannian metric, the flag curvature reduces
		to the sectional curvature. 
		
		A smooth vector field $J$ along a geodesic $\gamma$ is a  \textit{Jacobi field} if it satisfies 
		\begin{equation}\label{Jacobi-eqn} 
			D_{\dot \gamma}^{\dot \gamma}D_{\dot \gamma}^{\dot \gamma} J+R^{\dot \gamma}(J,\dot \gamma)\dot \gamma=0.
		\end{equation}

		The {\it Ricci curvature} of $y\in T_xM\backslash\{0\}$ is defined by
		\[
		 {\sf Ric}(y):= {F^2(x,y)}\sum_{i} {\sf K}(y;e_i),
		\]
		where $\{e_1,\ldots,e_{n-1}\}$ is a $g_y$-orthonormal basis of $y^\perp$ for $T_xM$.

		The {\it distortion} $\tau$ and the {\it S-curvature} ${\sf S}_{\sf m}$ of $\mathcal M=(M,F,{\sf m})$ are defined  by
		\begin{equation*}%\label{distsdef}
			\tau(x,y):= \log \Biggl[ \frac{\sqrt{\det g_{ij}(x,y)}}{\sigma(x)} \Biggr], \qquad
			{\sf S}_{\sf m}(x,y):=\left.\frac{\dd}{{\dd}t}\right|_{t=0}\tau\bigl( \gamma_y(t), \dot{\gamma}_y(t) \bigr), \quad y\in T_xM \setminus \{0\},
		\end{equation*}
		respectively, where  $t\mapsto \gamma_y(t)$ is the geodesic with $\dot \gamma _y(0)=y$.
		
		According to Ohta \cite{Ohta-book}, the {\it weighted $N$-Ricci curvature} ${\sf Ric}_N$ is defined as follows:
		given $N\in [n,\infty]$, for any unit vector $y\in SM:=\bigcup_{x\in M}S_xM$,
		\[
		{\sf Ric}_N(y):=\begin{cases}
			{\sf Ric}(y)+\left.\frac{\dd}{{\dd}t}\right|_{t=0}{\sf S}_{\sf m}(\dot{\gamma}_y(t))-\frac{{\sf S}_{\sf m}^2(y)}{N-n} & \text{ for }N\in (n, \infty),\\
			\displaystyle\lim_{L\downarrow n} {\sf Ric}_L(y) & \text{ for }N=n, \\
			{\sf Ric}(y)+\frac{\dd}{{\dd}t}\big|_{t=0} {\sf S}_{\sf m}(\dot{\gamma}_y(t)) & \text{ for }N=\infty.
		\end{cases}
		\]
	 When $N=n$, 	${\sf Ric}_n>-\infty$ forces ${\sf S}_{\sf m}=0.$
	 
	 	When $M$ is without boundary,  ${\sf Ric}_N(y) \ge KF^2(y)$ for all $y \in TM$ is equivalent to the synthetic \emph{curvature-dimension condition} ${\sf CD}(K,N)$, see Ohta \cite{Ohta-Calc},  Lott--Villani \cite{LV} and Sturm \cite{Sturm}. When ${\sf Ric}_N\geq 0$ and $\partial M\neq \emptyset$ is locally convex, then   $(M, d_F , {\sf m})$ is a
	 	${\sf CD}(0, N )$ space; in this case, we use the synthetic version of  Bishop--Gromov comparison, etc.   
	 	
	  If ${\sf Ric}_N(y) \ge 0$ (and $\partial M\neq \emptyset$ is locally convex), the Bishop--Gromov comparison principle states that for every $x_0\in M$ and $R>r>0$ we have
	  \begin{equation}\label{BG-inequality}
	  	\frac{{\sf m}(B^\pm_{x_0}(R))}{ R^N}\leq \frac{{\sf m}(B^\pm_{x_0}(r))}{ r^N}.
	  \end{equation}
	  	   The equality case in the latter relation is characterized as follows: 
	  
	  \begin{proposition}\label{N-volume-cone-prop}  Let $\mathcal M=(M,F,{\sf m})$ be a not necessarily reversible $n$-dimensional forward complete Finsler space $($possibly with locally convex boundary$)$ such that ${\sf Ric}_N\geq 0$ for some $N\in [n,\infty)$ and  ${\sf AVR}_{\mathcal M}^+>0$, and  $x_0\in M$ $($possibly $x_0\in\partial M$$).$ Then the following statements are equivalent$:$  
	  		\begin{itemize}
	  			\item[(a)] $\mathcal M=(M,F,{\sf m})$  is a  forward  $N$-volume cone  with  pole $x_0$, i.e., $ {{\sf m}(B^+_{x_0}(R))}={\sf AVR}_{\mathcal M}^+ \omega_N { R^N}$ for every $R>0;$
	  		\item[(b)] $\mathrm{Cut}(x_0)=\emptyset$, that is, every forward geodesic $\gamma$ starting from $x_0$
	  		is defined and minimizing on $[0,\infty)$, and along any of them, the flag curvature vanishes, ${\sf K}(\dot \gamma(t),\cdot)=0,$  while the $S$-curvature verifies ${\sf S}_{\sf m}(\gamma(t),\dot \gamma(t))=\frac{n-N}{t}$  for every $t>0.$
	  			  	\end{itemize}
	  \end{proposition}

%  \begin{proposition}\label{prop:cone}
%  	Let $\mathcal M=(M,F,\mathsf m)$ be a not necessarily reversible $n$-dimensional forward complete
%  	Finsler space, possibly with locally convex boundary, such that $\mathsf{Ric}_N\ge0$ for some
%  	$N\in[n,\infty)$ and $\mathsf{AVR}^+_{\mathcal M}>0$, and let $x_0\in M$ (possibly $x_0\in\partial M$).
%  	Then the following statements are equivalent:
%  	\begin{itemize}
%  		\item[(a)] $\mathsf m(B^+_{x_0}(R))=\mathsf{AVR}^+_{\mathcal M}\,\omega_N R^N$ for every $R>0$;
%  		\item[(b)] $\mathrm{Cut}(x_0)=\emptyset$, i.e. every forward geodesic $\gamma$ emanating from $x_0$
%  		is defined and minimizing on $[0,\infty)$, and along each of them
%  		$\mathsf K(\dot\gamma(t);\cdot)=0$ and $\mathsf S_{\mathsf m}(\gamma(t),\dot\gamma(t))=\frac{n-N}{t}$ for all $t>0$.
%  	\end{itemize}
%  \end{proposition}
	
	 		\begin{proof} (a) $\implies$ (b).  The proof is based on the identification of the equality in the estimates  of the Bishop--Gromov volume comparison \eqref{BG-inequality}; to do this, we recall the crucial steps from its proof, following Ohta \cite[Lemma 7.2 \& Theorem 7.3]{Ohta-Calc} and \cite[Theorem 9.15]{Ohta-book}. In the sequel, we assume that $N>n,$ since the case $N=n$ is classical (when ${\sf S}_{\sf m}=0$, by definition).  
	 			
	 		Let $x_0\in M$,  fix a unit vector $v\in S_{x_0}M,$ and let $\gamma:[0,l)\to M$ be defined by $\gamma(t):=\gamma_v(t)=\exp_{x_0}(tv)$. Assume that $l\leq \rho(v)$, where $\rho(v)=\sup\{t>0:d_F(x_0,\gamma_v(t))=t \}$  is the cut-time in direction $v$. We consider a $g_v$-orthonormal basis $\{e_i\}_{i=1}^{n-1}\cup \{v\}\subset T_{x_0}M$ and the (transversal) Jacobi vector fields $E_i(t)=d(\exp_{x_0})_{tv}(te_i)\in T_{\gamma(t)}M$ along the geodesic $\gamma$ for $i\in \{1,...,n-1\}.$ Let the $(n-1)\times (n-1)$ matrix $A:=A(t)=(a_{ij}(t))$ with $a_{ij}(t)=g_{\dot \gamma(t)}(E_i(t),E_j(t))$. Due to the fact that $l\leq \rho(v)$, for every $t\in (0,l)$ we have another $(n-1)\times (n-1)$ matrix $B:=B(t)=(b_{ij}(t))$ with the property that $E_i'(t )=\sum_{j=1}^{n-1}b_{ij}(t)E_j(t)$; hereafter, $E_i'(t):= D_{\dot \gamma}^{\dot \gamma}E_i(t)$. 
	 		 
	 			 		We consider the function 
	 		$$h(t)=({\rm det}\ A(t))^{1/2(n-1)},\quad t\in [0,l).$$
	 		It is clear that $h(0)=0$, $h(t)>0$ for every $t\in (0,l)$ and after a matrix calculus, it follows that $A'=2BA$ (note that in general, $A'=BA+AB^T$, but $BA$ is symmetric) and
	 		\begin{equation}\label{h-ra-vonatk0zo}
	 			(n-1)h' = \frac{h}{2}{\rm tr}(A'A^{-1}),\quad    (n-1)h''=\frac{h}{n-1}({\rm tr}\, B)^2-h\, {\sf Ric}(\dot \gamma(t))-h\,{\rm tr} (B^2).
	 		\end{equation}

	 		By the algebraic inequality 
	 	\begin{equation}\label{trace-ineq}
	 		({\rm tr}\, B)^2\leq (n-1){\rm tr} (B^2), 
	 	\end{equation}
	 		one has that
	 		\begin{equation}\label{11-trace-ineq-1}
	 			(n-1)h''(t)\leq -h(t)\, {\sf Ric}(\dot \gamma(t)). 
	 		\end{equation}
	 		
	 	We consider the representation of the measure ${\sf m}=e^{-\Psi_\gamma}{\rm vol}_{g_{\dot \gamma}}$ along  $\gamma$ and let 
	 	$$f(t)=e^{-\Psi_\gamma(t)}h(t)^{n-1},\quad h_1(t)=e^{-\Psi_\gamma(t)/(N-n)},\quad h_2(t)=f(t)^{1/(N-1)}.$$	
	 	A formal computation yields that
	 	$$(N-1)\frac{h_2''}{h_2}=(N-n)\frac{h_1''}{h_1}+(n-1)\frac{h ''}{h }-\frac{{(N-n)}{(n-1)}}{N-1}\left(\frac{h_1' }{h_1}-\frac{h ' }{h }\right)^2.$$
	 		According to \eqref{11-trace-ineq-1} and using ${\sf Ric}_N\geq 0$, it follows that
	 	\begin{equation}\label{estimate-ricci-}
	 		(N-1)\frac{h_2''}{h_2}\leq -\Psi_\gamma''(t)+ \frac{(\Psi_\gamma'(t))^2}{N-n}-{\sf Ric}(\dot \gamma(t))=-{\sf Ric}_N(\dot \gamma(t))\leq 0.
	 	\end{equation}
	 		In particular, $(t\, h_2'(t)-h_2(t))'\leq 0$ and $t\, h_2'(t)-h_2(t)\leq \lim_{t\to 0^+}(t\, h_2'(t)-h_2(t))\leq 0$ for every $ t\in [0,l)$. Thus, the function $h_2(t)/t$ is non-increasing on $[0,l)$; by construction, the same holds for $f(t)/t^{N-1}$.   Gromov's elementary monotonicity result for ratios of positive functions (see e.g.\ Chavel \cite[Lemma III.4.1]{Chavel}) implies that $r\mapsto \frac{\int_0^r f(t){\rm d}t}{r^N}$ is non-increasing on $[0,l)$ as well. 
	 		
	 		If $f_v$ is the function from above for a fixed $v\in S_{x_0}M$, then one has for every $r>0$  that 
	 		$${\sf m}(B_{x_0}^+(r))=\int_{S_{x_0}M}\int_0^{\min(r,\rho(v))}f_v(t){\rm d}t A_{x_0}({\rm d}v),$$ where $A_{x_0}$ stands for the measure induced by the Riemannian metric $g_v$.  
	 		Accordingly, for every $R>r>0$, by the above monotonicity, we have that
	 		\begin{align}%\label{ball-measures}
	 		\nonumber {\sf m}(B_{x_0}^+(R))&= \displaystyle\int_{S_{x_0}M}\int_0^{\min(R,\rho(v))}f_v(t){\rm d}t A_{x_0}({\rm d}v)\\&\label{111-1}\leq \displaystyle\int_{S_{x_0}M}\frac{\min(R,\rho(v))^N}{\min(r,\rho(v))^N}\int_0^{\min(r,\rho(v))}f_v(t){\rm d}t A_{x_0}({\rm d}v)
	 		\\&\label{111-2} \leq\frac{R^N}{r^N} \displaystyle\int_{S_{x_0}M}\int_0^{\min(r,\rho(v))}f_v(t){\rm d}tA_{x_0}({\rm d}v)\\&\nonumber=\frac{R^N}{r^N} {\sf m}(B_{x_0}^+(r)).
	 		\end{align}
	 		
	 Now, let us assume that  $\mathcal M=(M,F,{\sf m})$  is a  forward  $N$-volume cone  with  pole $x_0$, i.e., $ {{\sf m}(B^+_{x_0}(R))}={\sf AVR}_{\mathcal M}^+ \omega_N { R^N}$ for every $R>0.$	In particular, it follows that both in \eqref{111-1}		 and  \eqref{111-2} we have equality, i.e., 
	 $\frac{\min(R,\rho(v)) }{\min(r,\rho(v)) }=\frac{R}{r}$ and $\frac{\int_0^{\min(r,\rho(v))} f_v(t){\rm d}t}{r^N}=c_v$ for  every  $R>r>0$  and $A_{x_0}$-a.e.\ $v\in S_{x_0}M,$ where $c_v>0$ depends only on $v$. By the former relation, it follows that $\rho(v)=+\infty$ for $A_{x_0}$-a.e.\ $v\in S_{x_0}M;$ in fact, by the continuity of $\rho$, this relation is valid for every   $v\in S_{x_0}M.$ Therefore, $\mathrm{Cut}(x_0)=\emptyset$. 
	 	 Moreover, the latter relation implies that $f_v(t)/t^{N-1}=Nc_v$ for every $t>0$ and $v\in S_{x_0}M.$ In addition, we should have equality in the estimates \eqref{trace-ineq}, \eqref{11-trace-ineq-1} and \eqref{estimate-ricci-}, respectively. 
	 	 
	 	 First, by the equality in \eqref{trace-ineq} it follows that $B(t)=\lambda(t)I_{n-1}$ for some function $\lambda:(0,\infty)\to \mathbb R;$ in particular, $b_{ij}=\lambda \delta_{ij}$, thus $E_i'=\lambda E_i$ for every $i\in \{1,...,n-1\}.$ This relation implies that $E_i''=(\lambda'+\lambda^2)E_i$, which combined with the Jacobi equation \eqref{Jacobi-eqn}, yields that 
	 	 $ R^{\dot \gamma}(E_i,\dot \gamma)\dot \gamma = -(\lambda'+\lambda^2)E_i$. Consequently, for every  $i\in \{1,...,n-1\} $, by the definition of the flag curvature, we have
	 	 \begin{equation}\label{flag-rigid}
	 	  	 	 	\mathbf{\sf K}(\dot \gamma ;E_i ):=
	 	 	\frac{g_{\dot \gamma } (R^{\dot \gamma }(E_i ,\dot \gamma )\dot \gamma ,E_i )}{g_{\dot \gamma }(\dot \gamma ,\dot \gamma )\,g_{\dot \gamma }(E_i ,E_i ) -g_{\dot \gamma }(\dot \gamma ,E_i )^2}=-(\lambda'+\lambda^2). 
	 	\end{equation}
	  Therefore, 
	  \begin{equation}\label{Ricci-rigid}
	  	{\sf Ric}(\dot \gamma)=-(n-1)(\lambda'+\lambda^2). 
	  \end{equation}
  Since we have equality  in \eqref{estimate-ricci-}, it follows that $\frac{h_1' }{h_1}-\frac{h ' }{h }=0$. On one hand, since $A'=2BA$ and $B =\lambda\, I_{n-1}$, on account of the first relation \eqref{h-ra-vonatk0zo},   we have that 
  $\frac{h ' }{h }=\frac{1}{2(n-1)}{\rm tr}(A'A^{-1})=\lambda.$ On the other hand, by the definition of $h_1$, one has $\frac{h_1' }{h_1}=-\frac{\Psi_\gamma' }{N-n}.$ Therefore, $\frac{\Psi_\gamma' }{N-n}=-\lambda$, and consequently, $\Psi_\gamma''=-(N-n)\lambda'.$ By the equality \eqref{estimate-ricci-} we have ${\sf Ric}_N(\dot \gamma )=0$; replacing the above expression and using \eqref{Ricci-rigid}, it turns out that 
  $$0={\sf Ric}_N(\dot \gamma )= \Psi_\gamma'' - \frac{(\Psi_\gamma' )^2}{N-n}+{\sf Ric}(\dot \gamma )=-(N-1)(\lambda'+\lambda^2),$$
  thus $\lambda'+\lambda^2=0$. In particular, by \eqref{flag-rigid}, it follows that $$\mathbf{\sf K}(\dot \gamma(t) ;\cdot )=0,\quad \forall t>0.$$  
  
  The latter fact together with \eqref{Ricci-rigid} implies that ${\sf Ric}(\dot \gamma)=0$; since we have equality in \eqref{11-trace-ineq-1}, it follows that $h ''=0.$ Since $h(0)=0$, it turns out there exists $c>0$ such that $h(t)=c\,t$ for every $t>0$. Therefore, $\lambda(t)=\frac{h ' (t)}{h(t) }=\frac{1}{t}$ for every $t>0$. Consequently, we have that 
  $$ {\sf S}_{\sf m}(\gamma(t),\dot \gamma(t))=\Psi_\gamma'(t)=-(N-n)\lambda(t)=\frac{n-N}{t},\quad \forall t>0.  $$
  
  (b) $\implies$ (a). Let  
   $v\in   S_{x_0}M$.   
  Since $\mathsf K(\dot\gamma_v;\cdot)=0$, where $\gamma_v(0)=x_0$ and $\dot \gamma_v(0)=v$,  it follows that $ R^{\dot \gamma_v} = 0$ on the orthogonal of $\dot \gamma_v$ and it results that every Jacobi field $E_i$ along $\gamma_v$ has the
  form   $E_i(t) = t P_{v,i}(t)$, where $P_{v,i}$ is a parallel vector field along $ \gamma_v$, i.e., $D_{\dot \gamma_v}^{\dot \gamma_v}P_{v,i}=0$ and $P_{v,i}(0)=e_i$. Thus, using the above notations, one has that  
  $A(t)=t^2I_{n-1}$
  and $h(t)=t$, $t>0$.  Since $\Psi_{\gamma_v}'(t)=\mathsf S_{\mathsf m}(\gamma_v(t),\dot\gamma_v(t))=\frac{n-N}{t}$, it yields that
  $e^{-\Psi_{\gamma_v}(t)}=C_v t^{N-n}$ for some $C_v>0$, i.e., $f_v(t)=C_v t^{N-1}$. As $\rho(v)=\infty$ (note that   $\mathrm{Cut}(x_0)=\emptyset$), it follows for every $R>0$ that 
  $$\mathsf m(B^+_{x_0}(R))=\displaystyle\int_{S_{x_0}M}\int_0^{R}f_v(t){\rm d}t A_{x_0}({\rm d}v)=\frac{R^N}{N}\int_{S_{x_0}M}C_v\,\mathrm dA_{x_0}(v).$$  Therefore, 
  $$\int_{S_{x_0}M}C_v\,\mathrm dA_{x_0}(v)=N \lim_{R\to \infty}\frac{\mathsf m(B^+_{x_0}(R))}{R^N}=N\omega_N \mathsf{AVR}^+_{\mathcal M}.$$
  Accordingly, we have $\mathsf m(B^+_{x_0}(R))=\mathsf{AVR}^+_{\mathcal M}\,\omega_N R^N$
for every $R>0,$ 
i.e., $\mathcal M=(M,F,{\sf m})$  is a  forward  $N$-volume cone  with  pole $x_0$.
  	 		\end{proof}

		\subsection{Optimal mass transport on Finsler manifolds}
		
\subsubsection{Brenier--McCann--Ohta theory}  Let $\mathcal M=(M,F,{\sf m})$ be a not necessarily reversible $n$-dimensional  complete Finsler space $($possibly with boundary$)$ such that ${\sf Ric}_N\geq 0$ for some $N\in [n,\infty)$. For further use, $c$ always stands for the quadratic cost $c(x,y)=d_F^2(x,y)/2$. 

Let $X,Y\subset M$ be two compact sets.  A function $\phi:X\to \mathbb R\cup \{-\infty\}$ is said to be $c$-\textit{concave} relative to $(X,Y)$ if it is not identically $-\infty$ and there exists $\psi:Y\to \mathbb R\cup \{-\infty\}$ such that 
$$\phi(x)=\inf_{y\in Y}(c(x,y)-\psi(y)).$$

  Let 	$\mathcal P(M)$ be the set of Borel probability measures on $M$ and  	$\mathcal P_c(M)\subset \mathcal P(M)$ be the subset of compactly supported measures. If $\mu_0,\mu_1\in \mathcal P_c(M)$, Monge's problem is to find a map $\mathcal F:M\to M$ which realizes the infimum of the expression $\int_M c(x,\mathcal F(x)){\rm d}\mu_0(x)$ among the maps pushing forward $\mu_0$ to $\mu_1$ (i.e., $\mathcal F_\#\mu_0=\mu_1$). 
  
  In the Euclidean setting, Brenier \cite{Brenier} guaranteed the existence and uniqueness of the solution to Monge's problem, which is the gradient of a convex function. Brenier's result has been extended by McCann \cite{McCann} to Riemannian manifolds, which in turn, has its Finslerian counterpart, due to Ohta \cite{Ohta-Calc}. Namely, it turns out that, if in addition,    $\mu_0$ is absolutely continuous with respect to the   reference measure ${\sf m}$, and taking a compact set $Y\supset {\rm supp}\, \mu_1$ and an open set $U\supset {\rm supp}\, \mu_0$ with $X=\overline{U}$ compact, then there exists a $c$-concave function $\phi:X\to \mathbb R\cup \{-\infty\}$ relative to $(X,Y)$ such that $\mathcal F(x)=\exp_x(\bm \nabla(- \phi)(x))$ solves Monge's problem. Moreover, if $\mu_0=\rho_0 {\sf m}$ and $\mu_1=\rho_1 {\sf m}\in \mathcal P_c(M)$ are absolutely continuous measures with respect to the measure ${\sf m}$, it turns out that $\mathcal F$ is injective on a $\mu_0$-full measure set and  \textit{Monge--Amp\`ere equation} holds $\mu_0$-a.e. $x\in U\setminus (D\phi)^{-1}(0)$: 
  \begin{equation}
  	\label{MA-initial}
  \rho_0(x)=\rho_1(\mathcal F (x))\cdot {\bf D}[d \mathcal F  _x],
 \end{equation}
    see Ohta \cite[Theorem 5.2]{Ohta-Calc}.
  	Here,  ${\bf D}[d \mathcal F _x]$ stands for the Jacobian of the linear operator $d \mathcal F _x:T_xM\to T_{\mathcal F (x)}M$ with respect to the measure ${\sf m}.$ More precisely, if $L=\{x\in M:D\phi(x)=0\}$, then one has that ${\bf D}[d \mathcal F _x]=1$ for every $x\in L$,
  %	 (in which case  $\mathcal F_1(x)=\mathcal F_t(x)=x$  and 
 % $\rho_0|_L{\sf m}=\rho_t|_L{\sf m}=\rho_1|_L{\sf m}>0$), 
  while 
  $${\bf D}[d \mathcal F _x]=\frac{{\sf m}_{\mathcal F (x)}(d \mathcal F _x(A))}{{\sf m}_x(A)}$$ for $x\in M\setminus L$, where $A\subset T_xM$ is any open bounded set  (with ${\sf m}_x$ and ${\sf m}_{\mathcal F  (x)}$ as the Lebesgue measures on $T_xM$ and $T_{\mathcal F (x)}M$, respectively, induced by ${\sf m}$). 
		
		\subsubsection{Functional inequality from displacement convexity}
		
		For every $N\in [1,\infty)$ we introduce the set ${\mathcal DC}_N$  of continuous convex functions $U:[0,\infty)\to \mathbb R$ such that $U(0)=0$ and  the function $s\mapsto s^NU(s^{-N})$ is convex on $(0,\infty)$; in a similar way, let ${\mathcal DC}_\infty$  be the set of continuous convex functions $U:[0,\infty)\to \mathbb R$ such that $U(0)=0$ and  the function $s\mapsto e^sU(e^{-s})$ is convex on $\mathbb R$. Moreover, for an absolutely continuous measure $\mu=\rho{\sf m}\in \mathcal P(M)$, we define the \textit{internal energy} functional $$\mathcal U_{\sf m}(\mu)=\int_M U(\rho){\rm d}{\sf m}.$$
		
		\begin{proposition}\label{displacement-convexity}
				Let $\mathcal M=(M,F,{\sf m})$ be an $n$-dimensional forward-complete Finsler space $($possibly with locally convex boundary$)$ such that ${\sf Ric}_N\geq 0$ for some $N\in [n,\infty)$. If $\mu_0=\rho_0 {\sf m}$,  $\mu_1=\rho_1 {\sf m}\in \mathcal P_c(M)$, and $\mathcal F_t(x)=\exp_x(\bm \nabla(-t\phi)(x))$, $t\in [0,1],$ is the unique optimal transport map between $\mu_0$ and $\mu_1$ for some $c$-concave function	$\phi:M\to \mathbb R$, then 
				\begin{itemize}
					\item[(i)] for ${\sf m}$-a.e.\ $x\in {\rm supp}\mu_0$, one has  
					$({\bf D}[d (\mathcal F_1) _x])^{1/N}\leq 1+\frac{1}{N}{\bm\Delta_A(-\phi)(x)};$
					\item[(ii)]
					$\mathcal U_{\sf m}(\mu_1)-\mathcal U_{\sf m}(\mu_0)\geq %\begin{cases}
				\displaystyle\int_M \left(U(\rho_0)-\rho_0U'(\rho_0)\right)\bm\Delta_A(-\phi) {\rm d}{\sf m} 	 
				$  for every $U\in {\mathcal DC}_N;$
					\item[(iii)]    $\mathcal U_{\sf m}(\mu_1)-\mathcal U_{\sf m}(\mu_0)\geq %\begin{cases}
				-N\log	\left(\displaystyle\int_M   \rho_0 ({\bf D}[d  \mathcal F_1  ])^{1/N} {\rm d}{\sf m}\right) 	 
					$ for $U(s)=s\log s$, $s\geq 0.$
					
						\end{itemize}
		\end{proposition}
		
		\begin{proof}  Let $\mu_t=(\mathcal F_t)_\#\mu_0$, $t\in [0,1],$ be the unique minimal geodesic from $\mu_0$ to $\mu_1$. We notice that $\mathcal F_t$ is injective on a set of $\mu_0$-full measure, and   $\mu_t$ is absolutely continuous with respect to ${\sf m};$ let  $\rho_t$ be its  density, i.e., $\mu_t=\rho_t {\sf m}.$
		
		Let $L=\{x\in M:d\phi_x=0\}$. The map $\mathcal F_t$ is of class $C^\infty$ on  $M\setminus L$, while  the Monge--Amp\`ere equation arising from  $\mu_t=(\mathcal F_t)_\#\mu_0$, $t\in [0,1],$ reads as 
		\begin{equation}\label{Monge-Ampere-t}
			\rho_0(x)=\rho_t(\mathcal F_t(x))\cdot {\bf D}[d(\mathcal F_t)_x]\ \ \ \ {\rm for}\ \mu_0-{\rm a.e.}\ x\in M;
		\end{equation}	 
	 note that ${\bf D}[d (\mathcal F_t) _x]=1$ for every $x\in L$ 
	 	 (in which case  $\mathcal F_1(x)=\mathcal F_t(x)=x$  and 
	  $\rho_0|_L{\sf m}=\rho_t|_L{\sf m}=\rho_1|_L{\sf m}>0$). 
	
		 (i) For $\mu_0$-a.e.\ $x\in M$ -- in fact, on points of twice differentiability -- one has that ${\bf D}[d(\mathcal F_0)_x]=1$, while due to Ohta--Sturm \cite[p.\ 1424]{Ohta-Sturm} or Ohta \cite[Lemma 12.6/(ii)]{Ohta-book}, it turns out that 
		 \begin{equation}\label{derivative}
		 	\frac{{\rm d}}{{\rm d}t}\Big|_{t=0}{\bf D}[d(\mathcal F_t)_x]=\bm\Delta_A(-\phi)(x).
		 \end{equation}
	 	 Fix a point  $x\in M$ from a $\mu_0$-full measure set where $\mathcal F_1$ is injective.   Since $t\mapsto J(t):={\bf D}[d(\mathcal F_t)_x]^{1/N}$ is concave on $[0,1]$, see Ohta \cite[rel.\ (8.6)]{Ohta-Calc}, we have that $J(1)\leq J(0)+J'(0)$. Due to  ${\bf D}[d(\mathcal F_0)_x]=1$ and \eqref{derivative}, the latter inequality completes the proof. 
		
		(ii)	Let $U\in {\mathcal DC}_N $ be fixed. According  to Ohta \cite[Theorem 1.2 \& Definition 8.1]{Ohta-Calc}, since  ${\sf Ric}_N\geq 0$, the internal energy $\mathcal U_{\sf m}$ is displacement convex on $\mathcal P_c(M)$; in particular, $$t\, \mathcal U_{\sf m}(\mu_1)+(1-t)\,\mathcal U_{\sf m}(\mu_0)\geq \mathcal U_{\sf m}(\mu_t),  \ \ t\in [0,1].$$ 
			 Therefore, we have that
			\begin{equation}\label{U-convex-deriv}
				\mathcal  U_{\sf m}(\mu_1)-\mathcal U_{\sf m}(\mu_0) \geq \frac{\rm d}{{\rm d}t}\Big|_{t=0}\mathcal U_{\sf m}(\mu_t).
			\end{equation}
				A change of variables  and the Monge--Amp\`ere equation \eqref{Monge-Ampere-t} imply that
		\begin{eqnarray*}
	\mathcal 	U_{\sf m}(\mu_t)&=&\int_M U(\rho_t){\rm d}{\sf m}=\int_M U(\rho_t(\mathcal F_t(x))){\bf D}[d(\mathcal F_t)_x]{\rm d}{\sf m}(x)\\&=&\int_M U\left(\rho_0(x){\bf D}[d(\mathcal F_t)_x]^{-1}\right){\bf D}[d(\mathcal F_t)_x]{\rm d}{\sf m}(x).
		\end{eqnarray*}
		The latter relation  combined with  \eqref{derivative} and ${\bf D}[d(\mathcal F_0)_x]=1$ provide  that
			\begin{eqnarray*}
		\frac{\rm d}{{\rm d}t}\Big|_{t=0}\mathcal U_{\sf m}(\mu_t)&=&\int_M \frac{\rm d}{{\rm d}t}\Big|_{t=0} U\left(\rho_0(x){\bf D}[d(\mathcal F_t)_x]^{-1}\right){\bf D}[d(\mathcal F_t)_x]{\rm d}{\sf m}(x)\\&=&\int_M \left(U(\rho_0)-\rho_0U'(\rho_0)\right)\bm\Delta_A(-\phi) {\rm d}{\sf m}.
			\end{eqnarray*}
		It remains to combine the latter relation with  \eqref{U-convex-deriv}. 	
	
	(iii) Since $\mu_0=\rho_0 {\sf m}\in \mathcal P_c(M)$,   by Jensen's inequality   it follows that
	$$N\log	\left(\displaystyle\int_M   \rho_0 ({\bf D}[d  \mathcal F_1  ])^{1/N} {\rm d}{\sf m}\right)\geq \displaystyle\int_M   \rho_0 \log ({\bf D}[d  \mathcal F_1  ])  {\rm d}{\sf m}.$$ Moreover, the Monge--Amp\`ere equation  \eqref{Monge-Ampere-t} and a change of variable implies that 
\begin{align*}
\displaystyle\int_M   \rho_0 \log ({\bf D}[d  \mathcal F_1  ])  {\rm d}{\sf m}&= \int_M   \rho_0 \log \left( \frac{\rho_0}{\rho_1 \circ \mathcal F_1  }\right)  {\rm d}{\sf m}=\int_M   \rho_0 \log \rho_0    {\rm d}{\sf m}-\int_M   \rho_0  \log \rho_1 \circ \mathcal F_1     {\rm d}{\sf m}\\&= \int_M   \rho_0 \log \rho_0    {\rm d}{\sf m}-\int_M   \rho_1  \log \rho_1      {\rm d}{\sf m},
\end{align*} 
 concluding the proof. 
		\end{proof}
	
	\begin{remark}\rm 
		The Euclidean form of Proposition \ref{displacement-convexity}/(ii)
	is due to Agueh--Ghoussoub--Kang \cite{AGK}; see also Villani \cite[Chapter 5]{Villani-book}. 
\end{remark}
		
 	\subsection{Rearrangements}\label{subsection-rearr} Let $\mathcal M=(M,F,{\sf m})$ be an $n$-dimensional  complete Finsler space $($possibly with boundary$)$ such that ${\sf Ric}_N\geq 0$ for some $N\in [n,\infty)$ and ${\sf AVR} _{\mathcal M}>0$. For a measurable function $u\colon M\to[0,\infty)$, let $u^*\colon[0,\infty)\to[0,\infty)$ be the \textit{non-increasing rearrangement} of $u$ defined on the 1-dimensional model space $(\mathbb R_+,\vert\cdot\vert,{\sf m}_N=N\omega_Nr^{N-1}{\mathcal L ^1}_{\llcorner{[0,\infty)}})$ so that
 \begin{align}\label{rearrangement}
 	{\sf m}\big(\mathcal N_t(u)\big)={\sf m}_N(\mathcal N_t(u^*)),\ \ \forall t>0,
 \end{align}
 where 
 \begin{align*}
 	\mathcal N_t(u)=\lbrace x\in M: u(x)>t\rbrace\;\text{ and }\;\mathcal N_t(u^*)=\lbrace y\in[0,\infty):u^*(y)>t\rbrace,
 \end{align*}
 whenever ${\sf m}\big(\mathcal N_t(u)\big)<\infty$ for every $t>0$. 
 In particular, the open set $\Omega\subset M$ can be rearranged into the interval $\Omega^*=[0,r)$ with $ {\sf m}(\Omega)={\sf m}_N([0,r))=\omega_N r^N,$ with the convention that $\Omega^*=[0,\infty)$ whenever $ {\sf m}(\Omega)=+\infty$. For every continuous function $G\colon[0,\infty)\to \mathbb R$ and $u\colon M\to[0,\infty)$, the layer cake representation implies  the \textit{Cavalieri principle}: 
 \begin{equation}\label{Cavallieri}
 	\int_\Omega	G(u){\rm d}{\sf m}= \int_{\Omega^*}	G(u^*){\rm d}{\sf m}_N,
 \end{equation}
 assuming that the integrals are well-defined. 

By using the optimal  isoperimetric inequality \eqref{Manini-ineq} and coarea formula, a slight modification of the argument from Nobili--Violo \cite{NV-Calc} (see also  Mondino--Semola \cite{Mondino-Semola}) implies the following  \textit{forward P\'olya--Szeg\H o inequality}:  for every $p>1$ and  non-negative $u\in W^{1,p}_{\rm loc} (\mathcal M)$, one has
 \begin{equation}\label{Polya-Szego}
 	  	\|F^*(-Du)\|_{L^p(M,{\sf m})}\geq   ( {\sf AVR}^+_{\mathcal M})^{1/N}\|(u^*)'\|_{L^p(\mathbb R_+,{\sf m}_N)}.
 \end{equation}
If the left hand side is finite then $u^*$ is locally absolutely continuous on its domain of definition. Clearly, similar version to \eqref{Polya-Szego} can be stated for the reverse metric $\overleftarrow{F}$ with ${\sf AVR}^-_{\mathcal M}$, which is the   \textit{backward P\'olya--Szeg\H o inequality}.  The Riemannian analogue of    \eqref{Polya-Szego}  has been established by Balogh--Krist\'aly  \cite{BK}. 

 Assume in addition that $r_F<\infty$ and $\partial M$ is locally   convex  if $M$ has boundary. A closer inspection of the proof of Nobili--Violo \cite[Theorem 1.3]{NV-Calc}   shows that if equality holds in \eqref{Polya-Szego} for some non-negative $u\in W^p_p (\mathcal M)\setminus \{0\}$, 
% 
% The description of the equality case in the P\'olya--Szeg\H o inequality \eqref{Polya-Szego} is more delicate, even in the smooth setting. By using fine analysis of sets of finite perimeter,  Brothers and Ziemer \cite{BZ} characterized the equality in \eqref{Polya-Szego} in the Euclidean setting. Under certain regularity assumptions, a similar result to \cite{BZ} is provided  in \cite{BK} in the setting of Riemannian manifolds with non-negative Ricci curvature. As far as we know -- by using a general rearrangement technique  and  the sharp isoperimetric inequality \eqref{eqn-isoperimetric} -- the most general setting is handled by Nobili and Violo \cite[Theorem 1.3]{NV-new}: 
% if $(X,{\sf d}, {\sf m})$ is an essentially non-branching ${\sf CD}(0,N)$ space and equality holds for some  $u\in W^{ 1,p}_{\rm loc}(X,{\sf m})$ and 
% both sides being non-zero and finite in \eqref{Polya-Szego}, 
 then there exists  $x_0\in M$ such that  $\mathcal M=(M,F,{\sf m})$  is a \textit{forward} $N$\textit{-volume cone} with  pole $x_0$, i.e., 
 \begin{equation}\label{volume-cone}
 	\frac{ {\sf m}(B_{x_0}^+(R))}{\omega_NR^N}={\sf AVR}_{\mathcal M}^+,\quad \forall R>0.
 \end{equation}
In fact, \eqref{volume-cone} follows by the equality case in the optimal  isoperimetric inequality \eqref{Manini-ineq}, see Manini \cite[Theorem 1.4]{Manini}.

 In addition, if $(u^*)'\neq 0$ a.e.\ on   $\{{\rm ess\, inf}\, u<u^*<{\rm ess\, sup}\, u\}$ then $u$ is $x_0$-\textit{radial}, i.e., 
 \begin{equation}\label{u=u^*}
 	u(x)=u^*\left({\sf AVR}^+_{\mathcal M})^{1/N} d_F(x_0,x)\right)\ \ {\sf m}{\rm -a.e.}\ x\in M.
 \end{equation}

			\section{Optimal Finsler--Sobolev inequality: proof of Theorems \ref{Thm-Sobolev-1} \& \ref{Thm-Sobolev-2-equality}}    \label{section-3}
			
			\subsection{The case $q\neq p_N^*$} Let us consider a smooth non-increasing cut-off function $\eta:[0,\infty)\to [0,1]$ such that $\eta(t)=1$ for $t\in [0,1]$,  $\eta(t)=0$ for $t\geq 2 ,$ and $|\eta'|\leq C$ for some $C>0$ large enough. 
			
			For $x_0\in M$ and $R>0$ fixed, let 
			$$u_R(x)=-\eta\left(\frac{d_F(x_0,x)}{R}\right),\ x\in M.$$
			It is clear that $u_R\equiv -1$ on $ B^+_{x_0}(R)$, $u_R\equiv 0$ on $M\setminus  B^+_{x_0}(2R)$, and according to the eikonal equation \eqref{eikonal}, $F^*(Du_R)\leq C/R$. Therefore, one has that
			$${\sf m}(B^+_{x_0}(R))\leq \|u_R\|_{L^q(M,{\sf m})}^q=\int_M |u_R|^q {\rm d}{\sf m} \leq {\sf m}(B^+_{x_0}(2R)),$$
			and 
		$$\|F^*(Du_R)\|_{L^p(M,{\sf m})}^p=\int_M F^*(Du_R)^p {\rm d}{\sf m}\leq \frac{C^p}{R^p}{\sf m}(B^+_{x_0}(2R)).$$
		Consequently, $u_R\in W_q^p(\mathcal M)$ and  by the Bishop--Gromov comparison principle, one has
		\begin{equation}\label{estimate-R}
				\mathcal R_q^p(\mathcal M)\leq \frac{\|F^*(Du_R)\|_{L^p(M,{\sf m})}}{\|u_R\|_{L^q(M,{\sf m})}}\leq C \frac{{\sf m}(B^+_{x_0}(2R))^\frac{1}{p}}{R\, {\sf m}(B^+_{x_0}( R))^\frac{1}{q}} \leq 2^\frac{N}{p}C \frac{{\sf m}(B^+_{x_0}( R))^{\frac{1}{p}-\frac{1}{q}}}{R}
		\end{equation}
				for every $R>0$. 
				
			\subsubsection{The case $q<p_N^*$}	 Since ${\sf AVR}_{\mathcal M}^+>0$ (recall that   ${\sf AVR}_{\mathcal M} >0$ by assumption), relation    \eqref{estimate-R} implies through $q<p_N^*=\frac{pN}{N-p}$ that
				\begin{align*}
				\mathcal R_q^p(\mathcal M)&\leq 2^\frac{N}{p}C \lim_{R\to \infty}\frac{{\sf m}(B^+_{x_0}( R))^{\frac{1}{p}-\frac{1}{q}}}{R}=2^\frac{N}{p}C\omega_N^{\frac{1}{p}-\frac{1}{q}}\lim_{R\to \infty}\left(\frac{{\sf m}(B^+_{x_0}( R))}{\omega_N R^N}\right)^{\frac{1}{p}-\frac{1}{q}}R^{\frac{N}{p}-\frac{N}{q}-1}\\&=2^\frac{N}{p}C\omega_N^{\frac{1}{p}-\frac{1}{q}} \left({\sf AVR}^+_{\mathcal M}\right)^{\frac{1}{p}-\frac{1}{q}}\lim_{R\to \infty} R^{\frac{N}{p}-\frac{N}{q}-1}=0.
				\end{align*}

		\subsubsection{The case $q>p_N^*$ and $N$ saturates}	Since $N$ saturates, there exists $x_0\in M$ $($or $x_0\in \partial M$ if $M$ has boundary$)$ such that   $\theta_N(x_0)<\infty$. Without loss of generality, we may assume that $\theta_N^+(x_0)<\infty$; otherwise, the above estimates are considered for the reverse Finsler metric $\overleftarrow{F}$  and backward balls $B^-_{x_0}(R)$. Thus, again using the estimate \eqref{estimate-R}, it turns out that 
				\begin{align*}
				\mathcal R_q^p(\mathcal M)&\leq 2^\frac{N}{p}C \lim_{R\to 0}\frac{{\sf m}(B^+_{x_0}( R))^{\frac{1}{p}-\frac{1}{q}}}{R}=2^\frac{N}{p}C\lim_{R\to 0}\left(\frac{{\sf m}(B^+_{x_0}( R))}{  R^N}\right)^{\frac{1}{p}-\frac{1}{q}}R^{\frac{N}{p}-\frac{N}{q}-1}\\&=2^\frac{N}{p}C  \left(\theta_N^+(x_0)\right)^{\frac{1}{p}-\frac{1}{q}}\lim_{R\to 0} R^{\frac{N}{p}-\frac{N}{q}-1}=0,
			\end{align*}
				which concludes the proof. 
	 			
			\subsection{The case $q= p_N^*$} 
			
			\subsubsection{Proof of $\mathcal R^{p}_{p_N^*}(\mathcal M)\geq {\sf AT}(p,N)\cdot   {\sf AVR} _{\mathcal M} ^\frac{ 1}{N}$}\label{ssss-us} This part is divided into two steps. 
			
		\textit{Step 1.} 	We   prove first that
			\begin{equation}\label{Sob-pozitiv}
				\|F^*(-Du)\|_{L^p(M,{\sf m})} \geq \left({\sf AVR}^+_{\mathcal M}\right)^\frac{1}{N}\cdot {\sf AT}(p,N)\, \|u\|_{L^{p_N^*}(M,{\sf m})} 
			\end{equation}
		for every \textit{non-negative} $u$ with  $-u\in W^{p}_{p_N^*}(\mathcal M)$ (the latter being equivalent to $ u\in W^{p}_{p_N^*}(\overleftarrow {\mathcal M})$). By density reasons, it is enough to consider non-negative, compactly supported   function $u\in C_0^\infty(M)$,   normalized as   
		 		\begin{equation}\label{u-normalization}
		 			\|u\|_{L^{p_N^*}(M,{\sf m})} =1.
		 		\end{equation}
		For further use, let $\Omega=\{x\in M:u(x)>0\}$; in particular, $\overline \Omega\subset \overline M$ is compact. 
		
	Through the proof, we use certain Talentian-type functions on $M,$ appearing in the target density in our optimal mass transport argument and plays a crucial role in the appearance of the asymptotic volume ratio.  	Fixing $s>\frac{N}{p'},$ $\lambda>0$ and $x_0\in M$    arbitrarily, we first observe  that $x\mapsto (\lambda+ d_F^{p'}(x_0,x))^{-s}$, $x\in M$, belongs to $L^1(M,{\sf m})$. Indeed, by the layer cake representation and Bishop--Gromov comparison principle, it turns out that      
		\begin{align}\label{layer-cake}
		\nonumber	I_{\lambda,s}^{x_0}&:=\int_M (\lambda+ d_F^{p'}(x_0,x))^{-s}{\rm d}{\sf m}(x)= \int_0^{\lambda^{-s}}{\sf m}\left(\left\{x\in M:(\lambda+ d_F^{p'}(x_0,x))^{-s}>t\right\}\right){\rm d}t\\&=sp'\int_0^\infty {\sf m}(B^+_{x_0}(\rho))\left(\lambda+\rho^{p'}\right)^{-s-1}\rho^{p'-1}{\rm d}\rho\\&\leq sp'{\sf m}(B^+_{x_0}(1))\left(\int_0^1 \left(\lambda+\rho^{p'}\right)^{-s-1}\rho^{p'-1}{\rm d}\rho+\int_1^\infty \left(\lambda+\rho^{p'}\right)^{-s-1}\rho^{N+p'-1}{\rm d}\rho\right)\nonumber\\&<+\infty.\nonumber
		\end{align}
Moreover, by a change of   variables $\rho =\lambda^{\frac{1}{p'}}r $ in \eqref{layer-cake} together with the monotone convergence theorem (applied to $t\mapsto {\sf m}(B^+_{x_0}(t))/t^N$ via the Bishop-Gromov comparison principle) one has that
\begin{align}\label{AVR-in-Limit}
\nonumber\lim_{\lambda\to \infty}\lambda^{s-\frac{N}{p'}}\, I_{\lambda,s}^{x_0}&=sp'\omega_N \lim_{\lambda\to \infty}\int_0^\infty
\frac{{\sf m}(B^+_{x_0}(\lambda^{\frac{1}{p'}}r))}{\omega_N \left(\lambda^{\frac{1}{p'}}r\right)^N}\left(1+r^{p'}\right)^{-s-1}r^{N+p'-1}{\rm d}\rho\\&\nonumber=sp'\omega_N {\sf AVR}^+_{\mathcal M} \int_0^\infty
\left(1+r^{p'}\right)^{-s-1}r^{N+p'-1}{\rm d}\rho
\\& =s\,\omega_N \,  {\sf AVR}^+_{\mathcal M} \, {\sf B}\left(\frac{N}{p'}+1,s-\frac{N}{p'}\right).
\end{align}

We now return  to the proof of \eqref{Sob-pozitiv}. For every $k\in \mathbb N$, let $\eta_k:M\to [0,1]$ be a  continuous function such that $\eta_k\equiv 1$ on $B_{x_0}^+(k)$ and $\eta_k\equiv 0$ on $M\setminus B_{x_0}^+(k+1)$; for instance, $\eta_k(x)=(\min(0,k-d_F(x_0,x))+1)_+$. Note that ${\rm supp}\, \eta_k\subseteq \overline B_{x_0}^+(k+1)$. 

Let us consider the measures  $\mu_0=\rho_0 {\sf m}$ and $\mu_1=\rho_1 {\sf m}$ with the densities 
$$ \rho_0(x)=u^{p_N^*}(x) \ \ \text{and}\ \ \rho_1(x)=(I_{\lambda,N,k}^{x_0})^{-1} \eta_k(x) (\lambda+ d_F^{p'}(x_0,x))^{-N},$$
where $$I_{\lambda,N,k}^{x_0}=\int_M \eta_k(x) (\lambda+ d_F^{p'}(x_0,x))^{-N}{\rm d}{\sf m}(x).$$
Relation  \eqref{u-normalization} and  the well-definiteness of $I_{\lambda,N}^{x_0}$  together with $0<I_{\lambda,N,k}^{x_0}\leq I_{\lambda,N}^{x_0}$ for every $k\in \mathbb N$, imply  that the measures $\mu_0$ and $\mu_1$ belong to $\mathcal P_c(M)$. 
				
Let	$U(s)=-s^{1-\frac{1}{N}},$ $s\geq 0;$ it is clear that $s\mapsto s^NU(s^{-N})=-s$ is convex, thus $U\in {\mathcal DC}_N$. Therefore, we may apply 	Proposition \ref{displacement-convexity}/(ii) for these choices of $U$,  $\mu_0$ and $\mu_1$, obtaining that 
\begin{equation}\label{disp-convexity-Sobolev}
	-\int_M \rho_1^{1-\frac{1}{N}}{\rm d}{\sf m}+\int_M \rho_0^{1-\frac{1}{N}}{\rm d}{\sf m}\geq %\begin{cases}
	-\frac{1}{N}\displaystyle\int_M \rho_0^{1-\frac{1}{N}}\bm\Delta_A(-\phi) {\rm d}{\sf m},
\end{equation}
 where $\phi:M\to \mathbb R$ is a  $c$-concave function with the property that $\mathcal F_t(x)=\exp_x(\bm \nabla(-t\phi)(x))$, $t\in [0,1],$ is the unique optimal transport map between the probability measures $\mu_0$ and $\mu_1$. We stress that $\rho_1$, $\phi$ and $\mathcal F_t$ all depend on $k$ and $\lambda$; to simplify the notation, 
 this dependence will be made explicit only when necessary.
 
 Since $\phi$ is a $c$-concave function, it turns out that $-\phi$   is semiconvex, thus 
 \begin{equation}\label{distributional-alexandrv}
 	 \bm\Delta_A(-\phi)\leq \bm\Delta_{\mathcal M}(-\phi)
 \end{equation}
  in the distributional sense. 
  Moreover, we also have the Monge--Amp\`ere equation \eqref{MA-initial}, i.e., 
 \begin{equation}\label{MA-Sob}
 	\rho_0(x)=\rho_1(\mathcal F_1(x))\cdot {\bf D}[d(\mathcal F_1)_x]\ \ \ \ {\rm for}\ \mu_0-{\rm a.e.}\ x\in \Omega.
 \end{equation}
 % Since $u\in C_0^\infty(M)$, we have that $$C_0:=\int_M \rho_0^{1-\frac{1}{N}}{\rm d}{\sf m}=\int_\Omega u^{(1-\frac{1}{N})p_N^*}{\rm d}{\sf m}$$
% is finite and does not depend on $\lambda>0$. 
  Reorganizing \eqref{disp-convexity-Sobolev} and applying \eqref{distributional-alexandrv} and  \eqref{Laplacian} together with \eqref{Cauchy-Schwarz}, it follows that 
	\begin{align}\label{rho-1-rho-2}
		\nonumber	\int_M \rho_1^{1-\frac{1}{N}}{\rm d}{\sf m}&\leq \int_M \rho_0^{1-\frac{1}{N}}{\rm d}{\sf m}+\frac{1}{N}\int_M u^{(1-\frac{1}{N})p_N^*}\bm\Delta_{\mathcal M}(-\phi) {\rm d}{\sf m}\\&	\nonumber = \int_\Omega \rho_0^{1-\frac{1}{N}}{\rm d}{\sf m}- \frac{1}{N}\left(1-\frac{1}{N}\right)p_N^*\int_\Omega u^{(1-\frac{1}{N})p_N^*-1}   Du (\bm \nabla (-\phi))  {\rm d}{\sf m}\\&\leq \int_\Omega \rho_0^{1-\frac{1}{N}}{\rm d}{\sf m}+ \frac{p_N^*(N-1)}{N^2}\int_\Omega u^{(1-\frac{1}{N})p_N^*-1}  F^*(- Du) F(\bm \nabla (-\phi))  {\rm d}{\sf m}.
%			\\&\leq C_0  + \frac{p_N^*(N-1)}{N^2} \|F^*(-Du)\|_{L^p(M,{\sf m})}\left(\int_M u^{\left( (1-\frac{1}{N})p_N^*-1\right)p'}F(\bm \nabla (-\phi))^{p'}{\rm d}{\sf m}\right)^\frac{1}{p'}.
	\end{align}
By the definition of the optimal transport map $\mathcal F_1$, one has    for $\mu_0$-a.e.\ $x\in \Omega$ that 
\begin{align}\label{F-finsler-estimate-distance}
\nonumber F(\bm \nabla (-\phi)(x))&=d_F(x,\exp_x(\bm \nabla(- \phi)(x)))\leq d_F(x,x_0)+d_F(x_0,\exp_x(\bm \nabla(- \phi)(x)))\\&\leq K+d_F(x_0,\mathcal F_1(x)), 
\end{align}
where $K=\max_{x\in \overline\Omega}d_F(x,x_0).$ If $$C_0:=\int_\Omega \rho_0^{1-\frac{1}{N}}{\rm d}{\sf m}+\frac{p_N^*(N-1)}{N^2}K\int_\Omega u^{(1-\frac{1}{N})p_N^*-1}     {\rm d}{\sf m},$$
by the estimates \eqref{rho-1-rho-2}  and \eqref{F-finsler-estimate-distance} 
%together with H\"older's inequality, 
one obtains that
\begin{align}\label{majdnemutolso}
	\int_M \rho_1^{1-\frac{1}{N}}{\rm d}{\sf m}&\leq C_0+ \frac{p_N^*(N-1)}{N^2}\int_\Omega F^*(- Du(x))\, u(x)^{(1-\frac{1}{N})p_N^*-1}  \, d_F(x_0,\mathcal F_1(x)) {\rm d}{\sf m}(x).
\end{align}
Note that the constant $C_0>0$ is independent of both $k$ and $\lambda$.
Since $\left( (1-\frac{1}{N})p_N^*-1\right)p'=p_N^*$, H\"older's inequality, the Monge--Amp\`ere equation \eqref{MA-Sob} and a change of variable imply 
\begin{align}\label{estimate-jobb}
	\nonumber	J&:=\int_\Omega F^*(- Du(x))\, u(x)^{(1-\frac{1}{N})p_N^*-1}  \, d_F(x_0,\mathcal F_1(x)) {\rm d}{\sf m}(x)\\&	\nonumber \leq \left(\int_\Omega F^*(-Du)^p {\rm d}{\sf m}  \right)^\frac{1}{p}\left(\int_\Omega u(x)^{\left( (1-\frac{1}{N})p_N^*-1\right)p'} d_F^{p'}(x_0,\mathcal F_1(x)){\rm d}{\sf m}(x)\right)^\frac{1}{p'}
	\\&	\nonumber = \|F^*(-Du)\|_{L^p(M,{\sf m})}\left(\int_\Omega \rho_0(x) d_F^{p'}(x_0,\mathcal F_1(x)){\rm d}{\sf m}(x)\right)^\frac{1}{p'}\\&	\nonumber = \|F^*(-Du)\|_{L^p(M,{\sf m})}\left(\int_\Omega \rho_1(\mathcal F_1(x))\, d_F^{p'}(x_0,\mathcal F_1(x))\, {\bf D}[d(\mathcal F_1)_x]{\rm d}{\sf m}(x)\right)^\frac{1}{p'}
	\\&  = \|F^*(-Du)\|_{L^p(M,{\sf m})}\left(\int_M \rho_1(y)\, d_F^{p'}(x_0,y)\, {\rm d}{\sf m}(y)\right)^\frac{1}{p'}
.	\end{align}
Due to the properties of the cut-off function $\eta_k$ appearing in the density $\rho_1$,  the   dominated convergence theorem implies  that 
$$
\lim_{k\to \infty}	\int_M \rho_1^{1-\frac{1}{N}}{\rm d}{\sf m}=(I_{\lambda,N}^{x_0})^{-1+\frac{1}{N}}\, I_{\lambda,N-1}^{x_0} 
$$
 and 
$$\lim_{k\to \infty}\int_M \rho_1(y)\, d_F^{p'}(x_0,y)\, {\rm d}{\sf m}(y)=(I_{\lambda,N}^{x_0})^{-1}\, I_{\lambda,N-1}^{x_0} -\lambda. $$
Thus, relations \eqref{majdnemutolso} and \eqref{estimate-jobb} give that
$$(I_{\lambda,N}^{x_0})^{-1+\frac{1}{N}}\, I_{\lambda,N-1}^{x_0} \leq C_0+\frac{p_N^*(N-1)}{N^2}\|F^*(-Du)\|_{L^p(M,{\sf m})}\left((I_{\lambda,N}^{x_0})^{-1}\, I_{\lambda,N-1}^{x_0} -\lambda\right)^\frac{1}{p'}.$$
Since $C_0>0$ does not depend on $\lambda>0$,  dividing the latter inequality by $\lambda^\frac{1}{p'}>0$ and taking the limit $\lambda\to \infty$, relation \eqref{AVR-in-Limit} implies that
$$ \left({\sf AVR}^+_{\mathcal M}\right)^\frac{1}{N}\cdot {\sf AT}(p,N)\leq \|F^*(-Du)\|_{L^p(M,{\sf m})},$$
which is   the claimed inequality \eqref{Sob-pozitiv} under the  normalization \eqref{u-normalization}.  

\textit{Step 2.} 	We prove that
\begin{equation}\label{Sob-minden}
	\|F^*(Du)\|_{L^p(M,{\sf m})} \geq  {\sf AVR} _{\mathcal M} ^\frac{1}{N}\cdot {\sf AT}(p,N)\, \|u\|_{L^{p_N^*}(M,{\sf m})},\ \ \forall u\in W^{p}_{p_N^*}(\mathcal M), 
\end{equation}
 where ${\sf AVR}_{\mathcal M}=\min\{{\sf AVR}_{\mathcal M}^+,{\sf AVR}_{\mathcal M}^-\}$; note that no sign assumption is imposed on $u$. 
 
 To do this, we note that relation \eqref{Sob-pozitiv} for the reverse Finsler metric $\overleftarrow F$ reads as follows:   for every \textit{non-negative} $u$ with $- u\in W^{p}_{p_N^*}(\overleftarrow {\mathcal M})$ (or equivalently,  $u\in W^{p}_{p_N^*}(\mathcal M)$) one has 
  \begin{equation}\label{Sob-reverse}
  	\|\overleftarrow F^*(-Du)\|_{L^p(M,{\sf m})} \geq \left({\sf AVR}^-_{\mathcal M}\right)^\frac{1}{N}\cdot {\sf AT}(p,N)\, \|u\|_{L^{p_N^*}(M,{\sf m})}.  
  \end{equation}
			Indeed, we recall that the bound ${\sf Ric}_N\geq 0$ for 
			 $\mathcal M=(M,F,{\sf m})$ is inherited also by the reverse  Finsler space  $\overleftarrow {\mathcal M}=(M,\overleftarrow F,{\sf m})$; moreover, the latter is also forward complete as   $\mathcal M=(M,F,{\sf m})$ is  complete. In addition, ${\sf AVR}^+$ is also replaced by  ${\sf AVR}^-$, as we   work with the reverse distance $\overleftarrow {d_F}$ instead of $d_F$. It remains to follow the arguments as in \textit{Step 1}, which yields relation \eqref{Sob-reverse}. 
			 
			 Now, we focus on the proof of \eqref{Sob-minden}. Let us fix $u\in W^{p}_{p_N^*}(\mathcal M)$ arbitrarily  and consider its positive and negative parts, i.e., $u_+=\max\{u,0\} \geq 0$ and $u_-=-\min\{u,0\}\geq 0$; it is clear that $u=u_+-u_- $. In addition,  up to an ${\sf m}$-null set, $u_+$ and $u_-$ have disjoint essential supports and     one has ${\sf m}$-a.e.\ that
			 \begin{equation}\label{pos-neg-splitting-F}
			 	 F^*( Du)^p= F^*( Du_+-Du_-)^p =F^*( Du_+)^p + F^*(  -Du_-)^p
			 \end{equation} 
				and 
				 \begin{equation}\label{pos-neg-splitting-u}
					|u|^{p_N^*}=   u_+ ^{p_N^*}+  u_- ^{p_N^*}.
				\end{equation} 
			In particular, it follows that  $u_+\in W^{p}_{p_N^*}(\mathcal M)$ and $u_-\in W^{p}_{p_N^*}(\overleftarrow{\mathcal M})$. 			
				Therefore, using \eqref{pos-neg-splitting-F} together with   \eqref{Sob-reverse} for $u_+$ and \eqref{Sob-pozitiv} for $u_-$, one has that 
				\begin{align*}
					\|F^*(Du)\|_{L^p(M,{\sf m})}^p&=\|F^*(Du_+)\|_{L^p(M,{\sf m})}^p+\|F^*(-Du_-)\|_{L^p(M,{\sf m})}^p\\&=\|\overleftarrow F^*(-Du_+)\|_{L^p(M,{\sf m})}^p+\|F^*(-Du_-)\|_{L^p(M,{\sf m})}^p\\&\geq {\sf AT}(p,N)^p\left(\left({\sf AVR}^-_{\mathcal M}\right)^\frac{p}{N}  \|u_+\|^p_{L^{p_N^*}(M,{\sf m})}+ \left({\sf AVR}^+_{\mathcal M}\right)^\frac{p}{N} \|u_-\|^p_{L^{p_N^*}(M,{\sf m})}\right)\\&\geq {\sf AT}(p,N)^p\cdot  {\sf AVR} _{\mathcal M} ^\frac{p}{N} \left( \|u_+\|^p_{L^{p_N^*}(M,{\sf m})}+   \|u_-\|^p_{L^{p_N^*}(M,{\sf m})}\right).
				\end{align*}
	Recall the elementary fact that the function $s\mapsto (a^s+b^s)^\frac{1}{s}$ is	decreasing on $(0,\infty)$ for every fixed $a,b>0 $ (and constant if either $a=0$ or $b=0$). Thus, by the   above estimate and \eqref{pos-neg-splitting-u},  we obtain that
	\begin{align}\label{utolso-becsles-a-sobolevben}
		\nonumber \|F^*(Du)\|_{L^p(M,{\sf m})}&\geq {\sf AT}(p,N) \cdot  {\sf AVR} _{\mathcal M} ^\frac{ 1}{N} \left( \|u_+\|^p_{L^{p_N^*}(M,{\sf m})}+   \|u_-\|^p_{L^{p_N^*}(M,{\sf m})}\right)^\frac{1}{p} \\&\nonumber \geq {\sf AT}(p,N) \cdot  {\sf AVR} _{\mathcal M} ^\frac{ 1}{N} \left( \|u_+\|^{p_N^*}_{L^{p_N^*}(M,{\sf m})}+   \|u_-\|^{p_N^*}_{L^{p_N^*}(M,{\sf m})}\right)^\frac{1}{p_N^*} \\& = {\sf AT}(p,N)\cdot   {\sf AVR} _{\mathcal M} ^\frac{ 1}{N}\,   \|u \| _{L^{p_N^*}(M,{\sf m})},
	\end{align}		
		which is the required inequality \eqref{Sob-minden}, equivalent to  $\mathcal R^{p}_{p_N^*}(\mathcal M)\geq {\sf AT}(p,N)\cdot   {\sf AVR} _{\mathcal M} ^\frac{ 1}{N}$.

			\subsubsection{Optimality: $\mathcal R^{p}_{p_N^*}(\mathcal M)= {\sf AT}(p,N)\cdot   {\sf AVR} _{\mathcal M} ^\frac{ 1}{N}$} Assume by contradiction that   there exists $C>0$ such that  $\mathcal R^{p}_{p_N^*}(\mathcal M)> C> {\sf AT}(p,N)\cdot   {\sf AVR} _{\mathcal M} ^\frac{ 1}{N};$ in particular,
			\begin{equation}\label{contradiction}
				\|F^*(Du)\|_{L^p(M,{\sf m})} \geq  C\, \|u\|_{L^{p_N^*}(M,{\sf m})},\ \ \forall u\in W^{p}_{p_N^*}(\mathcal M). 
			\end{equation}
	 
	 Let us assume first that ${\sf AVR} _{\mathcal M}={\sf AVR} _{\mathcal M}^+$. Then, for some $x_0\in M $ and every $\lambda>0$, we consider the function $u(x)=-(\lambda+ d_F^{p'}(x_0,x))^{\frac{p-N}{p}}$, $x\in M$. Accordingly, using the notation from \eqref{layer-cake}, one has
	 $$\|u\|_{L^{p_N^*}(M,{\sf m})}=\left(\int_M(\lambda+ d_F^{p'}(x_0,x))^{-N}{\rm d}{\sf m}(x)\right)^\frac{1}{p_N^*}=(I_{\lambda,N}^{x_0})^\frac{1}{p_N^*}. $$
		Since $Du(x)=\frac{N-p}{p}p'(\lambda+ d_F^{p'}(x_0,x))^{\frac{ -N}{p}} d_F^{p'-1}(x_0,x) Dd_F (x_0,x) $	for ${\sf m}$-a.e.\ $x\in M$, 
		the eikonal equation \eqref{eikonal} implies that 
		\begin{align*}
			 \|F^*(Du)\|_{L^p(M,{\sf m})}&=\frac{N-p}{p}p'\left(\int_M(\lambda+ d_F^{p'}(x_0,x))^{-N}d_F^{p'}(x_0,x){\rm d}{\sf m}(x)\right)^\frac{1}{p}\\&=\frac{N-p}{p}p'\left(I_{\lambda,N-1}^{x_0}-\lambda I_{\lambda,N}^{x_0}\right)^\frac{1}{p}.
		\end{align*}
		Inserting   these relations into  \eqref{contradiction}, we obtain for every $\lambda>0$ that
		$$\frac{N-p}{p}p'\left(I_{\lambda,N-1}^{x_0}-\lambda I_{\lambda,N}^{x_0}\right)^\frac{1}{p}\geq C\, (I_{\lambda,N}^{x_0})^\frac{1}{p_N^*}.$$
		Multiplying the latter inequality by $\lambda^{\frac{N-p}{p^2}}>0$ and taking the limit $\lambda\to \infty$, we obtain by \eqref{AVR-in-Limit} that  ${\sf AT}(p,N)\cdot   \left({\sf AVR}^+ _{\mathcal M}\right) ^\frac{ 1}{N}\geq C,$ 
		which contradicts the choice of $C>0.$
		
		When ${\sf AVR} _{\mathcal M}={\sf AVR} _{\mathcal M}^-$, we proceed in a similar way as above, considering in \eqref{contradiction} the function $u(x)= (\lambda+ {\overleftarrow {d_F}}^{p'}(x_0,x))^{\frac{p-N}{p}}$, $x\in M$. This ends the proof of Theorem \ref{Thm-Sobolev-1}. \hfill $\square$
		
	\begin{remark}\label{remark-locally-convex} \rm We prove that  when $\partial M\neq \emptyset$, its  local  convexity   is indispensable in Theorem \ref{Thm-Sobolev-1}. To do so, let $n\geq 3$ and the space $\mathcal M$ consisting of the set $M=\mathbb R^n\setminus \mathbb B^n$, endowed with the classical Euclidean metric and Lebesgue measure. It is clear that $\partial M=\mathbb S^{n-1}$ is not locally convex,   ${\sf Ric}=0$, $N=n$ and ${\sf AVR}_{\mathcal M}=1$. We consider the Talentian function $u(x)=(1+|x|^{p'})^\frac{p-n}{p}$,   $x_0\in \partial M=\mathbb S^{n-1}$ 	and $u_\epsilon(x)=u((x-x_0)/\epsilon)$ for $\epsilon>0.$ For convenience, let $M_\epsilon=(M-x_0)/\epsilon.$ By change of variables, one has that 
		$$\int_{M}|u_\epsilon|^{p*}=\epsilon^n\int_{M_\epsilon}|u |^{p*}\ \  {\rm and}\ \  \int_{M}|\nabla u_\epsilon|^{p}=\epsilon^{n-p}\int_{M_\epsilon}|\nabla u |^{p}.$$
	One can check that $M_\epsilon \to H_{x_0}=\{x\in \mathbb R^n:x\cdot x_0\geq 0\}$ as $\epsilon\to 0$. Therefore, since $n/p^*=(n-p)/p$, it follows that 
	\begin{align*}
		\mathcal R^{p}_{p^*}(\mathcal M)&\leq \lim_{\epsilon\to 0} \frac{\|\nabla u_\epsilon\|_{L^p(M,{\rm d}x)}}{\|  u_\epsilon\|_{L^{p^*}(M,{\rm d}x)}}=\lim_{\epsilon\to 0} \frac{\|\nabla u\|_{L^p(M_\epsilon,{\rm d}x)}}{\|  u\|_{L^{p^*}(M_\epsilon,{\rm d}x)}}=\frac{\|\nabla u\|_{L^p(H_{x_0},{\rm d}x)}}{\|  u\|_{L^{p^*}(H_{x_0},{\rm d}x)}}\\&=\frac{2^{-1/p}\|\nabla u\|_{L^p(\mathbb R^n,{\rm d}x)}}{2^{-1/p^*}\|  u\|_{L^{p^*}(\mathbb R^n ,{\rm d}x)}}=2^{-1/n} {\sf AT}(p,n)<{\sf AT}(p,n),
	\end{align*}
	which contradicts relation \eqref{optim-sob-ineq}.
	 	
	\end{remark}
		 
			\subsection{Rigidity: proof of Theorem \ref{Thm-Sobolev-2-equality}}

	(i)$\implies$(ii) We first recall the inequality of Bliss \cite{Bliss}, which states that on the    model space $(\mathbb R_+,\vert\cdot\vert,{\sf m}_N=N\omega_Nr^{N-1}{\mathcal L ^1}_{\llcorner{[0,\infty)}})$ one has for every locally absolutely continuous function $v:\mathbb R_+\to \mathbb R$    that 
	\begin{equation}\label{Bliss-ineq}
		\|v'\|_{L^p(\mathbb R_+,{\sf m}_N)}\geq {\sf AT}(p,N) \|v\|_{L^{p_N^*}(\mathbb R_+,{\sf m}_N)}.
	\end{equation}
	Moreover, equality holds in \eqref{Bliss-ineq} -- up to a multiplicative factor -- only for the family of functions $v_\lambda(r)=(\lambda+r^{p'})^\frac{p-N}{p}$, $\lambda>0.$

		Assume now that $\mathcal R^p_{p^*_N}(\mathcal M)={\sf AVR}_{\mathcal M}^\frac{1}{N}\cdot {\sf AT}(p,N)$ is achieved  by   $u\in   W^p_{p^*_N}(\mathcal M)\setminus\{0\}$, i.e., 
	\begin{equation}\label{equality-in sobolev}
		\|F^*(Du)\|_{L^p(M,{\sf m})}={\sf AT}(p,N)\cdot   {\sf AVR} _{\mathcal M} ^\frac{ 1}{N}\,   \|u \| _{L^{p_N^*}(M,{\sf m})}.
	\end{equation}
	 		 Let $u_+=\max\{u,0\} \geq 0$ and $u_-=-\min\{u,0\}\geq 0$ be the positive and negative parts of $u$, verifying the properties \eqref{pos-neg-splitting-F} and \eqref{pos-neg-splitting-u}. By \eqref{pos-neg-splitting-F}, the backward/forward   P\'olya--Szeg\H o inequalities \eqref{Polya-Szego},   the Bliss' inequality \eqref{Bliss-ineq} and Cavalieri principle \eqref{Cavallieri}, we have that 
		\begin{align*}
		\|F^*(Du)\|_{L^p(M,{\sf m})}^p&=\|F^*(Du_+)\|_{L^p(M,{\sf m})}^p+\|F^*(-Du_-)\|_{L^p(M,{\sf m})}^p\\&=\|\overleftarrow F^*(-Du_+)\|_{L^p(M,{\sf m})}^p+\|F^*(-Du_-)\|_{L^p(M,{\sf m})}^p\\&\geq  ( {\sf AVR}^-_{\mathcal M})^{p/N}\|(u_+^*)'\|^p_{L^p(\mathbb R_+,{\sf m}_N)} + ( {\sf AVR}^+_{\mathcal M})^{p/N}\left(\|(u_-^*)'\|^p_{L^p(\mathbb R_+,{\sf m}_N)}\right)\\&\geq {\sf AVR}_{\mathcal M}^{p/N}{\sf AT}(p,N)^p\left( \|u_+^*\|^p_{L^{p_N^*}(\mathbb R_+,{\sf m}_N)} + \|u_-^*\|^p_{L^{p_N^*}(\mathbb R_+,{\sf m}_N)}\right)
		\\&= {\sf AVR}_{\mathcal M}^{p/N}{\sf AT}(p,N)^p \left( \|u_+\|^p_{L^{p_N^*}(M,{\sf m})}+   \|u_-\|^p_{L^{p_N^*}(M,{\sf m})}\right).
	\end{align*}
	By the latter estimate, similarly as in \eqref{utolso-becsles-a-sobolevben}, we  obtain that 
	$$\|F^*(Du)\|_{L^p(M,{\sf m})}\geq {\sf AT}(p,N)\cdot   {\sf AVR} _{\mathcal M} ^\frac{ 1}{N}\,   \|u \| _{L^{p_N^*}(M,{\sf m})}.$$
	By \eqref{equality-in sobolev}  it follows that in the above inequalities, we have necessarily equality. In particular, the   decreasing character of $s\mapsto (a^s+b^s)^\frac{1}{s}$, $s>0,$ which is the key in \eqref{utolso-becsles-a-sobolevben}, we   have either $ u_+=0$ or $u_-=0$, i.e., $u$ is sign-preserving. 
	
	Let us assume that $u_+=0$; thus $u=-u_-$. Note that we also have equality in the forward P\'olya--Szeg\H o inequality. i.e., 
$$	\|F^*(-Du_-)\|_{L^p(M,{\sf m})}=   ( {\sf AVR}^+_{\mathcal M})^{1/N}\|(u_-^*)'\|_{L^p(\mathbb R_+,{\sf m}_N)},$$
and ${\sf AVR}_{\mathcal M}={\sf AVR}_{\mathcal M}^+$.
According to \S\ref{subsection-rearr},  there exists  $x_0\in M$ such that  $\mathcal M=(M,F,{\sf m})$  is a  forward  $N$-volume cone  with  pole $x_0$, i.e., 
	$
		\frac{ {\sf m}(B_{x_0}^+(R))}{\omega_NR^N}={\sf AVR}_{\mathcal M}^+$ for every $ R>0.$ Due to Proposition \ref{N-volume-cone-prop}, one has that for every forward   geodesic $\gamma:[0,\infty)\to M$ 
		starting from $\gamma(0)=x_0,$ we have ${\sf K}(\dot \gamma(t),\cdot)=0 $  and ${\sf S}_{\sf m}(\gamma(t),\dot \gamma(t))=\frac{n-N}{t}$  for every $t>0.$ 
		
		At the same time, we have equality also in the Bliss inequality, i.e., 
		$$\|(u^*_-)'\|_{L^p(\mathbb R_+,{\sf m}_N)}= {\sf AT}(p,N) \|u^*_-\|_{L^{p_N^*}(\mathbb R_+,{\sf m}_N)},$$
		which implies that $u^*_-(r)=(\lambda+\mu r^{p'})^\frac{p-N}{p}$ for some $\lambda,\mu>0.$
		Moreover, due to \eqref{u=u^*}, it turns out that for ${\sf m}{\rm -a.e.}\  x\in M$, 
	$u_-(x)=u_-^*\left({\sf AVR}^+_{\mathcal M})^{1/N} d_F(x_0,x)\right),$ which concludes the proof together with $u=-u_-$. Clearly, when $u_-=0$, we proceed in a similar manner. 
				
		(ii)$\implies$(i) Assume that ${\sf AVR}_{\mathcal M}={\sf AVR}_{\mathcal M}^+$. According to Proposition \ref{N-volume-cone-prop}, from (ii)/(a) it follows that $\mathcal M=(M,F,{\sf m})$  is a  forward  $N$-volume cone  with  pole $x_0$, i.e.,  $ {{\sf m}(B^+_{x_0}(R))}={\sf AVR}_{\mathcal M}^+ \omega_N { R^N}$ for every $R>0.$ If $u(x)=-(\lambda+\mu d_F(x_0,x)^{p'})^\frac{p-N}{p}$ for some $\lambda,\mu>0$, by the $N$-volume cone character of $\mathcal M=(M,F,{\sf m})$, the layer cake representation and eikonal equation \eqref{eikonal}, one has that 
		$$\|F^*(Du)\|_{L^p(M,{\sf m})}^p={\sf AVR}_{\mathcal M}^+\|u_0'\|^p_{L^{p }(\mathbb R_+,{\sf m}_N)},$$
		where $u_0(r)=(\lambda+\mu r^{p'})^\frac{p-N}{p}$. In a similar manner, one has that 
		$$\|u\|_{L^{p_N^*}(M,{\sf m})}^{p_N^*}={\sf AVR}_{\mathcal M}^+\|u_0 \|^{p_N^*}_{L^{p_N^*}(\mathbb R_+,{\sf m}_N)}.$$
		Since $u_0$ produces equality in the Bliss inequality \eqref{Bliss-ineq}, the function $u$ saturates the Rayleigh quotient   $\mathcal R^p_{p^*_N}(\mathcal M)$, i.e., $\|F^*(Du)\|_{L^p(M,{\sf m})}={\sf AT}(p,N)\cdot   {\sf AVR} _{\mathcal M} ^\frac{ 1}{N}\,   \|u \| _{L^{p_N^*}(M,{\sf m})},$ which concludes the proof. \hfill $\square$

			\section{Optimal logarithmic Finsler--Sobolev inequality: proof of Theorem \ref{log-Sobolev-main}} \label{section-4} 
		
		\subsection{Proof of \eqref{LSI}} 
		
		{\textit{Step 1.}} We first prove that for every  non-negative, compactly supported   function $u\in C_0^\infty(M)$ with $\|u\|_{L^{p }(M,{\sf m})} =1,$ one has that 
		\begin{equation}\label{LSI-positive}
				\int_{M} u ^p\log  u ^p{\rm d} {\sf m}\leq \frac{N}{p}\log\left(\mathcal
			L_{p,N}\left({\sf AVR}^+_{\mathcal M}\right)^{-\frac{p}{N}}\|F^*(-Du)\|_{L^p(M,{\sf m})}^p \right);
		\end{equation} 
	by density arguments, it can be extended to non-negative functions belonging to $   W^{p}_p(\mathcal M)$. 
		Let $\Omega=\{x\in M:u(x)>0\}$; in particular, $\overline \Omega\subset \overline M$ is compact.
		
		First, we establish some properties of  Gaussian-type functions on $M.$ 	Fixing   $\lambda>0$ and $x_0\in M$, let 
	\begin{equation}\label{notations-J-k}
			J_{1,\lambda}^{x_0}=\int_M e^{-\lambda d_F(x_0,x)^{p'}}{\rm d}{\sf m}(x)\ \ {\rm and}\ \ J_{2,\lambda}^{x_0}=\int_M d_F(x_0,x)^{p'} e^{-\lambda d_F(x_0,x)^{p'}}{\rm d}{\sf m}(x).
	\end{equation}
			In a similar way as in \eqref{layer-cake}, one can prove that $J_{i,\lambda}^{x_0}<+\infty$ for every $\lambda>0$ and $i\in \{1,2\}$, and $\frac{\rm d}{{\rm d}\lambda}J_{1,\lambda}^{x_0}=-J_{2,\lambda}^{x_0}.$ Moreover,   the layer cake representation implies that 
	\begin{equation}\label{limitek-J}
		\lim_{\lambda\to 0}\lambda^\frac{N}{p'}J_{1,\lambda}^{x_0}=\omega_N {\sf AVR}_{\mathcal M}^+ \Gamma \left(\frac{N}{p'}+1\right) \ \ {\rm and}\ \  \lim_{\lambda\to 0}\lambda^{\frac{N}{p'}+1}J_{2,\lambda}^{x_0}=\omega_N {\sf AVR}_{\mathcal M}^+ \frac{N}{p'}\Gamma \left(\frac{N}{p'}+1\right).
	\end{equation}
		
		We now consider the probability measures $\mu_0=\rho_0 {\sf m}$ and $\mu_1=\rho_1 {\sf m}$ with 
		$$ \rho_0(x)=u^{p}(x) \ \ \text{and}\ \ \rho_1(x)=(J_{\lambda, k}^{x_0})^{-1} \eta_k(x) e^{-\lambda d_F(x_0,x)^{p'}},$$
		where $\eta_k$ is the cut-off function from \S \ref{ssss-us} and $$J_{\lambda, k}^{x_0}=\int_M \eta_k(x) e^{-\lambda d_F(x_0,x)^{p'}}{\rm d}{\sf m}(x).$$
	In fact, we have that  $\mu_0, \mu_1\in \mathcal P_c(M)$. In particular, we have  the Monge--Amp\`ere equation 
	\begin{equation}\label{MA-logSob}
		\rho_0(x)=\rho_1(\mathcal F_1(x))\cdot {\bf D}[d(\mathcal F_1)_x]\ \ \ \ {\rm for}\ \mu_0-{\rm a.e.}\ x\in \Omega,
	\end{equation}
		where $\mathcal F_t(x)=\exp_x(\bm \nabla(-t\phi)(x))$, $t\in [0,1],$ is the unique optimal transport map between the probability measures $\mu_0$ and $\mu_1$, the function $\phi:M\to \mathbb R$ being $c$-concave.   As before, $\rho_1$, $\phi$ and $\mathcal F_t$ all depend on $k$ and $\lambda$, which will not be explicitly stated in the sequel.   
		
			Proposition \ref{displacement-convexity}/(i)\&(iii) and \eqref{distributional-alexandrv} imply that  
		\begin{equation}\label{disp-convexity-logSobolev}
			\int_M\rho_0 \log \rho_0 {\rm d}{\sf m}\leq \int_M\rho_1 \log \rho_1 {\rm d}{\sf m}+ 
			 N\log	\left(1+\frac{1}{N}\displaystyle\int_M    \rho_0{\bm\Delta_{\mathcal M}(-\phi) } {\rm d}{\sf m}\right). 	 
					\end{equation}
		According to   \eqref{Laplacian} and \eqref{Cauchy-Schwarz}, it follows that 		
				$$\displaystyle\int_M    \rho_0{\bm\Delta_{\mathcal M}(-\phi) } {\rm d}{\sf m}=-p\int_M    u^{p-1}Du(\bm \nabla(-\phi))   {\rm d}{\sf m}\leq p\int_M    u^{p-1}F^*(-Du)F(\bm \nabla(-\phi))   {\rm d}{\sf m}.$$
				Similarly as in \eqref{F-finsler-estimate-distance}, it follows that
				     for $\mu_0$-a.e.\ $x\in \Omega$,  
				\begin{align}\label{F-finsler-logestimate-distance}
					\nonumber F(\bm \nabla (-\phi)(x)) \leq K+d_F(x_0,\mathcal F_1(x)), 
				\end{align}
				where $K=\max_{x\in \overline\Omega}d_F(x,x_0).$ If $C:=1+\frac{Kp}{N}\displaystyle\int_{\Omega}u^{p-1}F^*(-Du){\rm d}{\sf m}>0$, we obtain by H\"older's inequality and Monge--Amp\`ere equation \eqref{MA-logSob} that  
				\begin{align*}
						1+\frac{1}{N}\displaystyle\int_M    \rho_0{\bm\Delta_{\mathcal M}(-\phi) } {\rm d}{\sf m}\leq& C+\frac{p}{N}\displaystyle\int_M     u(x)^{p-1}F^*(-Du(x))d_F(x_0,\mathcal F_1(x)){\rm d}{\sf m}(x)\\ \leq & C+\frac{p}{N}\displaystyle \|F^*(-Du )\|_{L^p (M,{\sf m})}\left(\int_M     \rho_1(y)d_F(x_0,y)^{p'}{\rm d}{\sf m}(y)\right)^{1/p'}.
				\end{align*}
			The latter estimate jointly with  \eqref{disp-convexity-logSobolev} gives that
			\begin{align*}
				\int_M u^p \log u^p {\rm d}{\sf m}\leq& \int_M   \frac{\eta_k(x) e^{-\lambda d_F(x_0,x)^{p'}}}{J_{\lambda, k}^{x_0}} \log \frac{\eta_k(x) e^{-\lambda d_F(x_0,x)^{p'}}}{J_{\lambda, k}^{x_0}} {\rm d}{\sf m}\\&+ 
				N\log	\Bigg(C+\frac{p}{N}\displaystyle \|F^*(-Du )\|_{L^p (M,{\sf m})}\Bigg.\times\\&\qquad \qquad  \times\left.\left(\int_M     \frac{\eta_k(y) e^{-\lambda d_F(x_0,y)^{p'}}}{J_{\lambda, k}^{x_0}}d_F(x_0,y)^{p'}{\rm d}{\sf m}(y)\right)^{1/p'}\right).
			\end{align*} 
				Letting $k\to \infty$, a    reorganization of the latter inequality yields through the dominated  convergence theorem and notations from \eqref{notations-J-k}
				 that 
		\begin{align*}
			\int_M u^p \log u^p {\rm d}{\sf m}\leq&-\lambda\frac{	J_{2,\lambda}^{x_0}}{	J_{1,\lambda}^{x_0}}   + 
			N\log	\Bigg(\frac{C}{	(J_{1,\lambda}^{x_0})^{1/N}}+\frac{p}{N}\displaystyle \|F^*(-Du )\|_{L^p (M,{\sf m})}\Bigg.\frac{(	J_{2,\lambda}^{x_0})^{1/p'}}{	(J_{1,\lambda}^{x_0})^{1/N+1/p'}}\Bigg).
		\end{align*} 
		Taking the limit $\lambda\to 0$, due to \eqref{limitek-J}, it yields exactly \eqref{LSI-positive}. 
		
		{\textit{Step 2.}} We now prove the general form  \eqref{LSI}, namely, for every $u\in W_p^p(\mathcal M)$ with $\|u\|_{L^p(M,{\sf m})}=1,$ one has
			\begin{equation}\label{LSI-szoveg}
			\int_{M}|u|^p\log |u|^p{\rm d} {\sf m}\leq \frac{N}{p}\log\left(\mathcal
			L_{p,N}{\sf AVR}_ {\mathcal M}^{-\frac{p}{N}}\|F^*(Du)\|_{L^p(M,{\sf m})}^p \right). 
		\end{equation}
	
	 Let   $u\in W^{p}_{p }(\mathcal M)$ with $\|u\|_{L^p(M,{\sf m})}=1,$ and consider its positive and negative parts, i.e., $u_+=\max\{u,0\} \geq 0$ and $u_-=-\min\{u,0\}\geq 0$; it is clear that $u=u_+-u_- $. Up to an ${\sf m}$-null set, $u_+$ and $u_-$ have disjoint essential supports and     one has ${\sf m}$-a.e.\ that
	\begin{equation}\label{pos-neg-splitting-F-log}
		F^*( Du)^p= F^*( Du_+-Du_-)^p =F^*( Du_+)^p + F^*(  -Du_-)^p
	\end{equation} 
	and 
	\begin{equation}\label{pos-neg-splitting-u-log}
		|u|^{p}=   u_+ ^{p}+  u_- ^{p}.
	\end{equation} 
Since    $u \in W^{p}_{p }(\mathcal M),$ then it follows that    $u_+\in W^{p}_{p }(\mathcal M)$ and   $u_-\in W^{p}_{p }(\overleftarrow{\mathcal M}).$
 If either $u_+=0$ or $u_-=0$, the proof of \eqref{LSI-szoveg} is trivial; thus, we assume $U_\pm:=\|u_\pm\|_{L^p(M,{\sf m})}>0$. In particular, we also have that
 $1=\|u \|_{L^p(M,{\sf m})}^p=\|u_+ \|_{L^p(M,{\sf m})}^p+\|u _-\|_{L^p(M,{\sf m})}^p=U_+^p+U_-^p.$  
	According to \eqref{LSI-positive}, one has that 
\begin{equation}\label{LSI--}
	\int_{M}  {u_- ^p}\log  \frac{u_- ^p}{U_-^p}{\rm d} {\sf m}\leq {U_-^p}\frac{N}{p}\log\left(\mathcal
	L_{p,N}\left({\sf AVR}^+_{\mathcal M}\right)^{-\frac{p}{N}}\frac{\|F^*(-Du_-)\|_{L^p(M,{\sf m})}^p}{ U_-^p}  \right), 
\end{equation} 
 and   
		\begin{equation}\label{LSI-+}
			\int_{M}  {u_+ ^p}\log  \frac{u_+ ^p}{U_+^p}{\rm d} {\sf m}\leq {U_+^p}\frac{N}{p}\log\left(\mathcal
			L_{p,N}\left({\sf AVR}^-_{\mathcal M}\right)^{-\frac{p}{N}}\frac{\|\overleftarrow F^*(-Du_+)\|_{L^p(M,{\sf m})}^p}{ U_+^p}  \right). 
		\end{equation} 
	
	On the one hand, by the convexity of $s\mapsto s\log s$, $s\geq 0$, relation \eqref{pos-neg-splitting-u-log} and  $U_+^p+U_-^p=1$, it follows that 
	\begin{align}\label{bal-oldal}
	\nonumber	{\rm LHS}&:=\int_{M}\left(  {u_+ ^p}\log  \frac{u_+ ^p}{U_+^p}+  {u_- ^p}\log  \frac{u_- ^p}{U_-^p}\right){\rm d} {\sf m}=\int_{M}\left(  {U_+ ^p}\frac{u_+ ^p}{U_+^p}\log  \frac{u_+ ^p}{U_+^p}+  {U_- ^p}\frac{u_- ^p}{U_-^p}\log  \frac{u_- ^p}{U_-^p}\right){\rm d} {\sf m}\\&\geq \int_M(u_+^p+u_-^p)\log(u_+^p+u_-^p){\rm d} {\sf m}=\int_M |u|^p \log |u|^p{\rm d} {\sf m}.
	\end{align}
	On the other hand, by the concavity of $s\mapsto  \log s$, $s\geq 0,$ and \eqref{pos-neg-splitting-F-log}, it follows that 
		\begin{align}\label{jobb-oldal}
			\nonumber	{\rm RHS} :=&{U_+^p}\frac{N}{p}\log\left(\mathcal
			L_{p,N}\left({\sf AVR}^-_{\mathcal M}\right)^{-\frac{p}{N}}\frac{\|\overleftarrow F^*(-Du_+)\|_{L^p(M,{\sf m})}^p}{ U_+^p}  \right) \\&+{U_-^p}\frac{N}{p}\log\left(\mathcal
			L_{p,N}\left({\sf AVR}^+_{\mathcal M}\right)^{-\frac{p}{N}}\frac{\|F^*(-Du_-)\|_{L^p(M,{\sf m})}^p}{ U_-^p}  \right)\nonumber \\ \leq& \frac{N}{p}\log\left(\mathcal
			L_{p,N} {\sf AVR}_{\mathcal M} ^{-\frac{p}{N}} { \| F^*( Du)\|_{L^p(M,{\sf m})}^p}   \right).
		\end{align}
	It remains to combine relations \eqref{LSI--}-\eqref{jobb-oldal}, which gives \eqref{LSI-szoveg}, equivalent to \eqref{LSI}.  
	
	\subsection{Optimality of \eqref{LSI}} We proceed by contradiction, i.e., we assume  that   there exists  $0<C<\mathcal
	L_{p,N} {\sf AVR}_{\mathcal M} ^{-\frac{p}{N}}$   such that for every 
	$ u\in W^{p}_p(\mathcal M)$ with $\|u\|_{L^p(M,{\sf m})}=1,$ one has   
	\begin{equation}\label{LSI-contra}
		\int_{M}|u|^p\log |u|^p{\rm d} {\sf m}\leq \frac{N}{p}\log\left(C\, \|F^*(Du)\|_{L^p(M,{\sf m})}^p \right). 
	\end{equation}
	We may assume, without loss of generality, that ${\sf AVR} _{\mathcal M}={\sf AVR} _{\mathcal M}^+$; the other case is similar under minimal changes. Fixing $x_0\in M $ and   $\lambda>0$, we consider the Gaussian bubble $u_\lambda(x)=-(J_{1,\lambda}^{x_0})^{-\frac{1}{p}}e^{-\frac{\lambda}{p}d_F(x_0,x)^{p'}} $, $x\in M$, keeping the notation from \eqref{notations-J-k}. By construction, it follows that $ u_\lambda\in W^{p}_p(\mathcal M)$ and  $\|u_\lambda\|_{L^p(M,{\sf m})}=1,$ thus this function can be used in \eqref{LSI-contra}.
	Note that  on account of \eqref{notations-J-k}, it follows that 
	$$	\int_{M}|u_\lambda|^p\log |u_\lambda|^p{\rm d} {\sf m}=-\log J_{1,\lambda}^{x_0} -\lambda (J_{1,\lambda}^{x_0})^{-1} J_{2,\lambda}^{x_0}.$$ Moreover, by the eikonal equation \eqref{eikonal} one has for ${\sf m}$-a.e.\ $x\in M$ that   $$F^*(Du_\lambda(x))=\frac{\lambda p'}{p} (J_{1,\lambda}^{x_0})^{-\frac{1}{p}}e^{-\frac{\lambda}{p}d_F(x_0,x)^{p'}}d_F(x_0,x)^{p'-1},$$
	therefore, 
$$	\|F^*(Du_\lambda)\|_{L^p(M,{\sf m})}^p=\left(\frac{\lambda p'}{p}\right)^p (J_{1,\lambda}^{x_0})^{-1} J_{2,\lambda}^{x_0}.$$
Replacing $u_\lambda$ into \eqref{LSI-contra}, a minor reorganization of the terms provide
$$-\lambda (J_{1,\lambda}^{x_0})^{-1} J_{2,\lambda}^{x_0}\leq \frac{N}{p}\log\left(C\, \left(\frac{\lambda p'}{p}\right)^p (J_{1,\lambda}^{x_0})^{-1+\frac{p}{N}} J_{2,\lambda}^{x_0} \right).$$
Taking the limit $\lambda\to 0$ and using   \eqref{limitek-J}, it follows that $ \mathcal
L_{p,N} ({\sf AVR}_{\mathcal M}^+) ^{-\frac{p}{N}}\leq C$, a contradiction. 

	\subsection{Equality in \eqref{LSI}} 
	We first recall the $L^p$-logarithmic Sobolev inequality on the model space $\mathcal M_N=(\mathbb R_+,\vert\cdot\vert,{\sf m}_N=N\omega_Nr^{N-1}{\mathcal L ^1}_{\llcorner{[0,\infty)}})$, i.e., for every   function $v:\mathbb R_+\to \mathbb R$  with $v\in W^p_p(\mathcal M_N)$ and $\|v\|_{L^p(\mathbb R_+,{\sf m}_N)}=1$ one has that   
	\begin{equation}\label{LSI-1dim-ineq}
	\int_{\mathbb R_+}|v|^p\log |v|^p{\rm d} {\sf m}_N\leq \frac{N}{p}\log\left(\mathcal
	L_{p,N}\, \|v'\|_{L^p(\mathbb R_+,{\sf m}_N)}^p \right).
		\end{equation}
	Inequality \eqref{LSI-1dim-ineq} is known by del Pino--Dolbeault  \cite{Dol-delPino} for $n\in \mathbb N$ and $p\in (1,n)$, and by Balogh--Don--Krist\'aly \cite{BDK} for the general case $N,p>1$.  
	Moreover, according to \cite{BDK}, equality holds in \eqref{LSI-1dim-ineq} for some $v\in W^p_p(\mathcal M_N)$ with $\|v\|_{L^p(\mathbb R_+,{\sf m}_N)}=1$ if and only if   $v(r)=\pm \lambda^\frac{N}{pp'}\left(\omega_N    \Gamma \left(\frac{N}{p'}+1\right)\right)^{-1/p}e^{-\frac{\lambda}{p}r^{p'}}$ for some $\lambda>0.$

	Let us assume that equality holds in \eqref{LSI}  for some   $u\in   W^p_{p}(\mathcal M)$  with $\|u\|_{L^p(M,{\sf m})}=1,$  i.e., 
	\begin{equation}\label{LSI-equil} 	\int_{M}|u|^p\log |u|^p{\rm d} {\sf m}=\frac{N}{p}\log\left(\mathcal
		L_{p,N}{\sf AVR}_ {\mathcal M}^{-\frac{p}{N}}\|F^*(Du)\|_{L^p(M,{\sf m})}^p \right). 
	\end{equation}
Let $u_+=\max\{u,0\} \geq 0$ and $u_-=-\min\{u,0\}\geq 0$ be the positive and negative parts of $u$; it is clear that both relations \eqref{pos-neg-splitting-F-log} and 
\eqref{pos-neg-splitting-u-log} hold. The backward/forward   P\'olya--Szeg\H o inequalities \eqref{Polya-Szego} and \eqref{pos-neg-splitting-F-log} imply that   
	\begin{align}\label{standard-PSZ}
	\nonumber \|F^*(Du)\|_{L^p(M,{\sf m})}^p&=\|F^*(Du_+)\|_{L^p(M,{\sf m})}^p+\|F^*(-Du_-)\|_{L^p(M,{\sf m})}^p\\&\nonumber =\|\overleftarrow F^*(-Du_+)\|_{L^p(M,{\sf m})}^p+\|F^*(-Du_-)\|_{L^p(M,{\sf m})}^p\\&\nonumber\geq  ( {\sf AVR}^-_{\mathcal M})^{p/N}\|(u_+^*)'\|^p_{L^p(\mathbb R_+,{\sf m}_N)} + ( {\sf AVR}^+_{\mathcal M})^{p/N}\left(\|(u_-^*)'\|^p_{L^p(\mathbb R_+,{\sf m}_N)}\right)\\&\geq {\sf AVR}_{\mathcal M}^{p/N} \left( \|(u_+^*)'\|^p_{L^{p }(\mathbb R_+,{\sf m}_N)} + \|(u_-^*)'\|^p_{L^{p}(\mathbb R_+,{\sf m}_N)}\right)
	.
\end{align}

Assume, by contradiction, that  $U_+:=\|u_+\|_{L^p(M,{\sf m})}=\|u_+^*\|_{L^p(\mathbb R_+,{\sf m}_N)}>0$ and $U_-:=\|u_-\|_{L^p(M,{\sf m})}=\|u_-^*\|_{L^p(\mathbb R_+,{\sf m}_N)}>0$. By the strict convexity of 
$s\mapsto s\log s$, $s\geq 0$, the strict concavity of 
$s\mapsto  \log s$, $s\geq 0$, relations \eqref{pos-neg-splitting-u-log},  \eqref{LSI-1dim-ineq} and \eqref{standard-PSZ}, the Cavalieri principle \eqref{Cavallieri} and  $U_+^p+U_-^p=1$, it follows, similarly to \eqref{bal-oldal}, that 
\begin{align}\label{bal-oldal-logg}
	\nonumber	\int_M |u|^p \log |u|^p{\rm d} {\sf m} &=\int_M(u_+^p+u_-^p)\log(u_+^p+u_-^p){\rm d} {\sf m}\\&\nonumber<  \int_{M}\left(  {u_+ ^p}\log  \frac{u_+ ^p}{U_+^p}+  {u_- ^p}\log  \frac{u_- ^p}{U_-^p}\right){\rm d} {\sf m}\\&\nonumber =\int_{\mathbb R_+}\left(  {(u_+^*) ^p}\log  \frac{(u_+^*) ^p}{U_+^p}+  {(u_-^*) ^p}\log  \frac{(u_-^*) ^p}{U_-^p}\right){\rm d} {\sf m}_N\\&\nonumber  \leq U_+^p \frac{N}{p}\log\left(\mathcal
	L_{p,N}\, \frac{\|(u_+^*)'\|_{L^p(\mathbb R_+,{\sf m}_N)}^p}{U_+^p} \right)+ U_-^p \frac{N}{p}\log\left(\mathcal
	L_{p,N}\, \frac{\|(u_-^*)'\|_{L^p(\mathbb R_+,{\sf m}_N)}^p}{U_-^p} \right)\\& \leq \nonumber \frac{N}{p}\log\left(\mathcal
	L_{p,N}\, \left( \|(u_+^*)'\|^p_{L^{p }(\mathbb R_+,{\sf m}_N)} + \|(u_-^*)'\|^p_{L^{p}(\mathbb R_+,{\sf m}_N)}\right) \right)\\&\nonumber \leq \frac{N}{p}\log\left(\mathcal
	L_{p,N}{\sf AVR}_ {\mathcal M}^{-\frac{p}{N}}\|F^*(Du)\|_{L^p(M,{\sf m})}^p \right),
\end{align}
which contradicts \eqref{LSI-equil}.

Thus, either $U_+=0$ or  $U_-=0$. In particular, the extremizer $u$ is sign-preserving. Without loss of generality, we may assume that $U_+=0$, i.e., $u=-u_-\leq 0$; the other case can be discussed analogously. In this way, every estimate in the above chain of inequality should be equality. Thus, by the equality in the forward P\'olya--Szeg\H o inequality \eqref{standard-PSZ}, i.e., 
$$	\|F^*(-Du_-)\|_{L^p(M,{\sf m})}=   ( {\sf AVR}^+_{\mathcal M})^{1/N}\|(u_-^*)'\|_{L^p(\mathbb R_+,{\sf m}_N)},$$ jointly with 
  ${\sf AVR}_{\mathcal M}={\sf AVR}_{\mathcal M}^+$, implies the existence of  $x_0\in M$ such that  $\mathcal M=(M,F,{\sf m})$  is a  forward  $N$-volume cone  with  pole $x_0$,  see  \S\ref{subsection-rearr}, which means that  
$
\frac{ {\sf m}(B_{x_0}^+(R))}{\omega_NR^N}={\sf AVR}_{\mathcal M}^+$ for every $ R>0.$  Proposition \ref{N-volume-cone-prop} implies now the first part of our conclusion. In addition,  
since $u_-^*$  produces equality in \eqref{LSI-1dim-ineq}, 
it turns out that    $u^*_-(r)=\mu^\frac{N}{pp'}\left(\omega_N    \Gamma \left(\frac{N}{p'}+1\right)\right)^{-1/p}e^{-\frac{\mu }{p}r^{p'}}$ for some $\mu>0.$
According to \eqref{u=u^*},  
$u_-(x)=u_-^*\left({\sf AVR}^+_{\mathcal M})^{1/N} d_F(x_0,x)\right) $  for ${\sf m}{\rm -a.e.}\  x\in M$,  which concludes the proof, recalling that $u=-u_-$. Clearly, when $U_-=0$, we proceed in a similar way. 

The converse is a simple computation. \hfill $\square$

\section{Optimal hypercontractivity  for Finslerian heat flow: proof of Theorem \ref{hypercontractivity-main-theorem}} \label{section-5}

Let $\mathcal M=(M,F,{\sf m})$ be 
 a not necessarily reversible $n$-dimensional  complete Finsler space  and consider  the heat equation   
 \begin{equation}\label{heat-flow-1}
	\partial_t u= \bm \Delta_{\mathcal M} u \quad {\rm on}\ [0,T]\times M,   
\end{equation} 
for some $T>0$. 
 According to Ohta--Sturm \cite[Theorem 3.4 \& p.\ 1402]{Ohta-Sturm} and Xia \cite{Xia1}, for every non-negative initial data $u(0,\cdot)=u_0\in L^2(M,{\sf m})$, there exists a non-negative global solution $u_t(x)=u(t,x)$ to the heat equation \eqref{heat-flow-1} 
 belonging to $L^2([0,T],H_0^1(M))\cap H^1([0,T],H^{-1}(M))$, 
and  $t\mapsto \mathcal E(u_t)$ is non-increasing, where 
  $$\mathcal E(u_t)=\frac{1}{2}\int_M F^*(Du_t)^2{\rm d}{\sf m}$$
  is the energy   of $u_t$. Note that when $\partial M\neq \emptyset$, the space $H_0^1(M)$ is the closure of compactly supported functions in $M$, not necessarily vanishing on $\partial M.$   Since $\partial _t \|u_t\|^2_{L^2(M,{\sf m})}=-4 \mathcal E(u_t)$, see Ohta--Sturm \cite[Proposition 3.5]{Ohta-Sturm}, it follows for every $0<t\leq T$ that 
  $$t \mathcal E(u_t)\leq \int_0^t \mathcal E(u_s){\rm d}s=\frac{1}{4}(\|u_0\|^2_{L^2(M,{\sf m})}-\|u_t\|^2_{L^2(M,{\sf m})}) \leq \frac{1}{4}\|u_0\|^2_{L^2(M,{\sf m})}.$$ In particular, 
  $u_t\in W_2^2(\mathcal M)$ for every $0<t\leq T$.  In addition, if $u_0\in L^1(M,{\sf m})\cap L^\infty(M)$, the same holds for $u_t$, and   $u_t^\alpha \in W_2^2(\mathcal M)$ for every $\alpha\geq 1$ and $t\in (0,T]$.  
 
\subsection{Proof of \eqref{hyper-estimate}}
In the sequel, we shall use   Theorem \ref{log-Sobolev-main}  for $p=2$, for the reverse metric $\overleftarrow F$ and for  non-negative functions; namely, in   non-normalized form, one has for every non-negative function $u\in W^2_2(\mathcal M)\setminus \{0\}$ that 
\begin{equation}\label{log-Sobolev}
	\frac{{{\sf Ent}_{\sf m}}(u^2)}{\|u\|^2_{L^2(M,{\sf m})}}\le \frac{N}{2}\log\left(\frac{2}{N\pi e}({\sf AVR}_{\mathcal M}^{-})^{-\frac{2}{N}}\frac{\ds\int_M \overleftarrow F^*(-D{u})^2  \, {\rm d}{\sf m}}{\|u\|^2_{L^2(M,{\sf m})}}\right), 
\end{equation}
where 
$${\sf Ent}_{\sf m}(u^2)=\int_M u^2\log u^2\, {\rm d}{\sf m} -\left(\int_M u^2\, {\rm d}{\sf m}\right)\log \left(\int_M u^2\, {\rm d}{\sf m}\right) $$
is the entropy of $u^2$ with respect to the measure ${\sf m}$. We recall that \eqref{log-Sobolev} is optimal and equality holds for some $u\in W^2_2(\mathcal M)\setminus \{0\}$ if and only if there exists $x_0\in M$ $($or $x_0\in \partial M$ if $M$ has boundary$)$  such that the statement -- concerning the flag and $S$-curvatures -- from Theorem \ref{Thm-Sobolev-2-equality} {\rm (ii)/(a)} holds  and $u(x)=c e^{-\frac{\lambda}{2} d_F(x,x_0)^{2}}$, $x\in M,$ for some $\lambda,c>0.$
 
 Let $u_t(x)=u(t,x)$ be the non-negative global solution  to the heat equation \eqref{heat-flow-1} on $[0,T]\times M$ with initial data $u_0\in   L^1(M,{\sf m})\cap L^\infty(M)$, $u_0\geq 0$. We fix 
$2\leq p< q\leq +\infty$ and $0\leq T_0<  T_1\leq T$; the case $p=q$ will be obtained by limiting. Let us consider an increasing  function $r:[T_0,T_1]\to [2,\infty)$ such that $r(T_0)=p$ and $r(T_1)=q$.

By the properties of $u_t$,   $t\mapsto \|u_t\|^{r(t)}_{L^{r(t)}(M,{\sf m})}$ is differentiable a.e.\ on $(T_0,T_1)$ and according to the heat equation  \eqref{heat-flow-1} and relation \eqref{Laplacian}, for a.e.\ $t\in (T_0,T_1)$ one has that
	\begin{align*}
		\frac{\rm d}{{\rm d}t}\|u_t\|^{r(t)}_{L^{r(t)}(M,{\sf m})}=&\frac{\rm d}{{\rm d}t}\int_M u_t^{r(t)}{\rm d}{\sf m}=\int_M u_t^{r(t)-1}\left(r'(t)u_t\log u_t+r(t)\partial_tu_t \right){\rm d}{\sf m}\\=&\int_M u_t^{r(t)-1}\left(r'(t)u_t\log u_t+r(t)\bm \Delta_{\mathcal M} u_t\right){\rm d}{\sf m}\\=&
		\frac{	r'(t)}{r(t)}\int_M u_t^{r(t)}\log u_t^{r(t)}{\rm d}{\sf m} -r(t)(r(t)-1)\int_M u_t^{r(t)-2}F^*({D u_t})^2 \, {\rm d}{\sf m}.
	\end{align*}
	Therefore, it turns out that 
	\begin{equation}\label{derivative-1}
		\frac{\rm d}{{\rm d}t}\log \|u_t\|_{L^{r(t)}(M,{\sf m})}=\frac{r'(t)}{r^2(t)}\frac{	{\sf Ent}_{\sf m}\left(u_t^{r(t)}\right)}{\|{u_t}\|_{L^{r(t)}(M,{\sf m})}^{  r(t)} }-4\frac{ r(t)-1 }{r^2(t)}\frac{\displaystyle\int_M     F^*\left({D u_t^{r(t)/2}}\right)^2 \, {\rm d}{\sf m}}{\|{u_t}\|_{L^{r(t)}(M,{\sf m})}^{ r(t)} } .
	\end{equation}
	Due to the above reasoning, since $r\geq 2$, the function $u:=u_t^{r(t)/2}$ belongs to $W_2^2(\mathcal M)\setminus \{0\}$ for every $t\in (0,T].$ Thus, we may apply the optimal $L^2$-logarithmic Sobolev inequality \eqref{log-Sobolev}, obtaining jointly with \eqref{derivative-1} that for a.e.\ $t\in (T_0,T_1)$, 
		\begin{align}\label{derivative-2}
	\nonumber	\frac{\rm d}{{\rm d}t}\log \|u_t\|_{L^{r(t)}(M,{\sf m})}\leq& \frac{N}{2}\frac{r'(t)}{r^2(t)}\log\left(\frac{2}{N\pi e}({\sf AVR}_{\mathcal M}^{-})^{-\frac{2}{N}}\frac{\displaystyle\int_M     F^*\left({D u_t^{r(t)/2}}\right)^2 \, {\rm d}{\sf m}}{\|{u_t}\|_{L^{r(t)}(M,{\sf m})}^{ r(t)} }\right)\\&-4\frac{ r(t)-1 }{r^2(t)}\frac{\displaystyle\int_M     F^*\left({D u_t^{r(t)/2}}\right)^2 \, {\rm d}{\sf m}}{\|{u_t}\|_{L^{r(t)}(M,{\sf m})}^{ r(t)} } .
	\end{align}
Given the numbers $A_i>0,$ $i\in \{1,2,3\}$, we have the elementary inequality 
	\begin{equation}\label{elemi-ineq}
		A_1 \log(A_2\,s )-A_3\,s\leq A_1\log(\frac{A_1A_2}{eA_3}),\quad s>0,
	\end{equation}
with equality if and only if $s=\frac{A_1}{A_3}$. Thus, the natural identification of the coefficients  from \eqref{derivative-2} implies  for a.e.\ $t\in (T_0,T_1)$ that 
	\begin{align}\label{derivative-3}
	\nonumber	\frac{\rm d}{{\rm d}t}\log \|u_t\|_{L^{r(t)}(M,{\sf m})}\leq& \frac{N}{2}\frac{r'(t)}{r^2(t)}\log\left(\frac{({\sf AVR}_{\mathcal M}^{-})^{-\frac{2}{N}}}{4\pi e^2}\frac{r'(t)}{r(t)-1}\right).
\end{align}
Integrating the latter inequality on $(T_0,T_1)$, it yields that 
$$	 	 \log \frac{\|u_{T_1}\|_{L^{r(T_1)}(M,{\sf m})}}{\|u_{T_0}\|_{L^{r(T_0)}(M,{\sf m})}}
\leq  \frac{N}{2}\int_{T_0}^{T_1}\frac{r'(t)}{r^2(t)}\log\left(\frac{({\sf AVR}_{\mathcal M}^{-})^{-\frac{2}{N}}}{4\pi e^2}\frac{r'(t)}{r(t)-1}\right){\rm d}t.$$
As we want to have minimal value for the right hand side with the conditions $r(T_0)=p$ and $r(T_1)=q$, the Euler--Lagrange equation gives the first order ordinary differential equation $-\frac{r'(t)}{r^2(t)}=C$ for $t\in  [T_0,T_1]$ and some $C\in \mathbb R$, which has the solution 
\begin{equation}\label{r-function}
	r(t)=\frac{pq(T_0-T_1)}{(q-p)t+pT_0-qT_1},\quad t\in [T_0,T_1].
\end{equation}
Replacing this expression into the latter estimate, an elementary computation yields that 
$$\|u_{T_1}\|_{L^q(M,{\sf m})}\leq \frac{{\sf M}(p,q)^\frac{N}{2}({\sf AVR}_{\mathcal M}^-)^{\frac{1}{q}-\frac{1}{p}}}{(4\pi (T_1-T_0))^{\frac{N}{2}(\frac{1}{p}-\frac{1}{q})}}\|u_{T_0}\|_{L^p(M,{\sf m})},$$
where ${\sf M}(p,q)$ is from \eqref{M-p-q}, which concludes the proof of
\eqref{hyper-estimate}. 

\subsection{Optimality in \eqref{hyper-estimate}} Assume by contradiction that there exists $C>{\sf AVR}_{\mathcal M}^-$ such that 
\begin{equation}\label{contradiction-hyper}
	\|u(T_1,\cdot)\|_{L^q(M,{\sf m})}\leq \frac{{\sf M}(p,q)^\frac{N}{2}C^{\frac{1}{q}-\frac{1}{p}}}{(4\pi (T_1-T_0))^{\frac{N}{2}(\frac{1}{p}-\frac{1}{q})}}\|u(T_0,\cdot)\|_{L^p(M,{\sf m})}
\end{equation}  
for every   $2\leq p\leq q\leq +\infty$ and $0\leq T_0<  T_1\leq T$. Let us choose $T_0=0$, $T_1=t<T$, $p=2$ and $q:=q(t)=2+yt$, where $y>0$ will be fixed later; in addition, we may assume that $u_0=u(0,\cdot)\in W_2^2(\mathcal M)$.  Relation \eqref{contradiction-hyper} implies that 
$$\lim_{t\to 0}\frac{\log \|u(t,\cdot)\|_{L^{q(t)}(M,{\sf m})}-\log \|u(0,\cdot)\|_{L^2(M,{\sf m})}}{t}\leq \lim_{t\to 0}\log \left(\frac{{\sf M}(2,q(t))^\frac{N}{2}C^{\frac{1}{q(t)}-\frac{1}{2}}}{(4\pi t)^{\frac{N}{2}(\frac{1}{2}-\frac{1}{q(t)})}}\right)^{1/t}.$$
Note that a similar argument as in \eqref{derivative-1} and a straightforward limiting imply that
 $$\frac{y}{4}\frac{	{\sf Ent}_{\sf m}\left(u_0^{2}\right)}{\|{u_0}\|_{L^2(M,{\sf m})}^{  2} }- \frac{\displaystyle\int_M     F^*\left({D u_0}\right)^2 \, {\rm d}{\sf m}}{\|{u_0}\|_{L^{2}(M,{\sf m})}^{2} } \leq \frac{y}{4}\left(\frac{N}{2}\log(\frac{y}{4\pi e^2})-\log C\right).$$
 Dividing by $y/4>0$ the latter inequality and taking the minimum with respect to $y$, which has the value $$y:=\frac{8}{N}\frac{\displaystyle\int_M     F^*\left({D u_0}\right)^2 \, {\rm d}{\sf m}}{\|{u_0}\|_{L^{2}(M,{\sf m})}^{2} } ,$$
 we obtain that 
$$	\frac{{{\sf Ent}_{\sf m}}(u_0^2)}{\|u_0\|^2_{L^2(M,{\sf m})}}\le \frac{N}{2}\log\left(\frac{2}{N\pi e}C^{-\frac{2}{N}}\frac{\ds\int_M  F^*(D{u_0})^2  \, {\rm d}{\sf m}}{\|u_0\|^2_{L^2(M,{\sf m})}}\right). $$
Since $C>{\sf AVR}_{\mathcal M}^-$, the latter inequality contradicts the optimality of the $L^2$-logarithmic Sobolev inequality   \eqref{log-Sobolev}.

\subsection{Equality  in \eqref{hyper-estimate}} (i) $\implies$ (ii) Assume that equality holds in \eqref{hyper-estimate}  for some $ p< q$,  $ T_0<  T_1\leq T $ and positive function $u\in AC([T_0,T_1], L^2(M,{\sf m}))$. Tracking back the proof of \eqref{hyper-estimate}, we should have equality in \eqref{derivative-2}, i.e., in the optimal $L^2$-logarithmic Sobolev inequality. Thus, according to Theorem \ref{log-Sobolev-main}, for a.e.\ $t\in [T_0,T_1]$   there exists $x_0^t\in M$ $($or $x_0^t\in \partial M$ if $M$ has boundary$)$  such that the statement of Theorem \ref{Thm-Sobolev-2-equality} {\rm (ii)/(a)} holds (for $x_0^t\in M$) and 
\begin{equation}\label{u-t-formaja}
	u_t(x)^{r(t)/2}=c_t e^{-\frac{\lambda_t}{2}d_F(x,x_0^t)^2},\quad {\sf m}-{\rm a.e.}\    x\in M,
\end{equation}
for some $c_t,\lambda_t>0;$ the function $r$ is from \eqref{r-function}. Moreover, equality should also hold in \eqref{elemi-ineq} with the appropriate choices, i.e., 
\begin{equation}\label{F-equal}
	\frac{\displaystyle\int_M     F^*\left({D u_t^{r(t)/2}}\right)^2 \, {\rm d}{\sf m}}{\|{u_t}\|_{L^{r(t)}(M,{\sf m})}^{ r(t)}}=\frac{N}{8}\frac{r'(t)}{r(t)-1} \quad {\rm for\ a.e.}\ t\in [T_0,T_1].
\end{equation}
Since $\mathcal M=(M,F,{\sf m})$ is a   backward  $N$-volume cone with pole $x_0^t$ -- see the statement of Theorem \ref{Thm-Sobolev-2-equality} {\rm (ii)/(a)} and  Proposition \ref{N-volume-cone-prop} for the backward case -- one has that  ${\sf m}(B_{x_0^t}^-(R))={\sf AVR}_{\mathcal M}^- \omega_N R^N$  for every $R>0$ and a.e.\ $t\in [T_0,T_1]$. Thus, the eikonal equation \eqref{eikonal} and \eqref{u-t-formaja} imply that for a.e.\ $t\in [T_0,T_1]$, 
\begin{equation}\label{ur-t-dif}
		\|{u_t}\|_{L^{r(t)}(M,{\sf m})}^{ r(t)}=c_t^2\lambda_t^{-\frac{N}{2}} \pi^{\frac{N}{2}} {\sf AVR}_{\mathcal M}^-   \ \ {\rm and}\ \  \displaystyle\int_M     F^*\left({D u_t^{r(t)/2}}\right)^2 \, {\rm d}{\sf m}=\frac{N}{2}c_t^2\lambda_t^{-\frac{N}{2}+1} \pi^{\frac{N}{2}} {\sf AVR}_{\mathcal M}^- .
\end{equation}
By \eqref{F-equal}, it follows that   for  a.e.\  $t\in [T_0,T_1],$
$$\lambda_t=\frac{1}{4}\frac{r'(t)}{r(t)-1}=\frac{r(t)}{4t+K_0},$$
where  $$K_0=4\frac{(p-1)qT_1-(q-1)pT_0}{q-p}.$$ 
Note that $K_0>-4T_0.$ Clearly, $\lambda_t$ can be continuously (in fact, smoothly) extended to the whole interval $[T_0,T_1]$. On the other hand, since $t\mapsto \|u_t\|^{r(t)}_{L^{r(t)}(M,{\sf m})}$ is absolutely continuous, by the first relation of \eqref{ur-t-dif}, the function $t\mapsto c_t$ is absolutely continuous on $[T_0,T_1]$. Furthermore, since $u\in AC([T_0,T_1], L^2(M,{\sf m}))$, by \eqref{u-t-formaja}, the function $t\mapsto d_F(x,x_0^t)^2$ is absolutely  continuous on   $[T_0,T_1]$ for {\sf m}-a.e.\ $x\in M$.  
Therefore, the latter regularity properties imply that we have a continuous representation in \eqref{u-t-formaja}, i.e.,  
\begin{equation}\label{u_t-continuous}
	u_t(x)=C_te^{-\frac{d_F(x,x_0^t)^2}{4t+K_0}} \quad {\rm for\ every}\ t\in [T_0,T_1],\ x\in M,
\end{equation}
where $C_t= c_t^{2/r(t)}>0.$

For simplicity, let $\rho_t(x)=d_F(x,x_0^t)^2=d_{\overleftarrow F}( x_0^t,x)^2.$ Since $\mathcal M=(M,F,{\sf m})$ is a   backward  $N$-volume cone with pole $x_0^t$, it turns out that on $M\setminus \{x_0^t\}$, we have $$\Delta_{\overleftarrow {\mathcal M}}\rho_t = 2N$$ in the distributional sense, and $\Delta_{  {\mathcal M}} (-d_{\overleftarrow F})=-\Delta_{\overleftarrow {\mathcal M}} d_{\overleftarrow F}$. These relations imply that 
\begin{equation}\label{Delta-pointwise}
	\Delta_{  {\mathcal M}} u_t=u_t\left(\frac{4\rho_t}{(4t+K_0)^2}-\frac{2N}{4t+K_0}\right);
\end{equation}
although this relation is understood in the  distributional   sense on $M$, the representation of $u_t$ from \eqref{u_t-continuous}, combined with the  above regularity properties allow us to identify its corresponding representative pointwise. 
Since $u_t>0$, the $\partial_t \log u_t =\frac{\partial_t   u_t}{u_t}$ can be interpreted locally in $M$, with the  distributional time derivative. Accordingly, by \eqref{Delta-pointwise} and the heat equation  \eqref{heat-flow-1}, it follows that
\begin{equation}\label{ut-tutww}
 \partial_t \log u_t = \frac{4\rho_t}{(4t+K_0)^2}-\frac{2N}{4t+K_0};
\end{equation}
thus for every $x,y\in M$, 
\begin{equation}\label{utttt}
	\partial_t \log \frac{u_t(x)}{u_t(y)}=\frac{4(\rho_t(x)-\rho_t(y))}{(4t+K_0)^2}.
\end{equation}
On the other hand, differentiating relation \eqref{u_t-continuous} in the distributional sense  yields that 
$$\partial_t \log \frac{u_t(x)}{u_t(y)}=\frac{4(\rho_t(x)-\rho_t(y))}{(4t+K_0)^2}-\frac{\partial_t(\rho_t(x)-\rho_t(y))}{ 4t+K_0 }.$$ The latter relation and \eqref{utttt} imply that $\partial_t(\rho_t(x)-\rho_t(y))=0$, thus $$\rho_t(x)-\rho_t(y)=\rho_s(x)-\rho_s(y)$$ 
for every $x,y\in M$ and $t,s\in (T_0,T_1).$ If we insert $x:=x_0^t$ and $y:=x_0^s$ in the latter relation, it follows that $\rho_t(x_0^s)+\rho_s(x_0^t)=0$, thus   $x_0^s=x_0^t$. In particular, $x_0^t=x_0$ for every $t\in [T_0,T_1],$ thus $\rho_t(x)=d_F(x,x_0^t)^2$ is time-independent and  \eqref{u_t-continuous} reduces to 
\begin{equation}\label{u_t-continuous-1}
	u_t(x)=C_te^{-\frac{d_F(x,x_0)^2}{4t+K_0}} \quad {\rm for\ every}\ t\in [T_0,T_1],\ x\in M.
\end{equation}
Now, \eqref{ut-tutww} and \eqref{u_t-continuous-1} imply that $\frac{C_t'}{C_t}=-\frac{2N}{4t+K_0}$  on $(T_0,T_1)$, thus for some $\lambda>0$, one has that $C_t=\lambda (4t+K_0)^{-\frac{N}{2}}$  for every $t\in [T_0,T_1]$, which provides the form of $u(t,x)$ from \eqref{u-t-x-precise}.  

(ii) $\implies$ (i)  By (ii)/(a) and Proposition \ref{N-volume-cone-prop}, it turns out that $\mathcal M=(M,F,{\sf m})$ is a   backward  $N$-volume cone with pole $x_0$, i.e.,   ${\sf m}(B_{x_0}^-(R))={\sf AVR}_{\mathcal M}^- \omega_N R^N$  for every $R>0$. Therefore, the precise ansatz for $u_t(x)=u(t,x)$ from \eqref{u-t-x-precise} and the layer cake representation imply that for every $t\in [T_0,T_1]$ and $r\geq 1$, 
$$\|{u_t}\|_{L^{r }(M,{\sf m})} =\lambda (4t+K_0)^{\frac{N}{2}(\frac{1}{r}-1)}  \left(    \frac{\pi}{r}\right)^{\frac{N}{2r}} ({\sf AVR}_{\mathcal M}^-)^\frac{1}{r}.$$
Having the latter expression together with 
$$4T_0+K_0=\frac{4q(p-1)}{q-p}(T_1-T_0)\ \ {\rm and}\ \ 4T_1+K_0=\frac{4p(q-1)}{q-p}(T_1-T_0),$$
it remains to directly check the equality in \eqref{hyper-estimate}. \hfill $\square$

		\section{Examples} \label{section-6}
		
			\subsection{Zermelo navigation problems}  We consider the space $\mathcal M_g=(M,g,{{\rm vol}_g})$, where $(M,g)$ is an $n$-dimensional complete Riemannian manifold (without boundary) with ${\sf Ric}_g\geq 0$ and ${\sf AVR}_{\mathcal M_g}>0$, and  ${{\rm vol}_g}$ is the canonical measure on $(M,g)$. Let $W$ be a Killing vector field of $g$ with the flow $\psi_t$, i.e., $\psi_t^* g =g$, and let $$a:=\sup_M|W|_g<1;$$ the vector field $W$ represents usually the wind or a current.  
			
			The \textit{Zermelo navigation metric}  $F:TM\to [0,\infty)$  associated to the pair $(g,W)$ on $M$ is given by  $$|y-F(x,y)W(x)|_g=F(x,y),\quad (x,y)\in TM,$$ i.e.,  the unit sphere with respect to $F$ is the unit sphere with respect to $g$ translated by the vector field $W$, see 		  Bao--Robles--Shen \cite{BRS}. The explicit form of $F$ is  the  Randers metric $F=\alpha+\beta$, given by $$\alpha=\frac{\sqrt{(1-|W|_g^2)g(y,y)+g(W,y)^2}}{1-|W|_g^2},\quad \beta=-\frac{g(W,y)}{1-|W|_g^2}.$$
	Having the navigation metric, the Finsler space $\mathcal M=(M,F,{\sf m}_{BH})$ is non-reversible -- unless $W=0$, when $F=|\cdot|_g$ reduces to the Riemannian norm -- and its reversibility is $r_F=\frac{1+a}{1-a},$ while the Busemann--Hausdorff measure is ${\sf m}_{BH}={{\rm vol}_g}$. Since $r_F<\infty$, the space $\mathcal M$ is complete (both forward and backward complete).

	According to Foulon--Matveev \cite[Theorem 1]{FM}, 	every geodesic of $(M,F)$ is precisely a curve  $ \gamma_F(t):= (\psi_t\circ \gamma) (t)$  with $\gamma$ being a geodesic in $(M,g)$, while  the flag curvature ${\sf K}$ on $(M,F)$ and sectional curvature ${\sf K}_g$ on $(M,g)$ are related by $${\sf K}_g(v,w)=\mathbf{\sf K}(v+|v|_gW;w)$$ for every  $v,w\in T_xM$ and $x\in M$ with $g(v,w)=0$, see also Huang--Krist\'aly--Zhao \cite{HKZ}; note that $v+|v|_gW$ and $w$ are linearly independent, see Javaloyes--Vit\'orio \cite{JV}. In particular,  ${\sf Ric}(\dot \gamma_F)={\sf Ric}_g(\dot \gamma)\geq 0.$  Moreover, since $W$ is a Killing vector field of $g$, it turns out that the $S$-curvature of $(M,F)$ vanishes with respect to ${\sf m}_{BH}$, i.e., ${\sf S}_{{\sf m}_{BH}}\equiv 0$, see Shen \cite[p.\ 325]{Shen-kotet}. 
%	
%	
%	
%	 $g(W,\dot\gamma)$ is constant along the geodesic $\gamma$ of $(M,g)$. Therefore,   the distortion is $$\tau(\gamma_F(t),\dot \gamma_F(t))=-\frac{n+1}{2}\log(1+g((W\circ \psi_t)(\gamma(t)),\dot\gamma_F(t)))=-\frac{n+1}{2}\log(1+g(W(\gamma(t)),\dot\gamma(t))),$$
%	which is constant, thus for every $(x,y)\in TM\setminus\{0\}$, one has that 
%	$${\sf S}_{{\sf m}_{BH}}(x,y) =\left.\frac{\dd}{{\dd}t}\right|_{t=0}\tau(\gamma_F(t),\dot \gamma_F(t))=0,$$
%	where $\gamma_F(=0)$ and $\dot \gamma_F(0)=y.$ 
	Summing up, it follows that
	$${\sf Ric}_n(\dot \gamma_F)={\sf Ric}(\dot \gamma_F)={\sf Ric}_g(\dot \gamma)\geq 0.$$

 Note that $d_F(x,y)$ is the unique root $t$ of the equation $d_g(\psi_t(x),y)=t$, thus $B_x^+(R)=\psi_R(B_x^g(R))$ for every $x \in M$ and $R>0$, where $B_x^g(R)=\{y\in M:d_g(x,y)<R\}$. In the same manner, we also have that $B_x^-(R)=\psi_{-R}(B_x^g(R))$. Consequently,  $${\sf m}_{BH}(B_x^\pm(R))={{\rm vol}_g}(\psi_{\pm R}(B_x^g(R)))={{\rm vol}_g}( B_x^g(R)) $$ for every $x\in M$ and $R>0$. Therefore,  
	$${\sf AVR}_{\mathcal M}^\pm=\lim_{R\to \infty}\frac{{\sf m}_{BH}(B_x^\pm(R))}{\omega_nR^n}=\lim_{R\to \infty}\frac{{{\rm vol}_g}( B_x^g(R))}{\omega_nR^n}={\sf AVR}_{\mathcal M_g}>0.$$
 Thus, Theorems \ref{Thm-Sobolev-1}-\ref{hypercontractivity-main-theorem} can be applied for the Zermelo navigation space  with Killing vector field. 
 
 Let us discuss the rigidity cases. Since $M$ is without boundary, it turns out that the pole $x_0\in M$ has the property that ${\sf m}_{BH}(B_{x_0}^\pm(R)) ={{\rm vol}_g}( B_{x_0}^g(R)) \sim \omega_n R^n$ as $R\to 0.$ Thus, the $n$-volume cone structure of the space gives that $1={\sf AVR}_{\mathcal M}^\pm={\sf AVR}_{\mathcal M_g}$; in particular,   ${{\rm vol}_g}( B_{x_0}^g(R))=\omega_n R^n$ for every $R>0$, which implies that $(M,g)$ is isometric to $\mathbb R^n$. On the other hand, it is easy to check that any Killing vector field  $W$ of $\mathbb R^n$ with $a:=\sup_{\mathbb R^n}|W|<1$ is constant. Accordingly, the existence of extremizer functions in Theorems \ref{Thm-Sobolev-1}-\ref{hypercontractivity-main-theorem} is supported only when the  Zermelo navigation metric $F$ is necessarily a  Randers--Minkowski metric of the form  $$F(x,y)=\frac{\sqrt{(1-|W|^2)|y|^2+\langle W,y\rangle^2}-\langle W,y\rangle}{1-|W|^2},$$
 where $W\in \mathbb R^n$ with $|W|<1$.

	\subsection{Convex cones in $\mathbb R^n$ with weights} \label{subsection-cone}
	    Let $M:=\Sigma\subseteq \mathbb R^n$ be an open convex cone and let $H:\mathbb R^n\to [0,\infty)$ be a Minkowski norm, i.e., $H\in C^\infty(\mathbb R^n\setminus\{0\})$ is   positive and strongly convex. 
%	    We consider the (not necessarily reversible) natural distance function on $\Sigma$ given by $d_H(x,y)=H(y-x)$. 
	    	    Let $w $ be a continuous function 
	    in $\overline \Sigma$, positive and locally Lipschitz in $\Sigma,$ and homogeneous of degree $\alpha$ such that   $w^{1/\alpha}$ is concave on $\Sigma$; when $\alpha=0$, $w$ is assumed to be constant. 
	     These properties imply that the Finsler space $\mathcal M=(\Sigma,H,w\mathcal L^n_{\llcorner{\Sigma}})$ is a complete, non-compact ${\sf CD}(0,N)$ space  with the synthetic dimension $N=n+\alpha$ and finite reversibility; in particular, its boundary (if $\Sigma\neq \mathbb R^n$) is locally convex and  one has that ${\sf Ric}_N\geq 0 $ on $\mathcal M$, see Cabr\'e--Ros-Oton--Serra \cite[Remark 1.4]{Cabre-JEMS}.   By the  homogeneity properties of $H$ and
	    $w,$ we have that
	    $${\sf AVR}_{\mathcal M}^+=\frac{\displaystyle \int_{\{x\in \Sigma:H(x)<1\} }w}{\omega_N}>0\ \ {\rm and}\  \ {\sf AVR}_{\mathcal M}^-=\frac{\displaystyle \int_{\{x\in \Sigma:H(-x)<1\} }w}{\omega_N}>0.$$
	    Accordingly, Theorems \ref{Thm-Sobolev-1}-\ref{hypercontractivity-main-theorem} apply for $\mathcal M=(  \Sigma,H,w\mathcal L^n_{\llcorner{\Sigma}})$. 
	    
	     In particular, the rigidity results in Theorems \ref{Thm-Sobolev-2-equality}--\ref{hypercontractivity-main-theorem} for various choices for $\mathcal M=(  \Sigma,H,w\mathcal L^n_{\llcorner{\Sigma}})$   provide  alternative proofs of the main results from Balogh--Guti\'errez--Krist\'aly \cite{BGK}, Balogh--Don--Krist\'aly \cite{BDK} and   Cabr\'e--Ros-Oton \cite{Cabre-JDE}; in the latter work the monomial weights are discussed.  Being a Minkowski space,   the flag curvature is identically zero on $\mathcal M=( \Sigma,H,w\mathcal L^n_{\llcorner{\Sigma}})$. The pole $x_0\in  \overline  \Sigma\cap (-\overline  \Sigma)$ -- which is not necessarily the origin --  has the property that 
	     \begin{equation}\label{w-TAMS}
	     	w(x_0+x)=w(x),\quad \forall x\in \Sigma,
	     \end{equation}
	      see e.g.\ \cite{BGK}. Moreover,   geodesics in $\Sigma$ emanating from $x_0$ are semi-lines  $\gamma(t)= x_0+te$, $t>0,$ with $e\in \Sigma\setminus\{0\} $. If we have the measure representation $w(\gamma(t))=e^{-\Psi_\gamma(t)}$ along the geodesic $\gamma$, then by \eqref{w-TAMS} and Euler's relation $\alpha w(y)=\langle \nabla w(y),y\rangle$ for every $y\in \Sigma$, the  $S$-curvature  is 
	     $${\sf S}_w(\gamma(t),\dot \gamma(t))=\Psi_\gamma'(t)=-\frac{\langle\nabla w(\gamma(t)), \dot\gamma(t)\rangle}{w(\gamma(t))}=-\frac{\langle\nabla w(te), e\rangle}{w(te)}=-\frac{\alpha }{t}=\frac{n-N}{t},$$
	     which agrees with the statements concerning the rigidity characterizations in Theorems \ref{Thm-Sobolev-2-equality}--\ref{hypercontractivity-main-theorem}.   
	    
%	     	\begin{wrapfigure}[5]{r}{0.45\textwidth}
%	     	\begin{center}
%	     		\includegraphics[width=0.40\textwidth]{Minkowski-figure.jpg}	\caption{Birds}
%	     	\end{center}
%	     
%	     \end{wrapfigure}

	    We notice that usually ${\sf AVR}_{\mathcal M}^+\neq {\sf AVR}_{\mathcal M}^-$. To see this, we consider a   particular case of the above setting; namely, let the open convex cone $$\Sigma=\{x=(x_1,x_2)\in \mathbb R^2: x_1>0 \}\subset \mathbb R^2,$$   the Minkowski/Randers norm $H:\mathbb R^2\to [0,\infty)$ given by
	    \begin{equation}\label{H-norm}
	    	   H(x_1,x_2)=\sqrt{x_1^2+x_2^2}-\frac{x_1}{2},
	    \end{equation}
	  	    and the monomial weight   $w(x_1,x_2)=x_1$ on $\overline \Sigma.$  
	  	    
	  	       \begin{figure}[h]\label{figure-1}
	  	    	\includegraphics[width=0.40\textwidth]{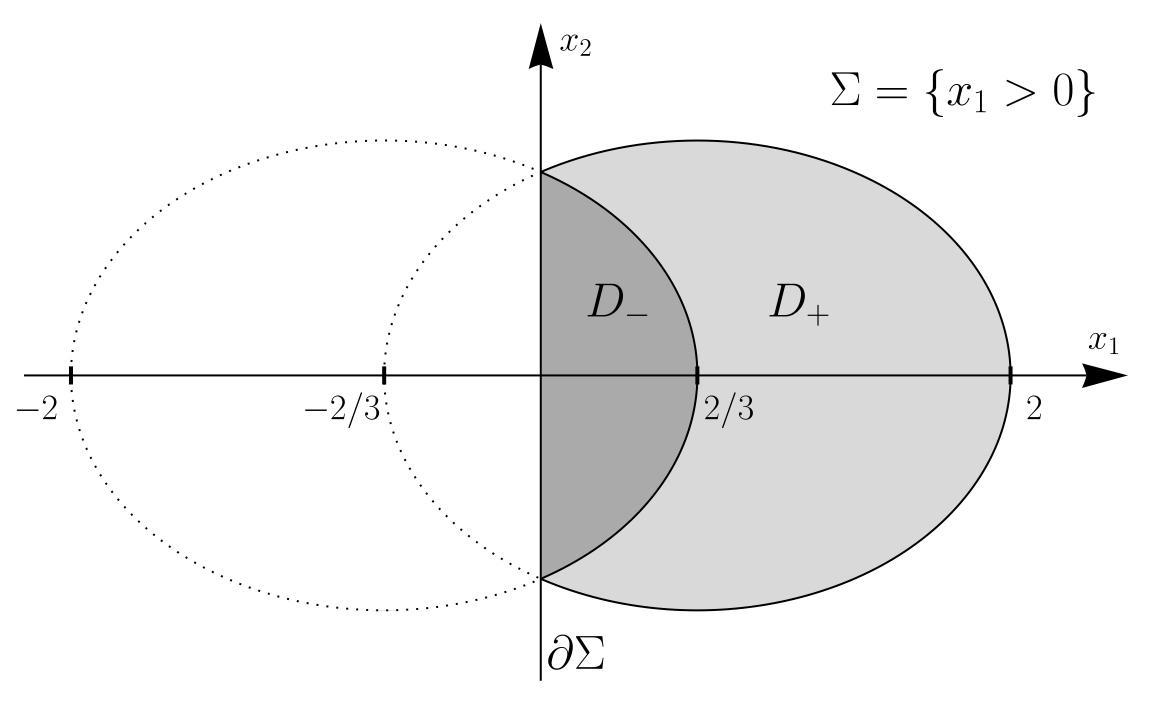}	\caption{The domains $D_\pm=\{x \in \Sigma:H(\pm x)<1\}$ for $H$ given by \eqref{H-norm}.}
	  	    \end{figure}
	    
%	    Note that $H$ is also a Randers metric, appearing in the Zermelo navigation problem, see Bao--Robles--Shen \cite{BRS}. 
	    It is clear that the synthetic dimension is $N=3$ and ${\sf Ric}_3\equiv  0$ on $\mathcal M=(  \Sigma,H,w\mathcal L^2_{\llcorner{\Sigma}})$. Moreover, the set of poles is the whole vertical line $\partial \Sigma=\{(0,x_2):x_2\in \mathbb R\}$.

	 If $D_\pm=\{x \in \Sigma:H(\pm x)<1\}$ as in Figure  1, an elementary computation  shows that 
	 $$\int_{D_+ }w=
	 %\int_{D^+}x_1{\rm d}x_1{\rm d}x_2=
	 \frac{4}{3}+\frac{32\pi}{27\sqrt{3}} \ \ {\rm and} \ \  \int_{D_-}w=
	   %\int_{D^-}x_1{\rm d}x_1{\rm d}x_2=
	   \frac{4}{3}-\frac{16\pi}{27\sqrt{3}}.$$
	   	   	   Therefore,  $${\sf AVR}_{\mathcal M}^+- {\sf AVR}_{\mathcal M}^-=\frac{16\pi}{9\sqrt{3}\omega_3}>0.$$

%	\subsection{Berwald spaces} If $\mathcal M=(M,F,{\sf m}_{BH})$ is an $n$-dimensional complete Berwald space without boundary having ${\sf Ric}_N\geq 0 $, it turns out that ${\sf S}_{{\sf m}_{BH}}\equiv 0$, see Shen \cite{Shen-konyv}, and  

%		\section{Questions}
%	
%		
%		1. Gagliardo-Nirenberg: ket ok miatt is: nem fog menni a skalazas, majd a sign-splitting modszer sem tunik mukodokepesnek (mert elveszti az "additivitasat" az egynlotlenseg).
		
%		2. $\Delta u = u^{2^*-1}$ megoldhato N-Ricci $\geq 0$ eseten, akkor mit tudunk mondani a sokasagrol?
	
	\vspace{0.5cm}
		
\noindent 	\textbf{Acknowledgment.} The author gratefully acknowledges the generosity of the organizers and the fruitful discussions held during the 2025 Geometry and Analysis Conference at the Department of Mathematics, Tsinghua University, Beijing, China, June 6-9, 2025, where this work was initiated.	
%He is  supported by the
%			Excellence Researcher Program \'OE-KP-2-2022 of \'Obuda University, Budapest, Hungary. 
		
		%%%%%%%%%%%%%%%%%%%%%%%%%%%%%%%%%%%%%%%%%%%%%%%%%%%%%%%%%%	

	\end{document}